\documentclass[reqno]{amsart}
\usepackage{amsmath}
\usepackage{mathtools}
\mathtoolsset{showonlyrefs}
\usepackage{amsthm}
\usepackage{amssymb}
\usepackage{float}
\usepackage{algorithm2e}
\usepackage{graphicx}
\usepackage[authoryear,sort]{natbib}

\makeatletter
\numberwithin{equation}{section}
\usepackage{paralist}

\usepackage{placeins}\usepackage{fancyhdr}\fancypagestyle{wmsappendix}
{\fancyhf{}\fancyhead[C]{\normalfont\bfseries Appendix}\fancyfoot[C]{\thepage}}
\AtBeginDocument{\setlength{\bibsep}{0pt}}

\numberwithin{equation}{section}%
\numberwithin{figure}{subsection}
\numberwithin{table}{subsection}
\usepackage{hyperref}
\usepackage{xcolor}
\usepackage{etoolbox}

\hypersetup{colorlinks=true,
linkcolor=blue,
anchorcolor=blue,
citecolor=black,
urlcolor=black}
\usepackage[letterpaper,margin=1in]{geometry}
\usepackage{setspace}
\usepackage{indentfirst}

\usepackage{fancyhdr}
\usepackage{caption}
\usepackage{booktabs}
\usepackage{siunitx}
\usepackage{iftex}
\ifXeTeX
  \usepackage{fontspec}
  \IfFontExistsTF{Times New Roman}
    {\setmainfont{Times New Roman}}
    {\setmainfont{TeX Gyre Termes}}
\else
  \usepackage[T1]{fontenc}
  \usepackage{times}
\fi
\usepackage{amsthm}
\newtheoremstyle{xxx}{}{}{\itshape}{}{\scshape\bfseries}{:}%
{\newline}{\thmname{#1}\thmnumber{ #2}\thmnote{(#3)}\smallskip}%
\theoremstyle{xxx}

\providecommand{\theoremname}{Theorem}
\ifdefined\AMFsubmissionlayout
\usepackage[nolists,tablesfirst,nomarkers]{endfloat}
\AtBeginFigures{}
\renewcommand{\l@figure}{\@tocline{0}{3pt plus2pt}{0pt}{2.2pc}{}}
\newcommand{\AMFpreservenumber}[1]{%
  \addtodelayedfloat*{#1}{%
    \string\setcounter{#1}{\number\value{#1}}%
    \string\def\expandafter\string\csname the#1\endcsname
      {\thesubsection.\string\arabic{#1}}%
    \string\def\expandafter\string\csname theH#1\endcsname
      {\thesubsection.\string\arabic{#1}}}%
  \stepcounter{#1}}
\AddToHook{env/figure/before}[amf-numbering]{\AMFpreservenumber{figure}}
\AddToHook{env/table/before}[amf-numbering]{\AMFpreservenumber{table}}
\AtBeginDelayedFloats{%
  \RemoveFromHook{env/figure/before}[amf-numbering]%
  \RemoveFromHook{env/table/before}[amf-numbering]}
\fi

\AtBeginDocument{\bibpunct{(}{)}{;}{a}{}{,}\hypersetup{urlcolor=black}\sloppy}

\makeatother

\providecommand\definitionname{Definition}
\theoremstyle{definition}
\newtheorem{defn}{\protect\definitionname}[section]
\providecommand\theoremname{Theorem}
\theoremstyle{plain}
\newtheorem{thm}{\protect\theoremname}[section]
\providecommand\propositionname{Proposition}
\newtheorem{prop}{\protect\propositionname}[section]
\providecommand\lemmaname{Lemma}
\newtheorem{lem}{\protect\lemmaname}[section]
\providecommand\corollaryname{Corollary}
\newtheorem{cor}{\protect\corollaryname}[section]
\providecommand\remarkname{Remark}
\theoremstyle{remark}
\newtheorem{rem}{\protect\remarkname}[section]
\usepackage[american]{babel}
\AtBeginDocument{\robustify\eqref}
\makeatletter
\renewcommand{\@setaddresses}{}
\apptocmd{\@setauthors}{\par\smallskip\begingroup
\noindent\begin{minipage}{\textwidth}\centering\normalfont\small
\textsuperscript{*}Department of Industrial Engineering and Economics,
Institute of Science Tokyo, 2-12-1 Ookayama, Meguro-ku, Tokyo 152-8552 Japan.\par
\textsuperscript{\dag}Department of Mathematics, Institute of Science Tokyo,
2-12-1 Ookayama, Meguro-ku, Tokyo 152-8551 Japan.\par
Corresponding author: Yuming Ma. Email: \href{mailto:ma.y.bb08@m.isct.ac.jp}{ma.y.bb08@m.isct.ac.jp}.\par
Syoiti Ninomiya: \href{mailto:syoiti.ninomiya@gmail.com}{syoiti.ninomiya@gmail.com}.\par
ORCIDs: Yuming Ma, \url{https://orcid.org/0009-0009-4407-3491}.
Syoiti Ninomiya, \url{https://orcid.org/0000-0003-4062-761X}.
\end{minipage}\endgroup\par}{}{\PackageError{authorinfo}{Could not add title-page author information}{}}
\newcommand{\printsubmissionmetadata}{\par\smallskip\begingroup\raggedright\noindent
\textbf{Keywords.} \submissionkeywords\par\endgroup}
\makeatother
\begin{document}
\title{A high-order recombination algorithm with adaptive scale selection
for weak approximation of stochastic differential equations}
\author[Y.~Ma]{Yuming Ma\textsuperscript{*}}
\address{
\textsuperscript{*}Department of Industrial Engineering and Economics,
Institute of Science Tokyo,
2-12-1 Ookayama, Meguro-ku, Tokyo 152-8552 JAPAN}
\email{ma.y.bb08@m.isct.ac.jp}
\urladdr[ORCID]{https://orcid.org/0009-0009-4407-3491}
\author[S.~Ninomiya]{Syoiti Ninomiya\textsuperscript{\dag}}      
\address{
\textsuperscript{\dag}Department of Mathematics,
Institute of Science Tokyo,
2-12-1 Ookayama, Meguro-ku, Tokyo 152-8551 JAPAN}
\email{syoiti.ninomiya@gmail.com}
\urladdr[ORCID]{https://orcid.org/0000-0003-4062-761X}

\begin{abstract}
  \noindent High-order recombination keeps weak approximation trees of diffusions small
  by a moment-preserving reduction of the measure on each patch of nodes.
  The weighted Litterer--Lyons (WLL) test of Ninomiya and Shinozaki
  sets the patch scale using a scalar \(\lambda>0\) chosen
  beforehand, whose appropriate value depends on coordinate units.
  First, the anisotropically weighted Litterer--Lyons (AWLL) rule
  selects the scale during recombination within a bounded factor of
  the measured increment amplitude. On a fixed set of nodes, expanding
  coordinate units leaves the AWLL patch-count bound unchanged but raises
  the WLL lower bound as a power until saturation. In Heston and shifted Stochastic-Alpha-Beta-Rho (sSABR)
  systems with initial price or forward rate 100, AWLL reduces
  the terminal supports of WLL with \(\lambda=1\) by factors of 137 and
  1488 at comparable errors, but enlarges them at smaller Heston initial
  prices. Second, \(\lambda\) is extended to a positive definite matrix
  with a joint error bound for prices and sensitivities, a bound on
  the gain of its shape, and patch-count bounds. Third, the elimination stage of batched basis
  exchange requires \(O(J)\)
  arithmetic operations for \(J\) nodes at fixed moment rank and panel width, against
  \(\Theta(J^{3})\) for the direct kernel update of Ninomiya and Shinozaki.
\end{abstract}


\makeatletter
\let\submissionkeywords\@keywords
\let\@keywords\@empty
\makeatother
\maketitle
\printsubmissionmetadata

\section{Introduction and background}

\label{sec:introduction}
\noindent
KLNV trees \citep{kusuoka_2013,morimoto_sasada_2017,ninomiya_shinozaki_2019}
compute expectations of
functionals of a diffusion, such as prices and parameter sensitivities
of derivatives, on finite trees whose branches follow the flows of
ordinary differential equations along the vector fields of the equation.
They are based on the approximation of \citet{kusuoka_2001}, which
S.~\citet{ninomiya_new_2003,ninomiya_partial_2003} showed to be effective
for problems in finance and whose theory \citet{kusuoka_2003,kusuoka_2004,kusuoka_2005,kusuoka_2013}
developed in parallel with cubature on Wiener space \citep{lyons_victoir_2004}.
The branches of a step are given by a high-order discretization method
based on discrete random variables, and with $d_{B}$ Brownian motions
each step multiplies the number of nodes by $2\times3^{d_{B}}$ for
the Ninomiya--Victoir (NV) method \citep{ninomiya_victoir_2008},
by $3^{2d_{B}}$ for the Ninomiya--Ninomiya method \citep{ninomiya_ninomiya_2009}
and, for $d_{B}=2$, by $5^{7}$ for the third-order method $Q^{(7,2)}_{(s)}$ of \citet{shinozaki_2017}.
The number of nodes therefore grows exponentially with the number
$n$ of time steps \citep[Table~1.1]{ninomiya_high-order_2025}, by
the factor 18 per step for the NV rules with $d_{B}=2$ used in this
paper. High-order recombination \citep{litterer_high_2012}
controls this growth without the sampling error of partial sampling
schemes such as the tree-based branching algorithm \citetext{\citealp{crisan_lyons_2002};
S.~\citealp{ninomiya_partial_2003}; M.~\citealp{ninomiya_2010}}.
At a reduction date the nodes are divided into
patches, and the probability measure on each patch is replaced by
a measure on at most as many of its nodes as the rank of its moment
matrix, with the same moments of degree at most $G$. Such a measure
exists by the finite form of Tchakaloff's theorem \citep[Theorem~2]{bayer_tchakaloff_2006},
and the local error of the reduction is of the order of the patch
radius raised to the power $G+1$.

Algorithms that simplify the support of a discrete measure were
applied to cubature on Wiener space by
\citet{litterer_lyons_2007} and \citet{litterer_2008},
and faster algorithms have since been developed \citep{tchernychova_caratheodory_2016,maalouf_jubran_feldman_2019,cosentino_oberhauser_abate_2020}.
Cubature formulas on Wiener space are constructed by reducing random
samples with linear programming \citep{hayakawa_tanaka_2022}, with
sample sizes bounded by estimates for random convex hulls \citep{hayakawa_lyons_oberhauser_2023,hayakawa_oberhauser_lyons_2023},
and formulas of higher degree or smaller support continue to be constructed
\citep{ferrucci_herschell_litterer_lyons_2026,litterer_2026,koepernik_coxon_foster_2026},
the last with a recombination step. Recombination also controls the
support in implementations of high-order weak approximation \citep{ninomiya_shinozaki_2021},
in a cubature filter \citep{lee_lyons_2016}, in Markov chain approximations
of diffusions on lattice trees \citep{cosentino_oberhauser_abate_2021}
and in the training of neural stochastic differential equations \citep{snow_krishnamurthy_2025}.
Outside these settings, it underlies convex kernel quadrature \citep{hayakawa_oberhauser_lyons_2022,hayakawa_oberhauser_lyons_2023b},
batch Bayesian quadrature and optimization \citep{adachi_2022,adachi_2023,adachi_2024},
stochastic gradient descent \citep{cosentino_oberhauser_abate_2020b}
and sparse approximation \citep{lyons_mcleod_2022,lyons_mcleod_2024}.

\citet{ninomiya_high-order_2025} allow each patch a radius that depends
on its mass, the step and the time to maturity, and divide the frontier
recursively until every patch passes this weighted Litterer--Lyons
(WLL) test. In their three-state Heston example with initial spot
one, the recombined support grows approximately as $n^{3.84}$ \citep[Section~4.2]{ninomiya_high-order_2025}.
The test contains a scaling parameter $\lambda>0$, which they vary
to control the degree of recombination \citep[Section~4.4]{ninomiya_high-order_2025}.
This scalar acts on all state directions and all dates alike and is
chosen in advance.

Its appropriate value depends on the units of the state coordinates,
since the WLL test measures every direction in the same Euclidean
coordinates, and the price, volatility, tangent and integral coordinates
of a model need not have comparable scales. In the six-state Heston
system of Section~\ref{sec:numerics} with initial spot 100, called
level 100 there, the number of terminal nodes of WLL with $\lambda=1$
grows nearly geometrically, by a factor of about 12.8 per added interval
up to $n=9$, and reaches $1.74\times10^{10}$ at nine intervals, where
the computation took about 477 calendar hours and 898 GiB of storage. A larger $\lambda$
relaxes all directions together. Choosing it by a sweep requires a
complete run for each value. The adaptive recombination of \citet[Eq.~(6.12) and Section~6.3.2]{lee_lyons_2016}
maps the particles affinely into $\left[0.5,1\right)^{N}$ and selects
the largest patch radius whose reduction error, measured with the
continuation function of the task, stays below a tolerance, which
requires continuation quantities that are costly when many prices
and sensitivities are computed together. Anisotropic approximation,
which adapts the shape of local cells to directional smoothness \citep{cohen_mirebeau_2011},
suggests replacing the scalar by a matrix.

This paper replaces the scalar $\lambda$ by a positive definite matrix
that can change at each reduction date and is selected from the tree
itself. Its contributions are threefold.

First, the anisotropically weighted Litterer--Lyons (AWLL) rule selects
the resolution scale of the WLL test during the recombination from
forward information, so that no value of $\lambda$ has to be chosen
in advance. Its matrix is computed in closed form from the covariance
of the current tree increments, the second moment of the increments
of the model, and a ridge (Section~\ref{sec:method}). We prove that
the matrix test then acts as a WLL test whose scalar is the measured
amplitude up to a bounded factor. On a fixed frontier, when the units
of a set of coordinates are expanded, a lower bound on the number of
WLL patches with a fixed scalar grows as a power of the expansion factor
until atomic saturation, whereas the AWLL upper bound does not depend
on that factor (Sections~\ref{sec:method} and \ref{sec:work}). In the six-state Heston and shifted Stochastic-Alpha-Beta-Rho (sSABR)
systems at level 100, AWLL reduces the terminal supports of WLL with
$\lambda=1$ by factors of 137 and 1488 at comparable reference errors,
and WLL runs at the scalars that AWLL selects confirm that this scale
accounts for the reduction (Sections~\ref{sec:numerics} and~\ref{sec:discussion}).

Second, the scalar $\lambda$ is extended to a positive definite matrix
with a joint error bound for prices and sensitivities in the metric
of the test, a bound on the gain of the shape of the matrix, and bounds
on the patch counts in terms of its resolution load (Sections~\ref{sec:framework}
and \ref{sec:work}). A synthetic system with explicit continuation
functions shows how the shape redistributes the error contributions
of tasks that vary along different axes (Section~\ref{sec:shape-test}).

Third, batched basis exchange computes each positive reduction with
$O(J)$ arithmetic operations after rank assessment and initialization,
where $J$ is the number of nodes of
the patch, for fixed moment rank and panel width, against $\Theta(J^{3})$
for the direct kernel update of \citet{ninomiya_high-order_2025} (Section~\ref{sec:solver}).
Its basis updates are the product-form updates of the revised simplex
method \citep{huangfu_hall_2018}, and the divide-and-conquer reduction
of \citet[Corollary~4]{litterer_high_2012} also uses work linear in
$J$.

Section~\ref{sec:framework} develops Taylor cancellation, metric localization
and the joint error theory. Section~\ref{sec:method} derives AWLL
and the positive reduction, and Section~\ref{sec:work} bounds the
patch counts and the retained populations. Section~\ref{sec:numerics}
compares support, accuracy and work, gives post hoc global WLL equivalents,
and tests the directional shape on a synthetic system. Section~\ref{sec:discussion}
shows that the resolution scale selected by AWLL accounts for the
observed reductions and discusses what the directional shape adds.
Section~\ref{sec:conclusion} states the contributions with their
limits and future directions. Appendix~\ref{app:geometry} gives the
unit response of the AWLL matrix, the amplitudes of the two models
and the action of the metric on the computed measures, and Appendix~\ref{app:algorithm}
gives the algorithm.

\FloatBarrier

\section{Theory}
\label{sec:framework}

\noindent
Throughout, inequalities between vectors hold componentwise. For a
matrix $X$, $\prescript{t}{}{X}$ denotes its transpose, $X_{i\cdot}$
and $X_{\cdot j}$ its $i$th row and $j$th column, and $X_{C\cdot}$,
$X_{\cdot C}$ and $X_{CC}$ its submatrices of the rows, of the columns,
and of the rows and columns with indices in $C$. Here $C$ is a set
of indices in increasing order or a tuple of indices in its own order,
and $x_{C}$ denotes the subvector of a vector $x$ with indices in
$C$. The identity matrix and the vector of ones are denoted by $I$
and $\mathbf{1}$, with their dimension as a subscript, and $e_{a}$
denotes the $a$th coordinate vector. The golden ratio $(1+\sqrt{5})/2$
is denoted by $\varphi$.

\subsection{Tree measures and derivative tasks}
\begin{defn}
\label{def:tree}Fix an open state domain $\mathcal{O}\subseteq\mathbb{R}^{N}$,
an initial state $Y_{0}\in\mathcal{O}$, a grid
$\Delta=\left\{ 0=t_{0}<\cdots<t_{n}=T\right\} $ with steps
$s_{i}=t_{i}-t_{i-1}$, and a set $\mathcal{R}\subseteq\left\{ 1,\ldots,n-1\right\} $
of reduction indices. For each $i\in\left\{ 1,\ldots,n\right\} $
let $\mathcal{B}_{i}$ be a finite nonempty set of branch labels, and
for $m\in\mathcal{B}_{i}$ let
$\pi_{i,m}:\mathcal{O}\to\left[0,1\right]$ and
$\Phi_{i,m}:\mathcal{O}\to\mathcal{O}$ satisfy
$\sum_{m\in\mathcal{B}_{i}}\pi_{i,m}\left(y\right)=1$
for every $y\in\mathcal{O}$. The superscript $\mathrm{Alg}$ of the operators below
labels the discretization method that produces these data (see the
paragraph after Remark~\ref{rem:lipschitz-semigroup}).
\begin{enumerate}
\item \emph{Branch space.} Put
$\widehat{\Omega}=\mathcal{B}_{1}\times\cdots\times\mathcal{B}_{n}$,
and let $\widehat{\mathcal{F}}_{i}$ be the $\sigma$-algebra generated
by the first $i$ coordinates, so that
$\widehat{\mathcal{F}}_{0}=\left\{ \emptyset,\widehat{\Omega}\right\} $
and $\widehat{\mathcal{F}}_{n}=2^{\widehat{\Omega}}$. The atoms
$\left[j\right]=\left\{ \omega\in\widehat{\Omega}:\left(\omega_{1},\ldots,\omega_{i}\right)=j\right\} $,
$j\in\mathcal{B}_{1}\times\cdots\times\mathcal{B}_{i}$, of
$\widehat{\mathcal{F}}_{i}$ are the nodes at date $i$. The state
maps $\widehat{Y}_{i}:\widehat{\Omega}\to\mathcal{O}$,
\begin{equation*}
\widehat{Y}_{0}\left(\omega\right)=Y_{0},\qquad \widehat{Y}_{i}\left(\omega\right)=\Phi_{i,\omega_{i}}\left(\widehat{Y}_{i-1}\left(\omega\right)\right),
\end{equation*}
are $\widehat{\mathcal{F}}_{i}$-measurable, and $\widehat{Y}_{i}\left(j\right)$
denotes the value of $\widehat{Y}_{i}$ on the node $\left[j\right]$. Distinct
nodes may have the same state.
\item \emph{Measures.} $\mathcal{M}_{i}$ is the space of finite
signed measures $\mu:\widehat{\mathcal{F}}_{i}\to\mathbb{R}$. The
support of $\mu$ is the set of nodes $\left[j\right]$ at date $i$
with $\mu\left(\left[j\right]\right)\ne0$. For an
$\widehat{\mathcal{F}}_{i}$-measurable $g:\widehat{\Omega}\to\mathbb{R}$
put $\mu\left(g\right)=\left\langle \mu,g\right\rangle =\int_{\widehat{\Omega}}g\,d\mu$,
and for $f:\mathcal{O}\to\mathbb{R}$ put
$\mu\left(f\right)=\left\langle \mu,f\right\rangle =\mu\left(f\circ \widehat{Y}_{i}\right)$.
\item \emph{Transition.} For $\omega\in\widehat{\Omega}$ and
$m\in\mathcal{B}_{i}$, let $\omega^{i\leftarrow m}$
be $\omega$ with its $i$th coordinate replaced by $m$. The
backward operator $Q^{\mathrm{Alg}}_{i-1\leftarrow i}$ maps
$\widehat{\mathcal{F}}_{i}$-measurable functions to
$\widehat{\mathcal{F}}_{i-1}$-measurable functions,
\begin{equation*}
\left(Q^{\mathrm{Alg}}_{i-1\leftarrow i}g\right)\left(\omega\right)=\sum_{m\in\mathcal{B}_{i}}\pi_{i,m}\left(\widehat{Y}_{i-1}\left(\omega\right)\right)g\left(\omega^{i\leftarrow m}\right),
\end{equation*}
and the forward operator
$K^{\mathrm{Alg}}_{i-1\to i}:\mathcal{M}_{i-1}\to\mathcal{M}_{i}$
is its adjoint,
$\left(K^{\mathrm{Alg}}_{i-1\to i}\mu\right)\left(F\right)=\mu\left(Q^{\mathrm{Alg}}_{i-1\leftarrow i}\mathbf{1}_{F}\right)$
for $F\in\widehat{\mathcal{F}}_{i}$, so that $\left\langle K^{\mathrm{Alg}}_{i-1\to i}\mu,g\right\rangle =\left\langle \mu,Q^{\mathrm{Alg}}_{i-1\leftarrow i}g\right\rangle $
for every $\widehat{\mathcal{F}}_{i}$-measurable $g$. It maps probability measures
to probability measures. The state operator acts on functions
$f:\mathcal{O}\to\mathbb{R}$ by
\begin{equation}
\left(Q^{\mathrm{Alg}}_{i}f\right)\left(y\right)=\sum_{m\in\mathcal{B}_{i}}\pi_{i,m}\left(y\right)f\left(\Phi_{i,m}\left(y\right)\right),\qquad\left(Q^{\mathrm{Alg}}_{i}f\right)\circ \widehat{Y}_{i-1}=Q^{\mathrm{Alg}}_{i-1\leftarrow i}\left(f\circ \widehat{Y}_{i}\right).\label{eq:ver1-007}
\end{equation}
\item \emph{Tree measures.} Let
$\mathbb{P}^{\mathrm{WT}}_{0}=\mathbb{P}^{R,+}_{0}$ be the unique
probability measure on $\widehat{\mathcal{F}}_{0}$, and define the whole
tree (WT) and recombined measures by
\begin{equation}
\mathbb{P}^{\mathrm{WT}}_{i}=K^{\mathrm{Alg}}_{i-1\to i}\mathbb{P}^{\mathrm{WT}}_{i-1},\qquad\mathbb{P}^{R,-}_{i}=K^{\mathrm{Alg}}_{i-1\to i}\mathbb{P}^{R,+}_{i-1},\qquad1\le i\le n.\label{eq:tree-recursion}
\end{equation}
For $i\notin\mathcal{R}$, $\mathbb{P}^{R,+}_{i}=\mathbb{P}^{R,-}_{i}$.
For $i\in\mathcal{R}$, $\mathbb{P}^{R,+}_{i}\in\mathcal{M}_{i}$
is a probability measure whose support is contained in that of
$\mathbb{P}^{R,-}_{i}$, produced by the positive reduction
\eqref{eq:moments}
below. Write
$p^{R,\pm}_{i,j}=\mathbb{P}^{R,\pm}_{i}\left(\left[j\right]\right)$,
and let $J^{\pm}_{i}$ be the number of nodes in the support of
$\mathbb{P}^{R,\pm}_{i}$. The measure $\mathbb{P}^{R,-}_{i}$ is called
the frontier at date $i$.
\end{enumerate}
\end{defn}
All measures of date $i$ live on the same measurable space
$\left(\widehat{\Omega},\widehat{\mathcal{F}}_{i}\right)$, so a reduction
changes the measure and not the space. The unrecombined measures are
consistent, $\mathbb{P}^{\mathrm{WT}}_{i}$ being the restriction
of $\mathbb{P}^{\mathrm{WT}}_{n}$ to $\widehat{\mathcal{F}}_{i}$.
The recombined ones are not: $\mathbb{P}^{R,-}_{i+1}$ restricts to
$\mathbb{P}^{R,+}_{i}$ on $\widehat{\mathcal{F}}_{i}$, which differs
from $\mathbb{P}^{R,-}_{i}$ at reduction dates. The image measures
$\widehat{Y}_{i\#}\mathbb{P}^{R,+}_{i}$ on $\mathcal{O}$ are the recombination
measures of \citet{ninomiya_high-order_2025}.

For a function $f:\mathcal{O}\to\mathbb{R}$, its discrete continuation
is
\begin{equation}
Q^{\mathrm{Alg}}_{(i,n]}f=Q^{\mathrm{Alg}}_{i+1}\cdots Q^{\mathrm{Alg}}_{n}f\quad\left(0\le i<n\right),\qquad Q^{\mathrm{Alg}}_{(n,n]}f=f.\label{eq:discrete-continuation}
\end{equation}
The signed recombination residual is
\begin{equation}
\mathbb{D}^{R}_{i}=\mathbb{P}^{R,+}_{i}-\mathbb{P}^{R,-}_{i}\in\mathcal{M}_{i},\qquad\mathbb{D}^{R}_{i}=0\quad\text{for }i\notin\mathcal{R}.\label{eq:residual}
\end{equation}

Let $B=\left(B^{1},\ldots,B^{d_{B}}\right)$ be a standard Brownian
motion on a filtered probability space satisfying the usual conditions,
with expectation $\mathbb{E}$, and put $B^{0}_{t}=t$. Let
$V_{0},\ldots,V_{d_{B}}$ be smooth tangent vector fields on
$\mathcal{O}$. For a tangent vector field $V$,
$VI_{N}:\mathcal{O}\to\mathbb{R}^{N}$ denotes its coefficient function,
so that $Vg\left(y\right)=Dg\left(y\right)\left[VI_{N}\left(y\right)\right]$
for $g\in C^{1}\left(\mathcal{O}\right)$, where $I_{N}$ is the identity
map of $\mathbb{R}^{N}$. We consider the stochastic differential equation
in Stratonovich form
\begin{equation}
Y_{t}=Y_{0}+\sum^{d_{B}}_{a=0}\int^{t}_{0}V_{a}I_{N}\left(Y_{v}\right)\circ dB^{a}_{v}.\label{eq:augmented}
\end{equation}
Its It\^o form $dY_{t}=b\left(Y_{t}\right)dt+\Sigma\left(Y_{t}\right)dB_{t}$
has
\begin{equation*}
b=V_{0}I_{N}+\frac{1}{2}\sum^{d_{B}}_{a=1}V_{a}\left(V_{a}I_{N}\right),\qquad\Sigma=\left(V_{1}I_{N},\ldots,V_{d_{B}}I_{N}\right),\qquad A=\Sigma\prescript{t}{}{\Sigma},
\end{equation*}
where $V_{a}$ acts on each component of $V_{a}I_{N}$. Assume that
for every $y\in\mathcal{O}$ the equation~\eqref{eq:augmented} with
$Y_{0}$ replaced by $y$ has a unique strong solution
$\left(Y^{y}_{t}\right)_{t\in\left[0,T\right]}$ that stays in $\mathcal{O}$.
As in \citet[Definition~2.1]{ninomiya_high-order_2025}, put
\begin{equation}
\left(P_{t}f\right)\left(y\right)=\mathbb{E}\left[f\left(Y^{y}_{t}\right)\right],\qquad y\in\mathcal{O},\ t\in\left[0,T\right],\label{eq:semigroup}
\end{equation}
for bounded Borel functions $f:\mathcal{O}\to\mathbb{R}$. By the
Markov property of the solution, $\left(P_{t}\right)_{t\in\left[0,T\right]}$
is a family of linear operators on the bounded Borel functions on
$\mathcal{O}$ with $P_{0}f=f$ and $P_{t}P_{u}=P_{t+u}$ for $t+u\le T$.
\begin{rem}
\label{rem:lipschitz-semigroup}\citet[Definition~2.1]{ninomiya_high-order_2025}
define the semigroup on Lipschitz functions, for vector fields on
$\mathbb{R}^{N}$ with coefficients in $C^{\infty}_{b}$. Here the
coefficients need not be bounded and the state domain is $\mathcal{O}$.
Formula \eqref{eq:semigroup} still defines $P_{t}f\left(y\right)$ for every
Borel $f$ with $\mathbb{E}\left[\left|f\left(Y^{y}_{t}\right)\right|\right]<\infty$,
the condition placed on task functions in Definition~\ref{def:tasks}.
If $\mathbb{E}\left[\left\lVert Y^{y}_{t}\right\rVert _{2}\right]<\infty$, this
includes every Lipschitz $f$, such as the call and put payoffs, since
$\left|f\left(Y^{y}_{t}\right)\right|\le\left|f\left(y\right)\right|+\left\lVert f\right\rVert _{\mathrm{Lip}}\left\lVert Y^{y}_{t}-y\right\rVert _{2}$.
\end{rem}
In this paper the data $\mathcal{B}_{i}$, $\pi_{i,m}$ and $\Phi_{i,m}$
of Definition~\ref{def:tree} are produced from $V_{0},\ldots,V_{d_{B}}$
and $s_{i}$ by a discretization method of \eqref{eq:augmented} on
$\Delta$ based on discrete random variables. As in
\citet[Section~1.2.2]{ninomiya_high-order_2025}, the label $\mathrm{Alg}$
names this method, for example $\mathrm{Alg}=\mathrm{EM}$, $\mathrm{NV}$
or $\mathrm{NN}$ for the Euler--Maruyama, Ninomiya--Victoir
\citep{ninomiya_victoir_2008} and Ninomiya--Ninomiya \citep{ninomiya_ninomiya_2009}
methods. The method therefore determines the operators $Q^{\mathrm{Alg}}_{i}$,
$Q^{\mathrm{Alg}}_{i-1\leftarrow i}$ and $K^{\mathrm{Alg}}_{i-1\to i}$, and the
branch numbers $b^{\mathrm{Alg}}_{i}$ of Section~\ref{sec:work}.
In the notation of \citet{ninomiya_high-order_2025},
$\widehat{Y}_{i\#}K^{\mathrm{Alg}}_{i-1\to i}\mu=\widehat{Q}^{(\mathrm{Alg},\Delta)}_{(s_{i})}\left(\widehat{Y}_{i-1\#}\mu\right)$
for $\mu\in\mathcal{M}_{i-1}$, where $\widehat{Q}^{(\mathrm{Alg},\Delta)}_{(s_{i})}$
is their one-step forward operator on discrete measures.
The state operator $Q^{\mathrm{Alg}}_{i}$ is the discrete counterpart
of $P_{s_{i}}$,
and the discrete continuation~\eqref{eq:discrete-continuation}
is that of $P_{T-t_{i}}f=P_{s_{i+1}}\cdots P_{s_{n}}f$.
Section~\ref{sec:numerics} uses $\mathrm{Alg}=\mathrm{NV}$.
\begin{defn}
\label{def:tasks}A task family consists of a finite label set $\mathcal{I}$
and, for each $\alpha\in\mathcal{I}$, a Borel function $f_{\alpha}:\mathcal{O}\to\mathbb{R}$
and a tolerance $\varepsilon_{\alpha}>0$, such that
$\mathbb{E}\left[\left|f_{\alpha}\left(Y^{y}_{t}\right)\right|\right]<\infty$
for every $y\in\mathcal{O}$ and $t\in\left[0,T\right]$. Its exact
continuations and targets are
\begin{equation}
h^{\alpha}_{i}=P_{T-t_{i}}f_{\alpha}\quad\left(0\le i\le n\right),\qquad h^{\alpha}_{0}\left(Y_{0}\right)=P_{T}f_{\alpha}\left(Y_{0}\right).\label{eq:tasks}
\end{equation}
\end{defn}
\begin{rem}
\label{rem:greeks}The targets include prices and parameter sensitivities.
Let the coefficients of a base equation in $\mathbb{R}^{N_{0}}$ depend
smoothly on a real parameter $\theta$, and augment its solution $X$
by the tangent process $\partial_{\theta}X$, which solves the linearized
equation, and by the running integrals a payoff needs. For a payoff
$g$ of the base state that is differentiable outside a set which $X_{T}$
avoids almost surely, and under an integrable dominating bound for
the difference quotients,
$\partial_{\theta}\mathbb{E}\left[g\left(X_{T}\right)\right]=\mathbb{E}\left[Dg\left(X_{T}\right)\partial_{\theta}X_{T}\right]$.
This is the target of the task $f\left(x,\xi\right)=Dg\left(x\right)\xi$
on the augmented state. Section~\ref{sec:numerics} uses prices, first
model vegas of this form and signed exercise probabilities, from which
the Deltas in \eqref{eq:numerical-delta} are formed.
\end{rem}

At $i\in\mathcal{R}$, let $U_{i,\kappa}$, $1\le\kappa\le K_{i}$,
be nonempty disjoint sets of nodes of date $i$ whose union is the
support of $\mathbb{P}^{R,-}_{i}$, let $\mathbb{P}^{R,-}_{i,\kappa}$
be the restriction of $\mathbb{P}^{R,-}_{i}$ to the nodes in $U_{i,\kappa}$,
and put $w_{i,\kappa}=\mathbb{P}^{R,-}_{i,\kappa}\left(1\right)>0$.
For an integer $G\ge1$, the \emph{recombination} at date $i$, also
called positive reduction below, replaces each
$\mathbb{P}^{R,-}_{i,\kappa}$ by a measure $\mathbb{P}^{R,+}_{i,\kappa}\in\mathcal{M}_{i}$
such that
\begin{equation}
\begin{aligned}\operatorname{supp}\left(\mathbb{P}^{R,+}_{i,\kappa}\right) & \subseteq U_{i,\kappa},\qquad\mathbb{P}^{R,+}_{i,\kappa}\ge0,\\
\left\langle \mathbb{D}^{R}_{i,\kappa},g\right\rangle  & =0,\qquad g\in\mathcal{P}_{G}\left(\mathbb{R}^{N}\right),
\end{aligned}
\label{eq:moments}
\end{equation}
where $\mathbb{D}^{R}_{i,\kappa}=\mathbb{P}^{R,+}_{i,\kappa}-\mathbb{P}^{R,-}_{i,\kappa}$
and $\mathcal{P}_{G}\left(\mathbb{R}^{N}\right)$ is the space of
real polynomials of total degree at most $G$. Put
$\mathbb{P}^{R,+}_{i}=\sum^{K_{i}}_{\kappa=1}\mathbb{P}^{R,+}_{i,\kappa}$.
Taking $g=1$ gives $\mathbb{P}^{R,+}_{i,\kappa}\left(1\right)=w_{i,\kappa}$,
so $\mathbb{P}^{R,+}_{i}$ is a probability measure,
and \eqref{eq:residual}
gives $\mathbb{D}^{R}_{i}=\sum^{K_{i}}_{\kappa=1}\mathbb{D}^{R}_{i,\kappa}$.
For a set $U$ of nodes of date $i$, radii, convex hulls and maxima
over $U$ below refer to its finite state set $\widehat{Y}_{i}\left(U\right)\subset\mathcal{O}$.

Condition \eqref{eq:moments} does not bound the number of retained
nodes, and $\mathbb{P}^{R,+}_{i,\kappa}=\mathbb{P}^{R,-}_{i,\kappa}$
satisfies it. In exact arithmetic the batched basis exchange of Section~\ref{sec:solver}
returns a solution supported on linearly independent moment columns,
that is, an extreme point of the feasible set, with at most
$M_{G}=\binom{N+G}{G}=\dim\mathcal{P}_{G}\left(\mathbb{R}^{N}\right)$
nodes (Proposition~\ref{prop:solver} and Theorem~\ref{thm:population}).
This is the finite form of Tchakaloff's theorem
\citep[Theorem~2]{bayer_tchakaloff_2006}. A patch whose moment columns
are already linearly independent has its input weights as its only
feasible point, and the implementation carries it forward unchanged.

Thus, on every patch, the image of $\mathbb{P}^{R,+}_{i,\kappa}/w_{i,\kappa}$
on $\mathcal{O}$ is a reduced measure of the image of
$\mathbb{P}^{R,-}_{i,\kappa}/w_{i,\kappa}$ with respect to the monomials
of degree at most $G$, in the sense of
\citet[Definition~1.2]{ninomiya_high-order_2025}. The map
$\mathbb{P}^{R,-}_{i}\mapsto\mathbb{P}^{R,+}_{i}$ is therefore the
measure reduction, or high-order recombination, of
\citet[Definitions~1 and~17]{litterer_high_2012}
and
\citet[Eq.~(1.4)]{ninomiya_high-order_2025}.
A failed numerical
acceptance test returns $\mathbb{P}^{R,-}_{i,\kappa}$ itself
(Algorithm~\ref{alg:bbe}), which keeps the moment equations but may
violate the cardinality condition of that definition.

The measure in~\eqref{eq:moments}
retains its patch mass rather than being divided
by it.
Both $\mathbb{P}^{R,\pm}_{i,\kappa}$ are nonnegative with mass $w_{i,\kappa}$,
so the total variation (TV) norm $\left\lVert \cdot\right\rVert _{\mathrm{TV}}$,
the sum of absolute atomic weights, gives $\left\lVert \mathbb{D}^{R}_{i,\kappa}\right\rVert _{\mathrm{TV}}\le2w_{i,\kappa}$,
and the same bound holds for the image measure $\widehat{Y}_{i\#}\mathbb{D}^{R}_{i,\kappa}$
on $\mathcal{O}$, since taking images does not increase the TV norm.

\subsection{From a scalar radius to directional resolution}\label{sec:directional}
\noindent
For a set $U$ of nodes, we denote by
$\operatorname{rad}_{I_{N}}\left(U\right)=\max_{x,y\in U}\left\lVert x-y\right\rVert _{2}$
its Euclidean pairwise radius. In the present notation, the WLL condition
of \citet[Eq.~(3.7)]{ninomiya_high-order_2025} requires every patch
at date $i$ to satisfy
\begin{equation*}
\operatorname{rad}_{I_{N}}\left(U_{i,\kappa}\right)\le\left(T-t_{i}\right)^{G/\left(2\left(G+1\right)\right)}s^{1/2}_{i}\left(K_{i}w_{i,\kappa}\right)^{-1/\left(G+1\right)}.
\end{equation*}
Under this condition, the date $i$ recombination term of
\citet[Theorem~3.1]{ninomiya_high-order_2025} is at most a constant
multiple of $s^{\left(G+1\right)/2}_{i}$, whatever $K_{i}$ is. Since
$K_{i}$ is not known until the subdivision is completed, the recursive
algorithm of \citet[Section~3.2]{ninomiya_high-order_2025} replaces
$K_{i}$ by one, and \citet[Eq.~(4.2)]{ninomiya_high-order_2025} multiplies
the right side by a scaling parameter $\lambda>0$, which is common
to all coordinates and reduction dates. We thus obtain the radius test
\begin{equation}
\operatorname{rad}_{I_{N}}\left(U_{i,\kappa}\right)\le\lambda\left(T-t_{i}\right)^{G/\left(2\left(G+1\right)\right)}s^{1/2}_{i}w^{-1/\left(G+1\right)}_{i,\kappa}\label{eq:raw-rule}
\end{equation}
and refer to it as the WLL test. Since the factor $K_{i}^{-1/(G+1)}$
is absent from \eqref{eq:raw-rule},
the patch count appears in \eqref{eq:joint-wms-budget}
and in the order condition \eqref{eq:matched-order-condition}.
As in \citet[Section~3.2 and Appendix~C]{ninomiya_high-order_2025}, the
subdivision is also stopped at a candidate with at most $M_{G}$ nodes,
whether or not the candidate passes the radius test. This rule is
referred to as the population rule, and it accompanies every radius
test below. A stopped patch is reduced like any other patch, and it is
carried forward unchanged when its moment columns are linearly independent
(see the remarks after \eqref{eq:moments}).

Suppose that the WLL test is applied in the computational coordinates
$y'=Dy$, obtained from the coordinates $y$ by a positive diagonal
change of units $D$.
Then the largest displacement in $y$ along the $j$th coordinate axis
that \eqref{eq:raw-rule} admits, which we call the allowance in $y$
along that axis, is the right side of \eqref{eq:raw-rule} divided
by $D_{jj}$. For the Heston unit map \eqref{eq:ver1-149}, for example,
$D_{jj}$ is the unit factor $c$ for the spot, its tangent and the
two running integrals, and $D_{jj}=1$ for the variance and its tangent.
A larger $\lambda$ multiplies all allowances in $y$ by the same factor,
so their ratios $D_{kk}/D_{jj}$ are fixed by the choice of units and
cannot be changed through $\lambda$.

Let $\Gamma_{i}$ be a symmetric positive definite matrix, and put $M_{i}=\Gamma^{-1/2}_{i}$.
We define
\begin{equation}
  \operatorname{rad}_{{M}_{i}}\left(U\right)=\max_{x,y\in U}\left\lVert {M}_{i}\left(x-y\right)\right\rVert _{2}\label{eq:rad}
\end{equation}
and call the following inequality the matrix test:
\begin{equation}
  \operatorname{rad}_{{M}_{i}}\left(U_{i,\kappa}\right)\le\lambda\left(T-t_{i}\right)^{G/\left(2\left(G+1\right)\right)}s^{1/2}_{i}w^{-1/\left(G+1\right)}_{i,\kappa}.\label{eq:partition-rule}
\end{equation}
It reduces to the WLL test when $\Gamma_{i}=I_{N}$,
and the positive reduction \eqref{eq:moments} is unchanged. The recursive
subdivision also computes its splitting directions in the coordinates
$M_{i}y$ (Section~\ref{sec:subdivision}). Suppose that $\Gamma_{i}$
transforms by congruence under an invertible linear state map $y'=Wy$,
that is, $\Gamma'_{i}=W\Gamma_{i}\prescript{t}{}{W}$. Then $M'_{i}=\left(\Gamma'_{i}\right)^{-1/2}$
satisfies $M'_{i}W=QM_{i}$, where $Q=\left(W\Gamma_{i}\prescript{t}{}{W}\right)^{-1/2}W\Gamma^{1/2}_{i}$
is orthogonal. Hence $\left\lVert M'_{i}Wz\right\rVert _{2}=\left\lVert M_{i}z\right\rVert _{2}$
for every $z\in\mathbb{R}^{N}$, and the matrix test accepts the same
patches in both coordinate systems. Since the splitting coordinates
differ only by $Q$, the subdivision is also the same up to the sign
and tie conventions of its directions. The WLL test, on the other
hand, uses $I_{N}$ in both systems, and $W\prescript{t}{}{W}=I_{N}$
holds only for orthogonal $W$. On a fixed frontier, the matrix selected
in Section~\ref{sec:method} transforms by congruence when $W$ is
a positive multiple of $I_{N}$ and no mixed block is zero (Theorem~\ref{thm:spectral-scale}),
but in general not under the diagonal
unit maps of Section~\ref{sec:numerics}. Appendix~\ref{sec:transport}
gives its response to such maps.

\begin{lem}
  \label{lem:metric-scalar}Fix a reduction date, and let $\Gamma\succ0$ and $\lambda>0$.
  \begin{enumerate}
  \item The pairs $\left(\lambda,\Gamma\right)$ and $\left(1,\lambda^{2}\Gamma\right)$
    have identical acceptance sets under \eqref{eq:partition-rule}.
  \item For every $z\in\mathbb{R}^{N}$,
    \begin{equation}
      \lambda_{\max}\left(\Gamma\right)^{-1/2}\left\lVert z\right\rVert _{2}\le\left\lVert \Gamma^{-1/2}z\right\rVert _{2}\le\lambda_{\min}\left(\Gamma\right)^{-1/2}\left\lVert z\right\rVert _{2},\label{eq:metric-sandwich}
    \end{equation}
    with equality along the corresponding extremal eigenvectors. Hence
    every patch accepted by the WLL test with scalar $\lambda\lambda_{\min}\left(\Gamma\right)^{1/2}$
    is accepted by the matrix test, and every patch accepted by the matrix
    test is accepted by the WLL test with scalar $\lambda\lambda_{\max}\left(\Gamma\right)^{1/2}$.
    In particular, for a common scalar, the matrix test accepts every
    candidate patch in $\mathbb{R}^{N}$ that the WLL test accepts if
    and only if $\Gamma\succeq I_{N}$, and it accepts only such patches
    if and only if $\Gamma\preceq I_{N}$.
  \end{enumerate}
\end{lem}
\begin{proof}
  Item~1 follows from $\left(\lambda^{2}\Gamma\right)^{-1/2}=\lambda^{-1}\Gamma^{-1/2}$.
  Inequality \eqref{eq:metric-sandwich} is the Rayleigh quotient bound
  for $\Gamma^{-1}$ \citep[Appendix~A.5.2]{boyd_vandenberghe_2004},
  and the inclusions follow pair by pair. Two-point patches show that
  the last two inclusions hold for every patch only if $\Gamma^{-1}\preceq I_{N}$
  and only if $\Gamma^{-1}\succeq I_{N}$, respectively.
\end{proof}
A patch passes the matrix test if and only if all its pairwise state
displacements lie in the ellipsoid whose principal axes are
the eigenvectors of $\Gamma$ and
whose semi-axes are the right side of \eqref{eq:partition-rule}
times the square roots of the corresponding eigenvalues. By item~1
of Lemma~\ref{lem:metric-scalar}, multiplying $\Gamma$ by a positive
number has the same effect as multiplying $\lambda$ by its square
root. Hereafter, we therefore regard the multiplication of $\Gamma$
by a positive number as a change of scale, and $\Gamma/\lambda_{\max}\left(\Gamma\right)$
as the shape of $\Gamma$. Item~2 of Lemma~\ref{lem:metric-scalar}
compares the two tests at one date, whereas the post hoc WLL equivalent
of Section~\ref{sec:global-lambda} is a single scalar fitted to the
patch counts, retained populations and patch sizes of a complete AWLL
history.
\subsection{Joint error and the first uncanceled tensor}
\noindent
Let $\Gamma_{i}$ and $M_{i}=\Gamma^{-1/2}_{i}$ be as in Section~\ref{sec:directional}.
For each patch $U_{i,\kappa}$, we fix a point $c_{i,\kappa}$ in the convex
hull $\mathrm{C.H.}(U_{i,\kappa})$ of $U_{i,\kappa}$ and call it the center
of the patch. Put
\begin{equation}
\begin{aligned}r^{{M}}_{i,\kappa} & =\max_{y\in U_{i,\kappa}}\left\lVert {M}_{i}\left(y-c_{i,\kappa}\right)\right\rVert _{2},\\
\mathcal{Q}^{{M}}_{i,\kappa,G+1} & =\begin{cases}
\dfrac{{\displaystyle \int\left({M}_{i}\left(y-c_{i,\kappa}\right)\right)^{\otimes\left(G+1\right)}\left(\widehat{Y}_{i\#}\mathbb{D}^{R}_{i,\kappa}\right)\left(dy\right)}}{w_{i,\kappa}\left(r^{{M}}_{i,\kappa}\right)^{G+1}}, & \text{if }r^{{M}}_{i,\kappa}>0,\\
0, & \text{if }r^{{M}}_{i,\kappa}=0.
\end{cases}
\end{aligned}
\label{eq:residual-tensor}
\end{equation}
The quantity $r^{{M}}_{i,\kappa}$ is called the centered metric radius
of the patch. Since \eqref{eq:moments}
cancels the moments of $\widehat{Y}_{i\#}\mathbb{D}^{R}_{i,\kappa}$
of degrees zero through $G$, the tensor $\mathcal{Q}^{{M}}_{i,\kappa,G+1}$
is referred to as the first uncanceled tensor. It does not change
when $M_{i}$ is multiplied by a positive number. The full componentwise
contraction of two tensors of equal order on $\mathbb{R}^{N}$ is denoted
by $\left\langle \cdot,\cdot\right\rangle _{\otimes}$. For a function
$h$ on $\mathcal{O}$, the $(G+1)$-linear form obtained by applying
${M}^{-1}_{i}$ to each argument of $D^{G+1}h\left(c_{i,\kappa}\right)$
is denoted by $D^{G+1}h\left(c_{i,\kappa}\right)\circ\left({M}^{-1}_{i}\right)^{\otimes\left(G+1\right)}$.

For $i\in\mathcal{R}$, the convex hull $\mathrm{C.H.}(\operatorname{supp}(\mathbb{P}^{R,-}_{i}))$
contains $\mathrm{C.H.}(U_{i,\kappa})$ for every patch
of every partition of the support, whichever matrix the test uses.
If the derivatives of each $h^{\alpha}_{i}$ up to order $G+1$ are
continuous on an open subset of $\mathcal{O}$ that contains this hull,
we put, for $\Gamma\succ0$,
\begin{equation}
L^{\left(G+1\right)}_{i}\left(\Gamma\right)=\max_{\alpha\in\mathcal{I}}\frac{1}{\varepsilon_{\alpha}}\sup_{\substack{x\in\mathrm{C.H.}(\operatorname{supp}(\mathbb{P}^{R,-}_{i}))\\
\left\lVert v\right\rVert _{2}=1
}
}\left|D^{G+1}h^{\alpha}_{i}\left(x\right)\left[\Gamma^{1/2}v,\ldots,\Gamma^{1/2}v\right]\right|.
\label{eq:derivative-envelope}
\end{equation}
The map $\Gamma\mapsto L^{\left(G+1\right)}_{i}\left(\Gamma\right)$
is homogeneous of degree $(G+1)/2$ and monotone in the Loewner order,
since $\Gamma^{1/2}$ maps the unit ball onto the ellipsoid $\left\{ z:\prescript{t}{}{z}\Gamma^{-1}z\le1\right\} $,
these ellipsoids are nested when the matrices are ordered, and the
supremum of a form of degree $G+1$ over a centered ellipsoid is attained
on its boundary. Consequently,
\begin{equation}
\lambda_{\min}\left(\Gamma\right)^{\left(G+1\right)/2}L^{\left(G+1\right)}_{i}\left(I_{N}\right)\le L^{\left(G+1\right)}_{i}\left(\Gamma\right)\le\lambda_{\max}\left(\Gamma\right)^{\left(G+1\right)/2}L^{\left(G+1\right)}_{i}\left(I_{N}\right).\label{eq:envelope-sandwich}
\end{equation}

The following theorem extends the weighted recombination estimate
of \citet[Theorem~3.1]{ninomiya_high-order_2025}, which refines
\citet[Theorem~19]{litterer_high_2012} by the patch mass factor retained
here. In comparison with that estimate, the theorem replaces the Euclidean
pairwise radius by the centered metric radius, which is at most $\operatorname{rad}_{M_{i}}\left(U_{i,\kappa}\right)$,
separates the matching degree $G$ from the weak order of the scheme,
and controls a finite family of tasks jointly through their individual
tolerances.
\begin{thm}
\label{thm:residual}Let the matrices $\Gamma_{i}$ and the centers $c_{i,\kappa}$
be as above.
\begin{enumerate}
\item For every $f:\mathcal{O}\to\mathbb{R}$,
\begin{equation}
\mathbb{P}^{R,-}_{n}\left(f\circ \widehat{Y}_{n}\right)-\mathbb{P}^{\mathrm{WT}}_{n}\left(f\circ \widehat{Y}_{n}\right)=\sum_{i\in\mathcal{R}}\left\langle \mathbb{D}^{R}_{i},\left(Q^{\mathrm{Alg}}_{(i,n]}f\right)\circ \widehat{Y}_{i}\right\rangle .\label{eq:discrete-residual}
\end{equation}
\item Let $i\in\mathcal{R}$, and let $h:\mathcal{O}\to\mathbb{R}$
have continuous derivatives up to order $G+2$ on an open subset of
$\mathcal{O}$ that contains $\mathrm{C.H.}(\operatorname{supp}(\mathbb{P}^{R,-}_{i}))$.
Then, for every patch $U_{i,\kappa}$,
\begin{equation}
\begin{aligned}\left\langle \mathbb{D}^{R}_{i,\kappa},h\right\rangle ={} & \frac{w_{i,\kappa}\left(r^{{M}}_{i,\kappa}\right)^{G+1}}{\left(G+1\right)!}\left\langle D^{G+1}h\left(c_{i,\kappa}\right)\circ\left({M}^{-1}_{i}\right)^{\otimes\left(G+1\right)},\mathcal{Q}^{{M}}_{i,\kappa,G+1}\right\rangle _{\otimes}\\
 & +O\left(w_{i,\kappa}\left(r^{{M}}_{i,\kappa}\right)^{G+2}\right),
\end{aligned}
\label{eq:first-tensor}
\end{equation}
where the absolute value of the remainder term is at most
\begin{equation}
\frac{2w_{i,\kappa}\left(r^{{M}}_{i,\kappa}\right)^{G+2}}{\left(G+2\right)!}\sup_{\substack{x\in\mathrm{C.H.}(\operatorname{supp}(\mathbb{P}^{R,-}_{i}))\\
\left\lVert v\right\rVert _{2}=1
}
}\left|D^{G+2}h\left(x\right)\left[{M}^{-1}_{i}v,\ldots,{M}^{-1}_{i}v\right]\right|.\label{eq:ver1-032}
\end{equation}
\item Let a task family be as in Definition~\ref{def:tasks}, and assume
that, for each $i\in\mathcal{R}$, the derivatives of each $h^{\alpha}_{i}$
up to order $G+1$ are continuous on an open subset of $\mathcal{O}$
that contains $\mathrm{C.H.}(\operatorname{supp}(\mathbb{P}^{R,-}_{i}))$.
Then the normalized joint error satisfies
\begin{equation}
\begin{aligned}\max_{\alpha\in\mathcal{I}}\frac{\left|\mathbb{P}^{R,-}_{n}\left(f_{\alpha}\circ \widehat{Y}_{n}\right)-h^{\alpha}_{0}\left(Y_{0}\right)\right|}{\varepsilon_{\alpha}}\le{} & \sum^{n}_{i=1}\max_{\alpha\in\mathcal{I}}\frac{1}{\varepsilon_{\alpha}}\left|\left\langle \mathbb{P}^{R,+}_{i-1},Q^{\mathrm{Alg}}_{i-1\leftarrow i}\left(h^{\alpha}_{i}\circ \widehat{Y}_{i}\right)-h^{\alpha}_{i-1}\circ \widehat{Y}_{i-1}\right\rangle \right|\\
 & +\frac{2}{\left(G+1\right)!}\sum_{i\in\mathcal{R}}L^{\left(G+1\right)}_{i}\left(\Gamma_{i}\right)\sum^{K_{i}}_{\kappa=1}w_{i,\kappa}\left(r^{{M}}_{i,\kappa}\right)^{G+1}.
\end{aligned}
\label{eq:joint-error}
\end{equation}
Moreover, \eqref{eq:joint-error} remains valid when the inner sum
is restricted to the patches $U_{i,\kappa}$ with $\widehat{Y}_{i\#}\mathbb{D}^{R}_{i,\kappa}\ne0$.
\end{enumerate}
\end{thm}
\begin{proof}
Let $v_{0},\ldots,v_{n}$ be arbitrary functions on $\mathcal{O}$.
The first identity below follows from \eqref{eq:tree-recursion} and
the duality in Definition~\ref{def:tree}~(3), and the second is the
definition~\eqref{eq:residual} of the residual:
\begin{equation}
\begin{aligned} & \mathbb{P}^{R,-}_{i}\left(v_{i}\circ \widehat{Y}_{i}\right)-\mathbb{P}^{R,+}_{i-1}\left(v_{i-1}\circ \widehat{Y}_{i-1}\right)\\
 & \quad=\left\langle \mathbb{P}^{R,+}_{i-1},Q^{\mathrm{Alg}}_{i-1\leftarrow i}\left(v_{i}\circ \widehat{Y}_{i}\right)-v_{i-1}\circ \widehat{Y}_{i-1}\right\rangle ,\\
 & \mathbb{P}^{R,+}_{i}\left(v_{i}\circ \widehat{Y}_{i}\right)-\mathbb{P}^{R,-}_{i}\left(v_{i}\circ \widehat{Y}_{i}\right)=\left\langle \mathbb{D}^{R}_{i},v_{i}\right\rangle .
\end{aligned}
\label{eq:ver1-034}
\end{equation}
Adding these equalities for $1\le i\le n$, and using $\mathbb{P}^{R,+}_{n}=\mathbb{P}^{R,-}_{n}$
and $\mathbb{D}^{R}_{i}=0$ for $i\notin\mathcal{R}$, we obtain
\begin{equation}
\begin{aligned}\mathbb{P}^{R,-}_{n}\left(v_{n}\circ \widehat{Y}_{n}\right)-v_{0}\left(Y_{0}\right)={} & \sum^{n}_{i=1}\left\langle \mathbb{P}^{R,+}_{i-1},Q^{\mathrm{Alg}}_{i-1\leftarrow i}\left(v_{i}\circ \widehat{Y}_{i}\right)-v_{i-1}\circ \widehat{Y}_{i-1}\right\rangle \\
 & +\sum_{i\in\mathcal{R}}\left\langle \mathbb{D}^{R}_{i},v_{i}\right\rangle .
\end{aligned}
\label{eq:ver1-035}
\end{equation}
For $v_{i}=Q^{\mathrm{Alg}}_{(i,n]}f$, each term of the first sum vanishes
by \eqref{eq:ver1-007} and \eqref{eq:discrete-continuation}, and
\eqref{eq:ver1-035} for the whole tree, where no reduction is made,
gives $v_{0}\left(Y_{0}\right)=\mathbb{P}^{\mathrm{WT}}_{n}\left(f\circ \widehat{Y}_{n}\right)$.
This proves item~1. For $v_{i}=h^{\alpha}_{i}$, \eqref{eq:ver1-035}
is the total error identity, which is the semigroup telescoping used
in \citet[Theorem~19]{litterer_high_2012}.

For item~2, let $y$ be a state of $U_{i,\kappa}$, and put $z=y-c_{i,\kappa}$.
Since $c_{i,\kappa}$ and $y$ lie in the convex set $\mathrm{C.H.}(U_{i,\kappa})$,
so does the segment between them, and Taylor's formula with integral
remainder gives
\begin{equation}
\begin{aligned}h\left(c_{i,\kappa}+z\right)={} & \sum^{G+1}_{q=0}\frac{1}{q!}D^{q}h\left(c_{i,\kappa}\right)\left[z,\ldots,z\right]\\
 & +\frac{1}{\left(G+1\right)!}\int^{1}_{0}\left(1-t\right)^{G+1}D^{G+2}h\left(c_{i,\kappa}+tz\right)\left[z,\ldots,z\right]dt.
\end{aligned}
\label{eq:ver1-036}
\end{equation}
By \eqref{eq:moments}, the terms of degrees zero through $G$ integrate
to zero against $\widehat{Y}_{i\#}\mathbb{D}^{R}_{i,\kappa}$, which
is supported by the states of $U_{i,\kappa}$, and multilinearity gives
\begin{equation}
\begin{aligned} & \int D^{G+1}h\left(c_{i,\kappa}\right)\left[z,\ldots,z\right]\left(\widehat{Y}_{i\#}\mathbb{D}^{R}_{i,\kappa}\right)\left(dy\right)\\
 & \quad=w_{i,\kappa}\left(r^{{M}}_{i,\kappa}\right)^{G+1}\left\langle D^{G+1}h\left(c_{i,\kappa}\right)\circ\left({M}^{-1}_{i}\right)^{\otimes\left(G+1\right)},\mathcal{Q}^{{M}}_{i,\kappa,G+1}\right\rangle _{\otimes}.
\end{aligned}
\label{eq:ver1-037}
\end{equation}
If $r^{{M}}_{i,\kappa}=0$, all states of the patch coincide with $c_{i,\kappa}$,
and both sides of \eqref{eq:ver1-037} vanish. Since
$c_{i,\kappa}+tz\in\mathrm{C.H.}(\operatorname{supp}(\mathbb{P}^{R,-}_{i}))$ for $0\le t\le1$ and $\left\lVert M_{i}z\right\rVert _{2}\le r^{{M}}_{i,\kappa}$,
the quantity $\left|D^{G+2}h\left(c_{i,\kappa}+tz\right)\left[z,\ldots,z\right]\right|$
is at most $\left(r^{{M}}_{i,\kappa}\right)^{G+2}$ times the supremum
in \eqref{eq:ver1-032}.
Together with $\int^{1}_{0}\left(1-t\right)^{G+1}dt=1/\left(G+2\right)$
and the bound $2w_{i,\kappa}$ on the total variation of $\widehat{Y}_{i\#}\mathbb{D}^{R}_{i,\kappa}$,
this shows that the integral of the remainder term of \eqref{eq:ver1-036}
against this measure is at most \eqref{eq:ver1-032} in absolute value,
which proves item~2.

For item~3, let $\alpha\in\mathcal{I}$. Taylor's formula is now applied
to $h^{\alpha}_{i}$ only through degree $G$, so that the integral
remainder involves $D^{G+1}h^{\alpha}_{i}$. The terms of degrees
zero through $G$ integrate to zero as before,
and the definition~\eqref{eq:derivative-envelope} gives
\begin{equation}
\begin{aligned}\left|\left\langle \mathbb{D}^{R}_{i,\kappa},h^{\alpha}_{i}\right\rangle \right| & \le\frac{\varepsilon_{\alpha}L^{\left(G+1\right)}_{i}\left(\Gamma_{i}\right)}{G!}\int\left\lVert {M}_{i}z\right\rVert ^{G+1}_{2}\left|\widehat{Y}_{i\#}\mathbb{D}^{R}_{i,\kappa}\right|\left(dy\right)\int^{1}_{0}\left(1-t\right)^{G}dt\\
 & \le\frac{2\varepsilon_{\alpha}w_{i,\kappa}}{\left(G+1\right)!}L^{\left(G+1\right)}_{i}\left(\Gamma_{i}\right)\left(r^{{M}}_{i,\kappa}\right)^{G+1}.
\end{aligned}
\label{eq:ver1-039}
\end{equation}
This is the local Taylor cancellation of \citet[Lemma~15 and Proposition~16]{litterer_high_2012},
with the patch mass factor of \citet[Theorem~3.1]{ninomiya_high-order_2025}.
Since $h^{\alpha}_{n}=P_{0}f_{\alpha}=f_{\alpha}$,
inserting \eqref{eq:ver1-039}
into \eqref{eq:ver1-035} with $v_{i}=h^{\alpha}_{i}$, dividing by
$\varepsilon_{\alpha}$ and taking the maximum over $\alpha\in\mathcal{I}$,
we obtain \eqref{eq:joint-error}.
A patch with $\widehat{Y}_{i\#}\mathbb{D}^{R}_{i,\kappa}=0$
contributes nothing to the last sum of \eqref{eq:ver1-035}, which
proves the last assertion.
\end{proof}

The first sum on the right side of \eqref{eq:joint-error} is called
the propagation term and is denoted by $E_{\mathrm{weak}}$. For each
patch $U_{i,\kappa}$, we call the term $2L^{\left(G+1\right)}_{i}\left(\Gamma_{i}\right)w_{i,\kappa}\left(r^{{M}}_{i,\kappa}\right)^{G+1}/\left(G+1\right)!$
of \eqref{eq:joint-error} the certificate of the patch. Unlike \citet[Theorem~3.1]{ninomiya_high-order_2025},
the bound \eqref{eq:joint-error} keeps the propagation term and the
envelope \eqref{eq:derivative-envelope} as they are, and does not
bound them by the Kusuoka estimates used there. Identity \eqref{eq:discrete-residual}
is signed, whereas \eqref{eq:joint-error} bounds the absolute values
of the propagation and recombination contributions separately. When
$\Gamma_{i}$ is replaced by $a\Gamma_{i}$ with $a>0$, the envelope
$L^{\left(G+1\right)}_{i}\left(\Gamma_{i}\right)$ is multiplied by
$a^{\left(G+1\right)/2}$ and $r^{{M}}_{i,\kappa}$ by $a^{-1/2}$,
so the certificate of each patch does not change. Like $\lambda$ (Lemma~\ref{lem:metric-scalar}),
the scale of $\Gamma_{i}$ therefore affects only which patches are
accepted.

Suppose now that the patches are produced by the recursive subdivision
with the matrix test \eqref{eq:partition-rule} and the population
rule, and that Theorem~\ref{thm:residual} is applied with the matrices
$\Gamma_{i}$ of the test. The right side of \eqref{eq:partition-rule}
is referred to as the allowance of the patch.
For a patch that passes \eqref{eq:partition-rule},
the choice $c_{i,\kappa}\in\mathrm{C.H.}(U_{i,\kappa})$
gives $r^{{M}}_{i,\kappa}\le\operatorname{rad}_{M_{i}}\left(U_{i,\kappa}\right)$,
and hence $w_{i,\kappa}\left(r^{{M}}_{i,\kappa}\right)^{G+1}\le\lambda^{G+1}\left(T-t_{i}\right)^{G/2}s^{\left(G+1\right)/2}_{i}$,
in which the patch mass has canceled.
A patch that the population rule stops without passing \eqref{eq:partition-rule}
satisfies $\widehat{Y}_{i\#}\mathbb{D}^{R}_{i,\kappa}=0$ whenever
the moment columns of its distinct states are linearly independent,
since every feasible weight vector then has the same total weight
at each state as the input. We assume this independence for every
such patch. It holds for states in general position, since the patch
has at most $M_{G}$ nodes, but it can fail when the states lie on
a lattice (Section~\ref{sec:shape-test}). By item~3 of Theorem~\ref{thm:residual},
the inner sum in \eqref{eq:joint-error} may then be restricted to
patches that pass \eqref{eq:partition-rule}. Since there are at most
$K_{i}$ such patches at date $i$, we obtain
\begin{equation}
\begin{aligned} & \max_{\alpha\in\mathcal{I}}\frac{\left|\mathbb{P}^{R,-}_{n}\left(f_{\alpha}\circ \widehat{Y}_{n}\right)-h^{\alpha}_{0}\left(Y_{0}\right)\right|}{\varepsilon_{\alpha}}\\
 & \quad\le E_{\mathrm{weak}}+\frac{2\lambda^{G+1}}{\left(G+1\right)!}\sum_{i\in\mathcal{R}}K_{i}L^{\left(G+1\right)}_{i}\left(\Gamma_{i}\right)\left(T-t_{i}\right)^{G/2}s^{\left(G+1\right)/2}_{i}.
\end{aligned}
\label{eq:joint-wms-budget}
\end{equation}
The bound \eqref{eq:joint-error} retains the actual patch radii. In
\eqref{eq:joint-wms-budget} they are replaced by the allowance, and
the recombination contribution is carried by the patch counts $K_{i}$
and the envelopes $L^{\left(G+1\right)}_{i}\left(\Gamma_{i}\right)$.

Suppose that, for some constants $C_{\mathrm{w}}<\infty$ and $p_{\mathrm{w}}>0$
independent of the grid, the normalized local weak defect $\max_{\alpha\in\mathcal{I}}\varepsilon^{-1}_{\alpha}\left|Q^{\mathrm{Alg}}_{i}h^{\alpha}_{i}-h^{\alpha}_{i-1}\right|$
is at most $C_{\mathrm{w}}s^{p_{\mathrm{w}}+1}_{i}$ at every state
in the support of $\widehat{Y}_{i-1\#}\mathbb{P}^{R,+}_{i-1}$ for $1\le i\le n$.
Since $\mathbb{P}^{R,+}_{i-1}$ is a probability measure,
\eqref{eq:ver1-007}
then gives $E_{\mathrm{weak}}\le C_{\mathrm{w}}T\left(\max_{i}s_{i}\right)^{p_{\mathrm{w}}}$.
If
\begin{equation}
\sum_{i\in\mathcal{R}}K_{i}L^{\left(G+1\right)}_{i}\left(\Gamma_{i}\right)\left(T-t_{i}\right)^{G/2}s^{\left(G+1\right)/2}_{i}\le C_{\mathrm{R}}\left(\max_{j}s_{j}\right)^{p_{\mathrm{w}}},\label{eq:matched-order-condition}
\end{equation}
for a constant $C_{\mathrm{R}}<\infty$ independent of the grid, then
the normalized joint error is at most
\begin{equation}
\left(C_{\mathrm{w}}T+\frac{2\lambda^{G+1}C_{\mathrm{R}}}{\left(G+1\right)!}\right)\left(\max_{i}s_{i}\right)^{p_{\mathrm{w}}}.\label{eq:ver1-043}
\end{equation}
Recombination therefore preserves the weak order $p_{\mathrm{w}}$ of
the propagation whenever \eqref{eq:matched-order-condition} holds.

For example, suppose that, for some constants $C_{KL}>0$ and $p_{K}\ge0$
independent of the grid,
\begin{equation}
\max_{i\in\mathcal{R}}K_{i}\left(T-t_{i}\right)^{G/2}L^{\left(G+1\right)}_{i}\left(\Gamma_{i}\right)\le C_{KL}\left(\max_{j}s_{j}\right)^{-p_{K}}.\label{eq:ver1-044}
\end{equation}
Then, since $\sum_{i}s_{i}=T$,
\begin{equation}
\begin{aligned} & \sum_{i\in\mathcal{R}}K_{i}L^{\left(G+1\right)}_{i}\left(\Gamma_{i}\right)\left(T-t_{i}\right)^{G/2}s^{\left(G+1\right)/2}_{i}\\
 & \qquad\le C_{KL}T\left(\max_{j}s_{j}\right)^{\left(G-1\right)/2-p_{K}}.
\end{aligned}
\label{eq:ver1-045}
\end{equation}
Hence, if $p_{\mathrm{w}}\le\left(G-1\right)/2-p_{K}$, then
\eqref{eq:matched-order-condition}
holds with $C_{\mathrm{R}}=C_{KL}T^{\left(G+1\right)/2-p_{K}-p_{\mathrm{w}}}$,
since $\max_{j}s_{j}\le T$. The exponent $p_{K}$ describes the combined
growth of the patch counts and of the continuation derivatives. Proposition~\ref{prop:wms-refinement}
gives the corresponding bound for the matrix selected in Section~\ref{sec:method},
under a weighted form of \eqref{eq:ver1-044}
with $L^{\left(G+1\right)}_{i}\left(I_{N}\right)$
in place of $L^{\left(G+1\right)}_{i}\left(\Gamma_{i}\right)$.

\subsubsection{Task-aligned certificate volume}

For \(i\in\mathcal{R}\) and \(z\in\mathbb{R}^{N}\), we put, under the continuity
assumption made for \eqref{eq:derivative-envelope},
\begin{equation}
\mathcal{N}_{i}\left(z\right)=\max_{\alpha\in\mathcal{I}}\frac{1}{\varepsilon_{\alpha}}
\sup_{x\in\mathrm{C.H.}(\operatorname{supp}(\mathbb{P}^{R,-}_{i}))}\left|D^{G+1}h^{\alpha}_{i}\left(x\right)\left[z,\ldots,z\right]\right|.
\label{eq:task-gauge}
\end{equation}
The function \(\mathcal{N}_{i}\) is continuous, even and absolutely homogeneous
of degree \(G+1\), and \eqref{eq:derivative-envelope} gives
\(L^{\left(G+1\right)}_{i}\left(\Gamma\right)=\sup_{\lVert v\rVert_{2}=1}\mathcal{N}_{i}\left(\Gamma^{1/2}v\right)\).

\begin{defn}
\label{def:task-aligned}Let \(w>0\) be a patch mass, and let \(\varepsilon_{\mathrm{loc}}>0\)
be a tolerance for the certificate of a patch. The certificate region is
\begin{equation}
\mathcal{C}_{i}=\left\{ z\in\mathbb{R}^{N}:
\mathcal{N}_{i}\left(z\right)\le\frac{\left(G+1\right)!\,\varepsilon_{\mathrm{loc}}}{2w}\right\},
\label{eq:certificate-region}
\end{equation}
and, for \(\Gamma\succ0\), the certificate volume \(\mathcal{V}_{i}\left(\Gamma\right)\)
is the largest volume of a centered ellipsoid of shape \(\Gamma\) contained
in \(\mathcal{C}_{i}\),
\begin{equation}
\mathcal{V}_{i}\left(\Gamma\right)=\sup\left\{ \operatorname{vol}
\left(E_{\Gamma}\left(r\right)\right):r>0,\;E_{\Gamma}\left(r\right)\subseteq\mathcal{C}_{i}\right\},
\qquad E_{\Gamma}\left(r\right)=\left\{ z:\left\lVert \Gamma^{-1/2}z\right\rVert _{2}\le r\right\} .
\label{eq:certificate-volume-def}
\end{equation}
\end{defn}
The certificate volume plays, for centered ellipsoids, the role of the
shape function of anisotropic piecewise polynomial approximation, which
is the smallest local error of a homogeneous polynomial over triangles
of unit area \citep{cohen_mirebeau_2011}.

\begin{prop}
\label{prop:task-aligned}Assume that \(\mathcal{N}_{i}\) does not vanish identically,
and let \(\Gamma\succ0\). Then \(L^{\left(G+1\right)}_{i}\left(\Gamma\right)>0\),
and \(E_{\Gamma}\left(r\right)\subseteq\mathcal{C}_{i}\) if and only if
\(r^{G+1}L^{\left(G+1\right)}_{i}\left(\Gamma\right)\le\left(G+1\right)!\varepsilon_{\mathrm{loc}}/\left(2w\right)\).
The supremum in \eqref{eq:certificate-volume-def} is attained, and
\begin{equation}
\mathcal{V}_{i}\left(\Gamma\right)=v_{N}\det\left(\Gamma\right)^{1/2}
\left[\frac{\left(G+1\right)!\varepsilon_{\mathrm{loc}}}{2wL^{\left(G+1\right)}_{i}\left(\Gamma\right)}\right]^{N/\left(G+1\right)},
\label{eq:certificate-volume}
\end{equation}
where \(v_{N}\) is the volume of the Euclidean unit ball. Moreover,
\(\mathcal{V}_{i}\left(a\Gamma\right)=\mathcal{V}_{i}\left(\Gamma\right)\) for
every \(a>0\), so the certificate volume depends only on the shape
of \(\Gamma\). If, in addition, \(\mathcal{N}_{i}\left(z\right)>0\) for every
\(z\neq0\), then \(\mathcal{C}_{i}\) is compact and \(\mathcal{V}_{i}\) attains
its maximum over \(\Gamma\succ0\). We call any maximizer an ideal task-aligned
shape and denote it by \(\Gamma^{\mathrm{TA}}_{i}\). The corresponding ellipsoid
is a maximum-volume centered ellipsoid contained in \(\mathcal{C}_{i}\).
\end{prop}

\begin{proof}
Let \(\mathcal{N}_{i}\left(z\right)>0\). For \(v=\Gamma^{-1/2}z/\lVert\Gamma^{-1/2}z\rVert_{2}\),
homogeneity gives \(\mathcal{N}_{i}\left(\Gamma^{1/2}v\right)=\mathcal{N}_{i}\left(z\right)/\lVert\Gamma^{-1/2}z\rVert^{G+1}_{2}>0\),
so \(L^{\left(G+1\right)}_{i}\left(\Gamma\right)>0\). Also by homogeneity,
\(\sup_{z\in E_{\Gamma}\left(r\right)}\mathcal{N}_{i}\left(z\right)
=r^{G+1}\sup_{\lVert v\rVert_{2}=1}\mathcal{N}_{i}\left(\Gamma^{1/2}v\right)
=r^{G+1}L^{\left(G+1\right)}_{i}\left(\Gamma\right)\), which gives the inclusion
criterion and shows that the largest admissible radius is
\begin{equation}
r=\left[\frac{\left(G+1\right)!\varepsilon_{\mathrm{loc}}}
{2wL^{\left(G+1\right)}_{i}\left(\Gamma\right)}\right]^{1/\left(G+1\right)}.
\label{eq:certificate-radius}
\end{equation}
The change of variables \(z=\Gamma^{1/2}u\) gives
\(\operatorname{vol}\left(E_{\Gamma}\left(r\right)\right)=v_{N}r^{N}\det\left(\Gamma\right)^{1/2}\),
and substitution gives \eqref{eq:certificate-volume}. Since
\(L^{\left(G+1\right)}_{i}\left(a\Gamma\right)=a^{\left(G+1\right)/2}
L^{\left(G+1\right)}_{i}\left(\Gamma\right)\) and
\(\det\left(a\Gamma\right)^{1/2}=a^{N/2}\det\left(\Gamma\right)^{1/2}\), the two
factors of \(a^{N/2}\) cancel. Finally, if \(\mathcal{N}_{i}\) is positive off
the origin, then it attains a positive minimum on the unit sphere, so
\(\mathcal{C}_{i}\) is bounded, and it is closed by continuity. Every ellipsoid
\(E_{\Gamma}\left(r\right)\subseteq\mathcal{C}_{i}\) then satisfies \(r\lambda_{\max}\left(\Gamma\right)^{1/2}\le R\) for a
radius \(R\) with \(\mathcal{C}_{i}\subseteq B\left(0,R\right)\). Scale
invariance permits the normalization \(\det\Gamma=1\), under which
\(\mathcal{V}_{i}\left(\Gamma\right)=v_{N}r^{N}\) with \(r\) from
\eqref{eq:certificate-radius}. Along a maximizing sequence eventually
\(\mathcal{V}_{i}\left(\Gamma\right)\ge\mathcal{V}_{i}\left(I_{N}\right)/2\), so \(r\)
is bounded below, \(\lambda_{\max}\left(\Gamma\right)\le\left(R/r\right)^{2}\) is
bounded above and \(\lambda_{\min}\left(\Gamma\right)\ge\lambda_{\max}\left(\Gamma\right)^{-\left(N-1\right)}\)
is bounded below. The sequence therefore stays in a compact subset of the
positive definite cone, and continuity of \(\Gamma\mapsto\mathcal{V}_{i}\left(\Gamma\right)\)
gives an attained maximum.
\end{proof}
By the first assertion of Proposition~\ref{prop:task-aligned}, a patch of
mass \(w\) whose centered metric radius for \(M_{i}=\Gamma^{-1/2}\) is \(r\)
has certificate at most \(\varepsilon_{\mathrm{loc}}\) if and only if
\(E_{\Gamma}\left(r\right)\subseteq\mathcal{C}_{i}\).

\begin{cor}
\label{cor:task-aligned-ellipsoid}Suppose that there exists \(\Gamma^{\ast}_{i}\succ0\)
such that
\begin{equation}
\mathcal{N}_{i}\left(z\right)=\left(\prescript{t}{}{z}\left(\Gamma^{\ast}_{i}\right)^{-1}z\right)^{\left(G+1\right)/2}
\qquad\text{for every }z\in\mathbb{R}^{N}.
\label{eq:task-scale-ellipsoid}
\end{equation}
Then \(\mathcal{C}_{i}\) is an ellipsoid, \(\Gamma\succ0\) maximizes \(\mathcal{V}_{i}\)
if and only if \(\Gamma\) is a positive multiple of \(\Gamma^{\ast}_{i}\), and
\begin{equation}
\frac{\mathcal{V}_{i}\left(\Gamma^{\mathrm{TA}}_{i}\right)}{\mathcal{V}_{i}\left(I_{N}\right)}
=\sqrt{\frac{\det\Gamma^{\ast}_{i}}{\lambda_{\min}\left(\Gamma^{\ast}_{i}\right)^{N}}}
=\prod^{N}_{j=1}\left(\frac{\lambda_{j}\left(\Gamma^{\ast}_{i}\right)}
{\lambda_{\min}\left(\Gamma^{\ast}_{i}\right)}\right)^{1/2}.
\label{eq:task-aligned-gain}
\end{equation}
The ratio equals one if and only if \(\Gamma^{\ast}_{i}\) is a multiple of the
identity, and it is at least \(\left(\lambda_{\max}\left(\Gamma^{\ast}_{i}\right)/\lambda_{\min}\left(\Gamma^{\ast}_{i}\right)\right)^{1/2}\).
\end{cor}

\begin{proof}
Put \(B=\left(\Gamma^{\ast}_{i}\right)^{-1/2}\Gamma^{1/2}\), an invertible
matrix with \(\det B>0\). Under
\eqref{eq:task-scale-ellipsoid},
\(L^{\left(G+1\right)}_{i}\left(\Gamma\right)
=\lambda_{\max}\left(\prescript{t}{}{B}B\right)^{\left(G+1\right)/2}
=\left\lVert B\right\rVert^{G+1}_{2}\), while
\(\det\left(\Gamma\right)^{1/2}=\det\left(\Gamma^{\ast}_{i}\right)^{1/2}\det B\).
Hence \eqref{eq:certificate-volume} gives
\begin{equation}
\mathcal{V}_{i}\left(\Gamma\right)=v_{N}
\left[\frac{\left(G+1\right)!\varepsilon_{\mathrm{loc}}}{2w}\right]^{N/\left(G+1\right)}
\det\left(\Gamma^{\ast}_{i}\right)^{1/2}\,
\frac{\det B}{\left\lVert B\right\rVert^{N}_{2}}.
\label{eq:certificate-volume-B}
\end{equation}
Since \(\det B=\prod_{j}\sigma_{j}\left(B\right)\le\left\lVert B\right\rVert^{N}_{2}\)
with equality if and only if every singular value equals
\(\left\lVert B\right\rVert_{2}\), that is,
\(\prescript{t}{}{B}B=\Gamma^{1/2}\left(\Gamma^{\ast}_{i}\right)^{-1}\Gamma^{1/2}=a^{2}I_{N}\)
for some \(a>0\), the volume \eqref{eq:certificate-volume-B} is maximal
if and only if \(\Gamma=a^{2}\Gamma^{\ast}_{i}\). Evaluating
\eqref{eq:certificate-volume-B} at \(\Gamma=\Gamma^{\ast}_{i}\), where
\(B=I_{N}\), and at \(\Gamma=I_{N}\), where
\(B=\left(\Gamma^{\ast}_{i}\right)^{-1/2}\), for which
\(\det B=\det\left(\Gamma^{\ast}_{i}\right)^{-1/2}\) and
\(\left\lVert B\right\rVert_{2}=\lambda_{\min}\left(\Gamma^{\ast}_{i}\right)^{-1/2}\),
gives \eqref{eq:task-aligned-gain}.
\end{proof}
\begin{rem}
\label{rem:shape-ceiling}Corollary~\ref{cor:task-aligned-ellipsoid} measures
the potential gain of shape alignment by
$\prod_{j}\left(\lambda_{j}(\Gamma^{\ast}_{i})/\lambda_{\min}(\Gamma^{\ast}_{i})\right)^{1/2}$,
which is not bounded over task families. The change produced by a given metric is limited by
its own spectrum for every task family. Let $\gamma_{1}\ge\cdots\ge\gamma_{N}>0$
be the eigenvalues of $\Gamma$. By \eqref{eq:envelope-sandwich},
$L^{(G+1)}_{i}(\Gamma)/L^{(G+1)}_{i}(I_{N})$ lies between $\gamma^{(G+1)/2}_{N}$
and $\gamma^{(G+1)/2}_{1}$, and \eqref{eq:certificate-volume} gives
\begin{equation}
\prod^{N}_{j=1}\left(\frac{\gamma_{j}}{\gamma_{1}}\right)^{1/2}\le\frac{\mathcal{V}_{i}(\Gamma)}{\mathcal{V}_{i}(I_{N})}\le\prod^{N}_{j=1}\left(\frac{\gamma_{j}}{\gamma_{N}}\right)^{1/2}.\label{eq:shape-ceiling}
\end{equation}
Equality holds on the right when $\mathcal{N}_{i}(z)$ depends only on the
projection of $z$ onto the eigenspace of $\gamma_{N}$, the narrowest axes
of $E_{\Gamma}(r)$, and on the left when it depends only on the projection
onto the eigenspace of $\gamma_{1}$, the widest axes. If
$\operatorname{cond}(\Gamma)=\gamma_{1}/\gamma_{N}\le\rho$, the ratio lies
between $\rho^{-(N-1)/2}$ and $\rho^{(N-1)/2}$, so such a metric captures
only part of a larger potential gain.
Under \eqref{eq:task-scale-ellipsoid}
the upper bound sharpens to
$\prod_{j}\min\{\lambda_{j}(\Gamma^{\ast}_{i})/\lambda_{\min}(\Gamma^{\ast}_{i}),\rho\}^{1/2}$,
attained by the coaxial matrix whose eigenvalues are those of
$\Gamma^{\ast}_{i}$ clipped at $\rho\lambda_{\min}(\Gamma^{\ast}_{i})$. Indeed,
since $\mathcal{V}_{i}$ is invariant under scaling, we may assume that
$\Gamma\preceq\Gamma^{\ast}_{i}$ with equality in some direction. Then
\eqref{eq:certificate-volume} gives
$\mathcal{V}_{i}(\Gamma)/\mathcal{V}_{i}(I_{N})=\prod_{j}(\gamma_{j}/\lambda_{\min}(\Gamma^{\ast}_{i}))^{1/2}$,
and Weyl monotonicity and the condition bound give
$\gamma_{j}\le\min\{\lambda_{j}(\Gamma^{\ast}_{i}),\rho\lambda_{\min}(\Gamma^{\ast}_{i})\}$.
In the one-block configuration, \eqref{eq:ridge-condition} gives
$\rho=\varphi^{2}$ with the golden ratio $\varphi$, so the AWLL shape changes the certificate volume of
the identity by at most the factor $\varphi^{N-1}$,
about $11.1$ for $N=6$, in either direction. With several blocks,
\eqref{eq:ridge-condition} bounds each block but not the ratios between blocks.
Section~\ref{sec:shape-scale} compares these bounds with the experiments.
\end{rem}
\FloatBarrier

\section{Method}

\label{sec:method}

\noindent
\label{sec:wms-method}This section constructs the matrix $\Gamma^{\mathrm{AWLL}}_{i}$
that AWLL uses in the matrix test \eqref{eq:partition-rule}, and
Section~\ref{sec:solver} computes the positive reduction
\eqref{eq:moments}.

\subsection{Forward references and closed-form metric selection}
\noindent
Fix a recombination date $i\in\mathcal{R}$ and nonempty coordinate
blocks $\left\{ {C}_{c}\right\} ^{L}_{c=1}$ partitioning $\left\{ 1,\ldots,N\right\} $.
The edges into date $i$ are the nodes $\left[e\right]$, $e=\left(j,m\right)$,
in the support of $\mathbb{P}^{R,-}_{i}$, where $\left[j\right]$
is a node of date $i-1$. By \eqref{eq:tree-recursion}, an edge has
probability
\begin{equation}
\varpi_{e}=p^{R,-}_{i,e}=p^{R,+}_{i-1,j}\pi_{i,m}\left(\widehat{Y}_{i-1}\left(j\right)\right),\qquad\sum_{e}\varpi_{e}=1,\label{eq:ver1-047}
\end{equation}
and its block increment is $\Delta_{e,c}=(\widehat{Y}_{i}(e)-\widehat{Y}_{i-1}(j))_{C_{c}}$.
The probability-weighted center of the increments is
\begin{equation}
\Delta_{i,c}=\sum_{e}\varpi_{e}\Delta_{e,c}.\label{eq:weighted-center}
\end{equation}
The covariance per unit time of the centered edge displacements is
\begin{equation}
\widehat{\Gamma}^{\mathrm{edge}}_{i,c}=\frac{1}{s_{i}}\sum_{e}\varpi_{e}\left(\Delta_{e,c}-\Delta_{i,c}\right)\prescript{t}{}{\left(\Delta_{e,c}-\Delta_{i,c}\right)}.\label{eq:edge-proxy}
\end{equation}
This matrix is called the edge reference.

The second reference, called the model reference, retains both diffusion
and predictable displacement. With $A=\Sigma\prescript{t}{}{\Sigma}$
from \eqref{eq:augmented}, it is defined by
\begin{equation}
\Gamma^{\mathrm{mod}}_{i,c}=\mathbb{P}^{R,+}_{i-1}\left(A_{C_{c}C_{c}}+s_{i}b_{C_{c}}\prescript{t}{}{b_{C_{c}}}\right).\label{eq:model-proxy}
\end{equation}
This is the second moment per unit time of the $C_{c}$ subvector of
the frozen increment $s_{i}b\left(\widehat{Y}_{i-1}\right)+\sqrt{s_{i}}\Sigma\left(\widehat{Y}_{i-1}\right)Z$,
where $Z$ is a standard Gaussian vector in $\mathbb{R}^{d_{B}}$ independent
of $\widehat{Y}_{i-1}$. For $q\in\mathbb{R}^{|C_{c}|}$,
the quadratic form of $\Gamma^{\mathrm{mod}}_{i,c}$ is written with
the parent measure as follows:
\begin{equation}
\prescript{t}{}{q}\Gamma^{\mathrm{mod}}_{i,c}q=\mathbb{P}^{R,+}_{i-1}\left(\left\lVert \prescript{t}{}{\Sigma_{C_{c}\cdot}}q\right\rVert ^{2}_{2}+s_{i}\left(\prescript{t}{}{q}b_{C_{c}}\right)^{2}\right).\label{eq:model-displacement-energy}
\end{equation}
Consequently, a direction $q$ with $\prescript{t}{}{q}\Gamma^{\mathrm{mod}}_{i,c}q=0$
has neither diffusion nor drift at any state of the support of $\widehat{Y}_{i-1\#}\mathbb{P}^{R,+}_{i-1}$.
The term $s_{i}b_{C_{c}}\prescript{t}{}{b_{C_{c}}}$ gives a positive
scale to directions with transport but no diffusion, such as running
integrals, also at $i=1$, where the only parent state is $Y_{0}$.
A centered model reference would subtract $s_{i}\bar{b}_{i,C_{c}}\prescript{t}{}{\bar{b}_{i,C_{c}}}$,
where $\bar{b}_{i}=\mathbb{P}^{R,+}_{i-1}\left(b\right)$, and would
discard the common transport, as the edge reference does by centering.
Invertibility is supplied by the ridge \eqref{eq:ridge-design} below.

For each edge, we put
\begin{equation}
\Psi_{e,c}=\frac{\left(\Delta_{e,c}-\Delta_{i,c}\right)\prescript{t}{}{\left(\Delta_{e,c}-\Delta_{i,c}\right)}}{s_{i}},\label{eq:ver1-053}
\end{equation}
and, for a tuple $X=\left(X_{c}\right)^{L}_{c=1}$ of symmetric blocks,
$\left\lVert X\right\rVert ^{2}_{\mathcal{H}}=\sum_{c}\left\lVert X_{c}\right\rVert ^{2}_{F}$.
A symbol without the index $c$, such as $\widehat{\Gamma}^{\mathrm{edge}}_{i}$,
denotes the tuple of its blocks.
The edge dispersion and the reference discrepancy are
\begin{equation}
\widehat{V}_{i}=\sum_{c}\sum_{e}\varpi^{2}_{e}\left\lVert \Psi_{e,c}-\widehat{\Gamma}^{\mathrm{edge}}_{i,c}\right\rVert ^{2}_{F},\qquad\widehat{D}_{i}=\sum_{c}\left\lVert \widehat{\Gamma}^{\mathrm{edge}}_{i,c}-\Gamma^{\mathrm{mod}}_{i,c}\right\rVert ^{2}_{F}.\label{eq:sensitivity-objective}
\end{equation}
For signs $\sigma_{e}=\pm1$, replacing $\varpi_{e}$ by
$\varpi_{e}\left(1+t\sigma_{e}\right)/\left(1+t\sum_{f}\varpi_{f}\sigma_{f}\right)$
changes $\widehat{\Gamma}^{\mathrm{edge}}_{i}$ at $t=0$ with derivative
$\sum_{e}\varpi_{e}\sigma_{e}\left(\Psi_{e}-\widehat{\Gamma}^{\mathrm{edge}}_{i}\right)$,
because the moving center contributes nothing by \eqref{eq:weighted-center}.
If the signs are independent with equal probabilities, then
$\mathbb{E}\sigma_{e}\sigma_{f}=\delta_{ef}$, and $\widehat{V}_{i}$
is the mean of the squared norm $\left\lVert \cdot\right\rVert ^{2}_{\mathcal{H}}$
of this derivative.
The mixing coefficient $\vartheta^{\mathrm{AWLL}}_{i}$ and the mixed
blocks $\overline{\Gamma}_{i,c}$ are defined by the shrinkage rule of
\citet[Section~3.3]{ledoit_wolf_2004}, with $\widehat{V}_{i}$ as the
dispersion estimate and $\Gamma^{\mathrm{mod}}_{i}$ in place of a multiple
of the identity,
\begin{equation}
\begin{aligned}\vartheta^{\mathrm{AWLL}}_{i} & =\left(1-\frac{\widehat{V}_{i}}{\widehat{D}_{i}}\right)_{+}\quad\text{if }\widehat{D}_{i}>0,\qquad\vartheta^{\mathrm{AWLL}}_{i}=0\quad\text{if }\widehat{D}_{i}=0,\\
\overline{\Gamma}_{i,c} & =\vartheta^{\mathrm{AWLL}}_{i}\widehat{\Gamma}^{\mathrm{edge}}_{i,c}+\left(1-\vartheta^{\mathrm{AWLL}}_{i}\right)\Gamma^{\mathrm{mod}}_{i,c},
\end{aligned}
\label{eq:mixing-coefficient}
\end{equation}
where $x_{+}=\max\left\{ x,0\right\} $. Since $\overline{\Gamma}_{i}-\widehat{\Gamma}^{\mathrm{edge}}_{i}=\left(1-\vartheta^{\mathrm{AWLL}}_{i}\right)\left(\Gamma^{\mathrm{mod}}_{i}-\widehat{\Gamma}^{\mathrm{edge}}_{i}\right)$,
for $\widehat{D}_{i}>0$
\begin{equation}
\left\lVert \overline{\Gamma}_{i}-\widehat{\Gamma}^{\mathrm{edge}}_{i}\right\rVert _{\mathcal{H}}=\min\left\{ \sqrt{\widehat{D}_{i}},\frac{\widehat{V}_{i}}{\sqrt{\widehat{D}_{i}}}\right\} \le\sqrt{\widehat{V}_{i}},\label{eq:bounded-model-injection}
\end{equation}
and the left side vanishes for $\widehat{D}_{i}=0$. However large
$\widehat{D}_{i}$ is, the mixture $\overline{\Gamma}_{i}$ differs from
the edge reference by at most $\sqrt{\widehat{V}_{i}}$ in $\left\lVert \cdot\right\rVert _{\mathcal{H}}$. Since $\varpi^{2}_{e}\le\left(\max_{f}\varpi_{f}\right)\varpi_{e}$,
\begin{equation*}
\widehat{V}_{i}\le\left(\max_{f}\varpi_{f}\right)\sum_{c}\sum_{e}\varpi_{e}\left\lVert \Psi_{e,c}-\widehat{\Gamma}^{\mathrm{edge}}_{i,c}\right\rVert ^{2}_{F}.
\end{equation*}
If $\max_{f}\varpi_{f}$ tends to zero while the sum on the right stays
bounded, then $\widehat{V}_{i}$ and the left side of \eqref{eq:bounded-model-injection}
tend to zero, and $\vartheta^{\mathrm{AWLL}}_{i}$ tends to one if, in
addition, $\widehat{D}_{i}$ stays bounded away from zero. The model
reference therefore acts chiefly at early dates, including $i=1$,
where the edge reference rests on few branches.

The reference discrepancy has an explicit leading term. Let
$\mathcal{L}=\sum_{k}b_{k}\partial_{k}+\frac{1}{2}\sum_{j,k}A_{jk}\partial_{j}\partial_{k}$
be the generator of \eqref{eq:augmented}, let $Db$ be the Jacobian
of $b$, and put
\begin{equation}
A^{(2)}=\frac{1}{2}\left(\left(Db\right)A+A\prescript{t}{}{\left(Db\right)}+\mathcal{L}A\right),\label{eq:covariance-correction}
\end{equation}
with $\mathcal{L}$ applied entrywise. Second-order Dynkin expansions
of $Y^{y}_{s}-y$ and $\left(Y^{y}_{s}-y\right)\prescript{t}{}{\left(Y^{y}_{s}-y\right)}$
give the conditional mean $sb\left(y\right)+O\left(s^{2}\right)$
and the conditional covariance $sA\left(y\right)+s^{2}A^{(2)}\left(y\right)+O\left(s^{3}\right)$,
since the terms $s^{2}b\prescript{t}{}{b}$ of the second moment and of the
squared mean cancel.
\begin{prop}
\label{prop:reference-discrepancy}Fix $i\in\mathcal{R}$ and suppose
that, from every parent state $y$ in the support of $\widehat{Y}_{i-1\#}\mathbb{P}^{R,+}_{i-1}$,
the branch increments have conditional mean $s_{i}b\left(y\right)+O\left(s^{2}_{i}\right)$
and conditional covariance $s_{i}A\left(y\right)+s^{2}_{i}A^{(2)}\left(y\right)+O\left(s^{3}_{i}\right)$.
Then
\begin{equation}
\begin{aligned}\widehat{\Gamma}^{\mathrm{edge}}_{i,c} & =\mathbb{P}^{R,+}_{i-1}\left(A_{C_{c}C_{c}}\right)+s_{i}\left(\operatorname{Cov}_{\mathbb{P}^{R,+}_{i-1}}\left(b_{C_{c}}\right)+\mathbb{P}^{R,+}_{i-1}\left(A^{(2)}_{C_{c}C_{c}}\right)\right)+O\left(s^{2}_{i}\right),\\
\widehat{\Gamma}^{\mathrm{edge}}_{i,c}-\Gamma^{\mathrm{mod}}_{i,c} & =s_{i}\left(\mathbb{P}^{R,+}_{i-1}\left(A^{(2)}_{C_{c}C_{c}}\right)-\bar{b}_{i,C_{c}}\prescript{t}{}{\bar{b}_{i,C_{c}}}\right)+O\left(s^{2}_{i}\right),\\
\widehat{D}_{i} & =s^{2}_{i}\sum_{c}\left\lVert \mathbb{P}^{R,+}_{i-1}\left(A^{(2)}_{C_{c}C_{c}}\right)-\bar{b}_{i,C_{c}}\prescript{t}{}{\bar{b}_{i,C_{c}}}\right\rVert ^{2}_{F}+O\left(s^{3}_{i}\right).
\end{aligned}
\label{eq:edge-expansion}
\end{equation}
If $A_{kk}\equiv0$, as for a running integral, then $A^{(2)}_{kk}\equiv0$,
the $\left(k,k\right)$ entry of the leading coefficient is $-\bar{b}^{2}_{i,k}$,
and
\begin{equation}
\widehat{D}_{i}\ge s^{2}_{i}\sum_{k:A_{kk}\equiv0}\bar{b}^{4}_{i,k}+O\left(s^{3}_{i}\right).\label{eq:transport-floor}
\end{equation}
\end{prop}
\begin{proof}
For a parent $\left[j\right]$, let $\overline{\Delta}_{j}$ and $\operatorname{Cov}_{j}$
be the conditional mean and covariance of its branch increments. The
law of total covariance over the parent measure gives
\[
s_{i}\widehat{\Gamma}^{\mathrm{edge}}_{i,c}=\sum_{j}p^{R,+}_{i-1,j}\left(\operatorname{Cov}_{j}\right)_{C_{c}C_{c}}+\sum_{j}p^{R,+}_{i-1,j}\left(\overline{\Delta}_{j}-\Delta_{i}\right)_{C_{c}}\prescript{t}{}{\left(\overline{\Delta}_{j}-\Delta_{i}\right)_{C_{c}}},
\]
with $\Delta_{i}=\sum_{j}p^{R,+}_{i-1,j}\overline{\Delta}_{j}$. The
first sum is $s_{i}\mathbb{P}^{R,+}_{i-1}\left(A_{C_{c}C_{c}}\right)+s^{2}_{i}\mathbb{P}^{R,+}_{i-1}\left(A^{(2)}_{C_{c}C_{c}}\right)+O\left(s^{3}_{i}\right)$,
and the second is $s^{2}_{i}\operatorname{Cov}_{\mathbb{P}^{R,+}_{i-1}}\left(b_{C_{c}}\right)+O\left(s^{3}_{i}\right)$.
Dividing by $s_{i}$ gives the first line of \eqref{eq:edge-expansion}.
Since $\mathbb{P}^{R,+}_{i-1}\left(b\prescript{t}{}{b}\right)=\operatorname{Cov}_{\mathbb{P}^{R,+}_{i-1}}\left(b\right)+\bar{b}_{i}\prescript{t}{}{\bar{b}_{i}}$,
\eqref{eq:model-proxy} gives the second line, and squaring gives the
third. If $A_{kk}\equiv0$, positive semidefiniteness gives $A_{k\cdot}\equiv0$,
so every term of $A^{(2)}_{kk}$ vanishes. The Frobenius norm dominates
the diagonal entries, which proves \eqref{eq:transport-floor}. The
parent support is finite, so all remainders are uniform.
\end{proof}
The model reference therefore differs from the edge reference at first
order in $s_{i}$ by two terms: the second-order covariance correction
$A^{(2)}$, and the mean transport $\bar{b}_{i}\prescript{t}{}{\bar{b}_{i}}$
that $\Gamma^{\mathrm{mod}}_{i}$ retains by not centering.
Both NV rules of Section~\ref{sec:numerics} satisfy the hypothesis,
since for either rule $Q^{\mathrm{NV}}_{i}g(y)=P_{s_{i}}g(y)+O(s^{3}_{i})$
at every $y\in\mathcal{O}$ for the test functions $g(z)=z-y$ and
$g(z)=(z-y)\prescript{t}{}{(z-y)}$.

The mixed block is regularized by a ridge, that is, by the shift
$\overline{\Gamma}_{i,c}+\varepsilon I_{|C_{c}|}$ with $\varepsilon\ge0$,
as in ridge regression. The shift keeps the eigenvectors of
$\overline{\Gamma}_{i,c}$ and raises every eigenvalue by $\varepsilon$;
unlike the mixing step, it is not a convex combination with a target.
Let $\gamma^{+}_{i,c}$ and $\gamma^{-}_{i,c}$ be the largest and smallest
eigenvalues of $\overline{\Gamma}_{i,c}$. For a nonzero mixed block,
the shift is selected by
\begin{equation}
\varepsilon^{\mathrm{AWLL}}_{i,c}=\operatorname*{arg\,min}_{\varepsilon\ge0}\max\left\{ \frac{\gamma^{+}_{i,c}+\varepsilon}{\gamma^{+}_{i,c}},\frac{\gamma^{+}_{i,c}}{\gamma^{-}_{i,c}+\varepsilon}\right\} ,\label{eq:ridge-design}
\end{equation}
where the second term is $+\infty$ at $\gamma^{-}_{i,c}+\varepsilon=0$.
Relative to $\gamma^{+}_{i,c}$, the first term is the largest eigenvalue
of the shifted block, and the second is the reciprocal of its smallest
eigenvalue.
\begin{thm}
\label{thm:wms-derivation}\label{thm:spectral-scale}For a nonzero
mixed block, the minimizer in \eqref{eq:ridge-design} is unique and
is given by
\begin{equation}
\varepsilon^{\mathrm{AWLL}}_{i,c}=\frac{\sqrt{\left(\gamma^{+}_{i,c}-\gamma^{-}_{i,c}\right)^{2}+4\left(\gamma^{+}_{i,c}\right)^{2}}-\gamma^{+}_{i,c}-\gamma^{-}_{i,c}}{2}.\label{eq:closed-coefficients}
\end{equation}
For a zero mixed block, put $\varepsilon^{\mathrm{AWLL}}_{i,c}=1$,
and define the AWLL matrix and metric by
\begin{equation}
\begin{aligned}\Gamma^{\mathrm{AWLL}}_{i} & =\operatorname{blockdiag}_{c=1,\ldots,L}\left(\overline{\Gamma}_{i,c}+\varepsilon^{\mathrm{AWLL}}_{i,c}I_{\left|{C}_{c}\right|}\right),\qquad{M}^{\mathrm{AWLL}}_{i}=\left(\Gamma^{\mathrm{AWLL}}_{i}\right)^{-1/2}.\end{aligned}
\label{eq:wms-covariance}
\end{equation}
Then they are positive definite, and, with the golden ratio $\varphi$, for every nonzero mixed block
$\overline{\Gamma}_{i,c}$,
\begin{equation}
\varphi^{-1}\gamma^{+}_{i,c}I_{|{C}_{c}|}\preceq\overline{\Gamma}_{i,c}+\varepsilon^{\mathrm{AWLL}}_{i,c}I_{|{C}_{c}|}\preceq\varphi\gamma^{+}_{i,c}I_{|{C}_{c}|},\label{eq:scale-transmission}
\end{equation}
with both constants attained when $\gamma^{-}_{i,c}=0$, and consequently
\begin{equation}
\frac{\gamma^{+}_{i,c}+\varepsilon^{\mathrm{AWLL}}_{i,c}}{\gamma^{-}_{i,c}+\varepsilon^{\mathrm{AWLL}}_{i,c}}\le\varphi^{2}.\label{eq:ridge-condition}
\end{equation}
The constants depend neither on $N$ nor on $d_{B}$.
\end{thm}
\begin{proof}
We write $\gamma^{\pm}$ for $\gamma^{\pm}_{i,c}$, so that $\gamma^{+}>0$ and $\gamma^{-}\ge0$.
In \eqref{eq:ridge-design}, the first term increases strictly with
$\varepsilon$, the second decreases strictly, and at $\varepsilon=0$ the
second is at least the first. The unique minimizer is therefore their
crossing $\left(\gamma^{+}+\varepsilon\right)\left(\gamma^{-}+\varepsilon\right)=\left(\gamma^{+}\right)^{2}$,
whose nonnegative root is \eqref{eq:closed-coefficients}. Dividing
by $\left(\gamma^{+}\right)^{2}$ and putting $u=1-\gamma^{-}/\gamma^{+}\in\left[0,1\right]$
gives
\begin{equation}
\frac{\gamma^{+}+\varepsilon}{\gamma^{+}}=\frac{u+\sqrt{u^{2}+4}}{2}\in\left[1,\varphi\right],\qquad\frac{\gamma^{-}+\varepsilon}{\gamma^{+}}=\frac{\gamma^{+}}{\gamma^{+}+\varepsilon},\label{eq:ver1-067}
\end{equation}
with the upper value at $\gamma^{-}=0$. These are the extreme eigenvalues
of the regularized block divided by $\gamma^{+}$, which proves
\eqref{eq:scale-transmission},
and the condition ratio in \eqref{eq:ridge-condition} is the square
of the first ratio. In particular $\gamma^{-}+\varepsilon>0$.
A zero mixed block becomes
the identity. Hence $\Gamma^{\mathrm{AWLL}}_{i}$ and ${M}^{\mathrm{AWLL}}_{i}$
are positive definite.
\end{proof}

The ridge preserves the principal axes of each mixed block, and,
by \eqref{eq:ridge-condition}, the ratio of the longest to the shortest
of the corresponding semi-axes of the ellipsoid of the matrix test (Section~\ref{sec:directional}) is at most $\varphi$.
By the lower bound in \eqref{eq:scale-transmission}, each of these
semi-axes is at least $\varphi^{-1/2}$ times the longest one
before the shift, which limits the subdivision that a nearly singular
mixed block would cause. With the single block $L=1$ used in the
experiments, these bounds apply to the whole six-dimensional ellipsoid.

Since the partial derivatives of \eqref{eq:closed-coefficients} satisfy
$\partial_{\gamma^{+}}\varepsilon>0>\partial_{\gamma^{-}}\varepsilon$ and
$\partial_{\gamma^{+}}\varepsilon-\partial_{\gamma^{-}}\varepsilon\le3/\sqrt{5}$,
the mean value theorem and Weyl's inequality show that the ridge map
$H\mapsto H+\varepsilon\left(H\right)I_{|C_{c}|}$ is Lipschitz on nonzero
positive semidefinite blocks with constant $1+3/\sqrt{5}$ in operator
norm, sharp for $|C_{c}|\ge2$. With \eqref{eq:bounded-model-injection},
it follows that, whenever the edge and mixed blocks are nonzero, $\Gamma^{\mathrm{AWLL}}_{i}$
lies within $(1+3/\sqrt{5})\sqrt{\widehat{V}_{i}}$ in
operator norm of the matrix that \eqref{eq:wms-covariance} gives with
$\vartheta^{\mathrm{AWLL}}_{i}=1$.

\subsection{State covariance and grid refinement}
\noindent
Instead of $\Gamma^{\mathrm{AWLL}}_{i}$, we consider in
\eqref{eq:partition-rule}
the covariance of $\widehat{Y}_{i}$ under $\mathbb{P}^{R,-}_{i}$ divided
by $s_{i}$. This matrix has the units of a diffusion covariance, but
its numerator records the state dispersion accumulated up to $t_{i}$
rather than the next increment. For a probability measure $\mu$ on
$\widehat{\mathcal{F}}_{i}$, the covariance of $\widehat{Y}_{i}$ under
$\mu$ is denoted by $\operatorname{Cov}_{\mu}\left(\widehat{Y}_{i}\right)=\mu\left(\left(\widehat{Y}_{i}-\mu\left(\widehat{Y}_{i}\right)\right)\prescript{t}{}{\left(\widehat{Y}_{i}-\mu\left(\widehat{Y}_{i}\right)\right)}\right)$.
For a nonzero positive semidefinite $N\times N$ matrix $H$, let $\operatorname{Reg}\left(H\right)$
be the matrix that \eqref{eq:wms-covariance} gives with $L=1$ and
$\overline{\Gamma}_{i,1}=H$.
\begin{prop}
\label{prop:state-covariance-cancellation}Let $i\in\mathcal{R}$, suppose
that $\operatorname{Cov}_{\mathbb{P}^{R,-}_{i}}\left(\widehat{Y}_{i}\right)\ne0$,
and put $\Gamma_{i}=\operatorname{Reg}\left(s^{-1}_{i}\operatorname{Cov}_{\mathbb{P}^{R,-}_{i}}\left(\widehat{Y}_{i}\right)\right)$.
Then a patch $U_{i,\kappa}$ passes \eqref{eq:partition-rule} if and
only if
\begin{equation}
\operatorname{rad}_{\left[\operatorname{Reg}\left(\operatorname{Cov}_{\mathbb{P}^{R,-}_{i}}\left(\widehat{Y}_{i}\right)\right)\right]^{-1/2}}\left(U_{i,\kappa}\right)\le\lambda\left(T-t_{i}\right)^{G/\left(2\left(G+1\right)\right)}w^{-1/\left(G+1\right)}_{i,\kappa}.\label{eq:state-covariance-cancellation}
\end{equation}
If $\mathbb{P}^{R,-}_{i}$, $\widehat{Y}_{i}$, $t_{i}$ and $\lambda$
are fixed, then neither the acceptance decisions nor the patches produced
by the recursive subdivision of Section~\ref{sec:subdivision}, with
its population rule, depend on $s_{i}$. Moreover, the date $i$ term
of the recombination contribution in \eqref{eq:joint-wms-budget} equals
$2\lambda^{G+1}K_{i}L^{\left(G+1\right)}_{i}\left(\operatorname{Reg}\left(\operatorname{Cov}_{\mathbb{P}^{R,-}_{i}}\left(\widehat{Y}_{i}\right)\right)\right)\left(T-t_{i}\right)^{G/2}/\left(G+1\right)!$,
which contains no power of $s_{i}$.
\end{prop}
\begin{proof}
For $a>0$, the extreme eigenvalues of $aH$ are $a$ times those
of $H$, so \eqref{eq:closed-coefficients} gives
$\operatorname{Reg}\left(aH\right)=a\operatorname{Reg}\left(H\right)$.
Hence $\Gamma^{-1/2}_{i}=\sqrt{s_{i}}\left[\operatorname{Reg}\left(\operatorname{Cov}_{\mathbb{P}^{R,-}_{i}}\left(\widehat{Y}_{i}\right)\right)\right]^{-1/2}$,
and canceling $\sqrt{s_{i}}$ in \eqref{eq:partition-rule}
proves \eqref{eq:state-covariance-cancellation}.
Compared with $\left[\operatorname{Reg}\left(\operatorname{Cov}_{\mathbb{P}^{R,-}_{i}}\left(\widehat{Y}_{i}\right)\right)\right]^{-1/2}$,
the metric $\Gamma^{-1/2}_{i}$ multiplies every projection used by
the splitter by $\sqrt{s_{i}}$, which leaves its principal axes, order
and median unchanged, and the population rule does not involve $s_{i}$.
Induction over rejected candidates therefore proves the assertion on
the subdivision. Since $L^{\left(G+1\right)}_{i}\left(a\Gamma\right)=a^{\left(G+1\right)/2}L^{\left(G+1\right)}_{i}\left(\Gamma\right)$
by \eqref{eq:derivative-envelope},
\begin{equation}
s^{\left(G+1\right)/2}_{i}L^{\left(G+1\right)}_{i}\left(\operatorname{Reg}\left(s^{-1}_{i}\operatorname{Cov}_{\mathbb{P}^{R,-}_{i}}\left(\widehat{Y}_{i}\right)\right)\right)=L^{\left(G+1\right)}_{i}\left(\operatorname{Reg}\left(\operatorname{Cov}_{\mathbb{P}^{R,-}_{i}}\left(\widehat{Y}_{i}\right)\right)\right).\label{eq:state-covariance-error-factor}
\end{equation}
Substitution into \eqref{eq:joint-wms-budget} gives the last assertion.
\end{proof}

The test \eqref{eq:state-covariance-cancellation} can therefore accept
a patch $U_{i,\kappa}$ for which $\left\langle \mathbb{D}^{R}_{i,\kappa},h^{\alpha}_{i}\right\rangle $
does not tend to zero as $s_{i}\to0$. Let $N=1$, so that $M_{G}=G+1$,
and let the patch consist of the $G+2$ states $x+j\delta$, $j=0,\ldots,G+1$,
with weights $w2^{-G-1}\binom{G+1}{j}$, where $w>0$ is its mass.
Since the patch has more than $M_{G}$ nodes, the population rule does
not stop it. The moment equations of \eqref{eq:moments} for this patch
have a one-dimensional kernel, spanned by the alternating weights
$\left(-1\right)^{j}\binom{G+1}{j}$, so that its feasible set is
a segment. The two endpoints double the even weights and delete the
odd ones, or conversely, and every extreme-point reduction, in particular
the output of Proposition~\ref{prop:solver}, is one of them. If $h$
has continuous derivatives up to order $G+1$ on an open interval containing
$\left[x,x+\left(G+1\right)\delta\right]$, both endpoints change $\left\langle \cdot,h\right\rangle $
in absolute value by $w2^{-G-1}$ times the $\left(G+1\right)$st difference
of $h$ with step $\delta$, that is, by $w2^{-G-1}\delta^{G+1}\left|h^{\left(G+1\right)}\left(\xi\right)\right|$
for some $\xi\in\left(x,x+\left(G+1\right)\delta\right)$. Let the
remaining nodes of date $i$ carry a fixed measure of positive variance,
so that the state covariance does not shrink with $\delta$, and fix
a sufficiently small $\delta>0$. Then \eqref{eq:state-covariance-cancellation}
accepts the patch at every step size. Since $h^{\alpha}_{i}$ in \eqref{eq:tasks}
depends on $t_{i}$ but not on the grid, $\left|\left\langle \mathbb{D}^{R}_{i,\kappa},h^{\alpha}_{i}\right\rangle \right|$
stays bounded away from zero as $s_{i}\to0$ whenever $D^{G+1}h^{\alpha}_{i}$
does not vanish on $\left[x,x+\left(G+1\right)\delta\right]$.
\begin{prop}
\label{prop:wms-refinement}Let $L=1$ and $\Gamma_{i}=\Gamma^{\mathrm{AWLL}}_{i}$
in \eqref{eq:partition-rule}, fix $\lambda>0$ and $G\ge1$, and consider
a family of grids with $s_{\max}=\max_{i}s_{i}\to0$.
Suppose that a constant $K_{0}>0$, independent of the grid, satisfies,
at every reduction date,
\begin{equation}
\begin{aligned}\sum_{e}\varpi_{e}\frac{\left\lVert \Delta_{e}-\Delta_{i}\right\rVert ^{2}_{2}}{s_{i}} & \le K_{0},\\
\mathbb{P}^{R,+}_{i-1}\left(\left\lVert \Sigma\right\rVert ^{2}_{F}+s_{i}\left\lVert b\right\rVert ^{2}_{2}\right) & \le K_{0}.
\end{aligned}
\label{eq:wms-energy-assumptions}
\end{equation}
Put $K_{\Gamma}=\max\left\{ 1,\varphi K_{0}\right\} $ with the golden ratio $\varphi$.
Then, for every $i\in\mathcal{R}$,
\begin{equation}
\begin{aligned}\lambda_{\max}\left(\Gamma^{\mathrm{AWLL}}_{i}\right) & \le K_{\Gamma},\\
L^{\left(G+1\right)}_{i}\left(\Gamma^{\mathrm{AWLL}}_{i}\right) & \le K^{\left(G+1\right)/2}_{\Gamma}L^{\left(G+1\right)}_{i}\left(I_{N}\right).
\end{aligned}
\label{eq:wms-refinement-bounds}
\end{equation}
Every patch $U_{i,\kappa}$ that passes \eqref{eq:partition-rule}
satisfies $\operatorname{rad}_{I_{N}}\left(U_{i,\kappa}\right)\le\lambda K^{1/2}_{\Gamma}\left(T-t_{i}\right)^{G/\left(2\left(G+1\right)\right)}s^{1/2}_{i}w^{-1/\left(G+1\right)}_{i,\kappa}$.
Suppose, moreover, that, for constants $C_{T}>0$ and $p_{K}\ge0$
independent of the grid,
\begin{equation}
\sum_{i\in\mathcal{R}}K_{i}\left(T-t_{i}\right)^{G/2}s_{i}L^{\left(G+1\right)}_{i}\left(I_{N}\right)\le C_{T}s^{-p_{K}}_{\max}.\label{eq:wms-weighted-regularity}
\end{equation}
Then
the recombination contribution in \eqref{eq:joint-wms-budget} is at
most
\begin{equation}
\frac{2\lambda^{G+1}K^{\left(G+1\right)/2}_{\Gamma}C_{T}}{\left(G+1\right)!}\,s^{\left(G-1\right)/2-p_{K}}_{\max}.\label{eq:wms-refinement-error}
\end{equation}
\end{prop}
\begin{proof}
The two left sides of \eqref{eq:wms-energy-assumptions} are the traces
of $\widehat{\Gamma}^{\mathrm{edge}}_{i}$ and $\Gamma^{\mathrm{mod}}_{i}$.
Both matrices are positive semidefinite, so the mixture $\overline{\Gamma}_{i}$,
a convex combination of them, has $0\le\gamma^{-}_{i}\le\gamma^{+}_{i}\le K_{0}$.
If the mixture is nonzero, \eqref{eq:closed-coefficients} yields
\begin{equation}
\begin{aligned}\lambda_{\max}\left(\Gamma^{\mathrm{AWLL}}_{i}\right) & =\frac{\gamma^{+}_{i}-\gamma^{-}_{i}+\sqrt{\left(\gamma^{+}_{i}-\gamma^{-}_{i}\right)^{2}+4\left(\gamma^{+}_{i}\right)^{2}}}{2}\\
 & \le\varphi\gamma^{+}_{i}\le\varphi K_{0}.
\end{aligned}
\label{eq:wms-ridge-upper-proof}
\end{equation}
If it is zero, $\Gamma^{\mathrm{AWLL}}_{i}=I_{N}$ by Theorem~\ref{thm:spectral-scale},
and $1\le K_{\Gamma}$. Since $\lambda_{\max}\left(\Gamma^{\mathrm{AWLL}}_{i}\right)\le K_{\Gamma}$,
\eqref{eq:metric-sandwich} and \eqref{eq:partition-rule} give the
bound on $\operatorname{rad}_{I_{N}}\left(U_{i,\kappa}\right)$, and
\eqref{eq:envelope-sandwich} gives the second inequality in \eqref{eq:wms-refinement-bounds}.
Apply these bounds to \eqref{eq:joint-wms-budget}. Since $s^{\left(G+1\right)/2}_{i}\le s^{\left(G-1\right)/2}_{\max}s_{i}$,
the recombination contribution is at most
\begin{equation}
\frac{2\lambda^{G+1}K^{\left(G+1\right)/2}_{\Gamma}}{\left(G+1\right)!}s^{\left(G-1\right)/2}_{\max}\sum_{i\in\mathcal{R}}K_{i}\left(T-t_{i}\right)^{G/2}s_{i}L^{\left(G+1\right)}_{i}\left(I_{N}\right).\label{eq:wms-error-sum-proof}
\end{equation}
Equation \eqref{eq:wms-weighted-regularity} gives \eqref{eq:wms-refinement-error}.
\end{proof}

For $G=5$, the bound \eqref{eq:wms-refinement-error} is of order
$s^{2-p_{K}}_{\max}$. As in \eqref{eq:ver1-044}, $p_{K}$ describes
the combined growth of the patch counts and of the continuation derivatives.
It can be taken to be zero if both are bounded, whereas the counts
in Section~\ref{sec:numerics} grow with $n$. For the state covariance,
\eqref{eq:state-covariance-error-factor} cancels the factor $s^{\left(G+1\right)/2}_{i}$
of \eqref{eq:joint-wms-budget}. The proof of Proposition~\ref{prop:wms-refinement}
uses the mixture only through \eqref{eq:wms-energy-assumptions} and
convexity, and the ridge only through the upper bound in \eqref{eq:scale-transmission}.
Hence \eqref{eq:wms-refinement-error} holds for every mixing coefficient
in $\left[0,1\right]$, and its factor $s^{\left(G-1\right)/2}_{\max}$
comes from both references being second moments of increments per
unit time, not from the shrinkage rule \eqref{eq:mixing-coefficient}.

\subsection{Positive reduction in the original row space}
\label{sec:solver}
\noindent
Algorithm~\ref{alg:bbe} in Appendix~\ref{app:algorithm} gives the
complete procedure of this section, including its numerical acceptance
tests.

\subsubsection{Moment feasibility and numerical coordinates}

Fix $i\in\mathcal{R}$ and a patch index $\kappa$, number the nodes
of the patch as $1,\ldots,J$, and let $y_{j}\in\mathcal{O}$ and $p^{R,\pm}_{j}$
be the state of node $j$ and its weight under $\mathbb{P}^{R,\pm}_{i,\kappa}$.
The original monomial matrix of the patch is
\begin{equation}
A^{\mathrm{org}}_{\alpha j}=y^{\alpha}_{j}=\prod^{N}_{h=1}y^{\alpha_{h}}_{j,h},\qquad\alpha\in\mathbb{N}^{N}_{0},\quad\left|\alpha\right|\le G,\quad1\le j\le J.\label{eq:original-matrix}
\end{equation}
It has $M_{G}$ rows, including the constant row $\alpha=0$. Write
$k=\operatorname{rank}\left(A^{\mathrm{org}}\right)$ and $w=w_{i,\kappa}$.
Since the monomials span $\mathcal{P}_{G}\left(\mathbb{R}^{N}\right)$,
condition \eqref{eq:moments} is equivalent to
\begin{equation}
A^{\mathrm{org}}p^{R,+}=A^{\mathrm{org}}p^{R,-},\qquad p^{R,+}\ge0,\qquad\prescript{t}{}{\mathbf{1}_{J}}p^{R,+}=\prescript{t}{}{\mathbf{1}_{J}}p^{R,-}=w.\label{eq:feasible-reduction}
\end{equation}
The set of solutions of \eqref{eq:feasible-reduction} contains $p^{R,-}$
and is compact, since $0\le p^{R,+}_{j}\le w$. At an extreme point
$p^{R,+}$ of this set, the columns of $A^{\mathrm{org}}$ at its positive
entries are linearly independent. Otherwise a nonzero $v\in\ker A^{\mathrm{org}}$
supported on these entries gives the solutions $p^{R,+}\pm tv$ for
small $t>0$, whose midpoint is $p^{R,+}$. Hence an extreme point has
at most $k$ positive entries. As in the cubature reduction of \citet{tchernychova_caratheodory_2016},
\eqref{eq:feasible-reduction} is a feasibility linear program, and
the procedure below finds an extreme point.

Different monomial degrees and coordinate units can produce widely
separated row and column magnitudes. Positive diagonal equilibration
balances these magnitudes before numerical rank assessment \citep[Section~2]{takapoui_javadi_2016}.
With positive diagonal matrices $D_{\mathrm{row}}$ and $D_{\mathrm{col}}$,
compute
\begin{equation}
A^{\mathrm{eq}}=D_{\mathrm{row}}A^{\mathrm{org}}D_{\mathrm{col}},\qquad k_{\tau}=\operatorname{rank}_{\tau}\left(A^{\mathrm{eq}}\right).\label{eq:rank-scaling}
\end{equation}
Here $\operatorname{rank}_{\tau}$ counts the singular values above
a threshold $\tau$ (Appendix~\ref{app:algorithm}). Sinkhorn--Knopp
iterations on the nonzero rows of the matrix with entries $\left(A^{\mathrm{org}}_{\alpha j}/a_{\max}\right)^{2}$,
$a_{\max}=\max_{\alpha,j}\left|A^{\mathrm{org}}_{\alpha j}\right|$,
give the squares of the diagonal entries of $D_{\mathrm{col}}$ and
of the corresponding entries of $D_{\mathrm{row}}$. The nonzero rows
of $A^{\mathrm{eq}}$ thus have nearly equal Euclidean norms, and so
do its columns, while $A^{\mathrm{eq}}$ has the signs and the rank
of $A^{\mathrm{org}}$. The algebraic results below assume $k_{\tau}=k$.
The rank decision is validated afterwards: the computed factors of
\eqref{eq:rowspace} must pass orthogonality and reconstruction tests,
and the output must pass the moment residual tests below.

Modified Gram--Schmidt, applied to the rows of $A^{\mathrm{org}}$ beginning
with the constant row, gives a matrix $U\in\mathbb{R}^{k\times J}$
whose rows form an orthonormal basis of the row space of $A^{\mathrm{org}}$
and whose first row is $U_{1\cdot}=J^{-1/2}\prescript{t}{}{\mathbf{1}_{J}}$.
Since $\prescript{t}{}{U}U$ is the orthogonal projection onto this
row space, $R=A^{\mathrm{org}}\prescript{t}{}{U}\in\mathbb{R}^{M_{G}\times k}$
satisfies
\begin{equation}
A^{\mathrm{org}}=RU,\qquad U\prescript{t}{}{U}=I_{k},\qquad\ker U=\ker A^{\mathrm{org}}.\label{eq:rowspace}
\end{equation}
Indeed, $Uv=0$ implies $A^{\mathrm{org}}v=RUv=0$, and $A^{\mathrm{org}}v=0$
implies $Uv=0$, since every row of $U$ is a linear combination of
rows of $A^{\mathrm{org}}$. Elimination in $U$ therefore enforces the
moment equations of \eqref{eq:feasible-reduction}, while the columns
of $U$ still correspond to the nodes and the weights stay in the original
coordinates.

The acceptance test therefore uses the original residual. For computed
factors $R$ and $U$ and $r^{U}=U\left(p^{R,+}-p^{R,-}\right)$, the
original residual satisfies the identity
\begin{equation}
\begin{aligned}d^{\mathrm{org}} & =A^{\mathrm{org}}\left(p^{R,+}-p^{R,-}\right)\\
 & =Rr^{U}+\left(A^{\mathrm{org}}-RU\right)\left(p^{R,+}-p^{R,-}\right).
\end{aligned}
\label{eq:solver-residual-lift}
\end{equation}
The reconstruction term vanishes under \eqref{eq:rowspace}, but even
then the rows of $R$ scale the components of $r^{U}$ differently.
The algorithm therefore computes $d^{\mathrm{org}}$ from $A^{\mathrm{org}}$,
which also covers the reconstruction term, compares each $\left|d^{\mathrm{org}}_{\alpha}\right|$
with a tolerance $\tau^{\mathrm{org}}_{\alpha}$, and separately checks
the total mass and the nonnegativity of the weights.

\subsubsection{Prefix absorption and basis exchange}

The procedure starts from a tuple $B=\left(\beta_{1},\ldots,\beta_{k}\right)$
of column indices with $U_{\cdot B}$ nonsingular, the weights $\eta=p^{R,-}_{B}$
and $\tilde{p}_{j}=p^{R,-}_{j}$ for $j\notin B$, and the stream $\mathcal{J}_{\mathrm{rem}}$
of these $j$ in a fixed order. For a bound $b_{\max}\ge1$, the next
panel $J_{b}=\left(j_{1},\ldots,j_{m_{b}}\right)$
consists of the first $m_{b}\le b_{\max}$ entries of the stream, and
the coefficients of its columns in the basis are
\begin{equation}
H_{B}=U^{-1}_{\cdot B},\qquad C_{b}=H_{B}U_{\cdot J_{b}},\qquad c_{j}=H_{B}U_{\cdot j}\quad\left(j\in J_{b}\right).\label{eq:panel-coefficients}
\end{equation}
Since $U_{\cdot B}c_{j}=U_{\cdot j}$, adding $\tilde{p}_{j}c_{j}$ to $\eta$ and
setting $\tilde{p}_{j}=0$ leaves $U_{\cdot B}\eta+\sum_{j'\in\mathcal{J}_{\mathrm{rem}}}\tilde{p}_{j'}U_{\cdot j'}$
unchanged. This step is called the absorption of $j$. Consecutive nodes
of a panel are absorbed together as long as $\eta$ stays nonnegative.
Define
\begin{equation}
\eta^{\left(0\right)}=\eta,\qquad\eta^{\left(h\right)}=\eta+\sum^{h}_{q=1}\tilde{p}_{j_{q}}c_{j_{q}},\qquad1\le h\le m_{b},\label{eq:prefix}
\end{equation}
and let $h^{*}_{b}$ be the length of the longest nonnegative prefix,
\begin{equation}
h^{*}_{b}=\max\left\{ h\in\left\{ 0,\ldots,m_{b}\right\} \middle|\eta^{\left(q\right)}\ge0\text{ for every }1\le q\le h\right\} .\label{eq:prefix-length}
\end{equation}
Then $\eta$ is replaced by $\eta^{\left(h^{*}_{b}\right)}$, and the
absorbed nodes are removed from the stream with their weights set
to zero. If $h^{*}_{b}=m_{b}$, the whole panel has been consumed.

Otherwise the next column $j=j_{h^{*}_{b}+1}$ has $\eta+\tilde{p}_{j}c_{j}\not\ge0$,
and only the largest amount that keeps $\eta$ nonnegative is transferred.
The ratio test selects a leaving position $a^{*}$ and the transfer
amount $\delta_{j}$,
\begin{equation}
a^{*}\in\operatorname*{arg\,min}_{a:\left(c_{j}\right)_{a}<0}\frac{\eta_{a}}{-\left(c_{j}\right)_{a}},\qquad\delta_{j}=\frac{\eta_{a^{*}}}{-\left(c_{j}\right)_{a^{*}}}.\label{eq:basis-exchange}
\end{equation}
At this amount the weight in position $a^{*}$ reaches zero. The exchange
replaces $\beta_{a^{*}}$ by $j$ in the same position and gives $j$ the
untransferred weight,
\begin{equation}
B'=\left(\beta_{1},\ldots,\beta_{a^{*}-1},j,\beta_{a^{*}+1},\ldots,\beta_{k}\right),\qquad\eta'=\eta+\delta_{j}c_{j}+\left(\tilde{p}_{j}-\delta_{j}\right)e_{a^{*}}.\label{eq:ordered-basis-replacement}
\end{equation}
Then $j$ leaves the stream with $\tilde{p}_{j}=0$, the remaining columns
of the panel are expressed in the basis $B'$ by \eqref{eq:panel-suffix-update}
below, and the procedure, called batched basis exchange, continues
with $B'$ and $\eta'$. When the stream is empty, its output $p^{R,+}$
has $p^{R,+}_{B}=\eta$ and $p^{R,+}_{j}=0$ for $j\notin B$.
\begin{prop}
\label{prop:solver} Assume \eqref{eq:rowspace}. Then every absorption
and every exchange keeps $\eta\ge0$ and $U_{\cdot B}$ nonsingular, and preserves
\begin{equation}
U_{\cdot B}\eta+\sum_{j\in\mathcal{J}_{\mathrm{rem}}}\tilde{p}_{j}U_{\cdot j}=Up^{R,-}.\label{eq:solver-invariant}
\end{equation}
The procedure ends after finitely many steps, and its output $p^{R,+}$
is an extreme point of the set of solutions of \eqref{eq:feasible-reduction}
with at most $k$ positive entries.
\end{prop}

\begin{proof}
The initial weights and stream satisfy \eqref{eq:solver-invariant}.
For a complete absorption, 
\begin{equation}
U_{\cdot B}\left(\eta+\tilde{p}_{j}c_{j}\right)=U_{\cdot B}\eta+\tilde{p}_{j}U_{\cdot j},\label{eq:ver1-140}
\end{equation}
 so consuming node $j$ preserves the invariant. The prefix inequalities
guarantee nonnegativity at every such absorption. At a failed full
absorption, $\eta+\tilde{p}_{j}c_{j}$ has a negative component. Since
$\eta\ge0$, some $\left(c_{j}\right)_{a}<0$ and $\eta_{a}/\left(-\left(c_{j}\right)_{a}\right)<\tilde{p}_{j}$.
Hence the selected ratio satisfies $0\le\delta_{j}<\tilde{p}_{j}$, and
\begin{equation}
\eta_{a}+\delta_{j}\left(c_{j}\right)_{a}\ge0\quad\left(1\le a\le k\right),\qquad\eta_{a^{*}}+\delta_{j}\left(c_{j}\right)_{a^{*}}=0.\label{eq:ver1-141}
\end{equation}
 Let $B'$ replace $\beta_{a^{*}}$ by $j$. Multilinearity of the determinant
gives 
\begin{equation}
\det U_{\cdot B'}=\left(c_{j}\right)_{a^{*}}\det U_{\cdot B}\ne0.\label{eq:ver1-142}
\end{equation}
 Because the old weight at position $a^{*}$ is zero after the ratio transfer,
\begin{equation}
\begin{aligned}U_{\cdot B'}\left(\eta+\delta_{j}c_{j}+\left(\tilde{p}_{j}-\delta_{j}\right)e_{a^{*}}\right) & =U_{\cdot B}\left(\eta+\delta_{j}c_{j}\right)+\left(\tilde{p}_{j}-\delta_{j}\right)U_{\cdot j}\\
 & =U_{\cdot B}\eta+\tilde{p}_{j}U_{\cdot j}.
\end{aligned}
\label{eq:ver1-143}
\end{equation}
 Thus the new nonnegative weights preserve \eqref{eq:solver-invariant}.
Every absorption and every exchange removes a node from the stream,
and an exchange follows every prefix that stops before the end of
its panel. Hence
the procedure ends, and \eqref{eq:solver-invariant} then gives $Up^{R,+}=U_{\cdot B}\eta=Up^{R,-}$.
Multiplication by $R$ in \eqref{eq:rowspace} gives the moment equations
of \eqref{eq:feasible-reduction}, the constant row giving the mass.
If $p^{R,+}$ is the midpoint of two solutions, both vanish off $B$,
so their difference is a kernel vector supported on independent columns
and vanishes. The output is therefore an extreme point, and it has
at most $k$ positive entries.
\end{proof}

\subsubsection{Basis reuse and elimination work}

An exchange changes one column of $U_{\cdot B}$. With $c=c_{j}$ and the
leaving position $a^{*}$, write
\begin{equation}
E=I_{k}+\left(c-e_{a^{*}}\right)\prescript{t}{}{e_{a^{*}}},\qquad U_{\cdot B'}=U_{\cdot B}E.\label{eq:basis-eta-update}
\end{equation}
The ratio test gives $c_{a^{*}}<0$, and direct multiplication verifies
\begin{equation}
E^{-1}=I_{k}-\frac{\left(c-e_{a^{*}}\right)\prescript{t}{}{e_{a^{*}}}}{c_{a^{*}}},\qquad H_{B'}=E^{-1}H_{B}.\label{eq:basis-inverse-update}
\end{equation}
Let $C_{\mathrm{tail}}$ consist of the columns of $C_{b}$ for the
nodes of the panel that remain in the stream. In the new basis they
become
\begin{equation}
C'_{\mathrm{tail}}=E^{-1}C_{\mathrm{tail}}=C_{\mathrm{tail}}-\frac{\left(c-e_{a^{*}}\right)\left(\prescript{t}{}{e_{a^{*}}}C_{\mathrm{tail}}\right)}{c_{a^{*}}}.\label{eq:panel-suffix-update}
\end{equation}
This tail update costs $O\left(km_{\mathrm{tail}}\right)$ operations
for $m_{\mathrm{tail}}$ columns and permits further absorptions and
exchanges within the panel without updating $H_{B}$. Algorithm~\ref{alg:bbe}
therefore recomputes $H_{B}$ by a dense LU factorization of $U_{\cdot B}$
only once, at the end of a panel in which an exchange occurred.

The elementary matrix in \eqref{eq:basis-eta-update} is the product-form
basis update reviewed by \citet[Section~2.4]{huangfu_hall_2018}.
\citet{forrest_tomlin_1972} update triangular basis factors to preserve
sparsity. Huangfu and Hall's parallel dual revised simplex method
groups minor iterations and collective updates \citep[Sections~2.3 and~3]{huangfu_hall_2018}.
Here the objective is constant, so the columns enter in the order of
the stream. The coefficients of a panel are computed by one matrix
product, and different patches are reduced independently.
\begin{prop}
\label{prop:kernel-work}Let $r=J-k\ge1$, let $N_{\mathrm{ex}}$ be
the number of exchanges, and let $Q_{\mathrm{ref}}$ be the number
of panels in which an exchange occurs. Then $N_{\mathrm{ex}}\le r$
and $Q_{\mathrm{ref}}\le N_{\mathrm{ex}}$, and, excluding the computation
of $k$, $U$ and the initial basis, batched basis exchange uses
\begin{equation}
O\left(k^{3}\left(1+Q_{\mathrm{ref}}\right)+k^{2}r+kb_{\max}N_{\mathrm{ex}}\right)\label{eq:panel-work-bound}
\end{equation}
arithmetic operations, which is $O\left(J\right)$ for fixed $k$
and $b_{\max}$, and $O\left(k^{2}+kb_{\max}\right)$ entries of auxiliary
storage besides $U$ and the weights. The direct kernel update of
\citet[Appendix~A, Eqs.~(A.1) and (A.2)]{ninomiya_high-order_2025}
uses $\Theta\left(kr^{2}+r^{3}\right)$ operations, which is $\Theta\left(J^{3}\right)$
for fixed $k$, excluding its initial singular value decomposition,
and stores $\Theta\left(Jr\right)$ kernel entries.
\end{prop}
\begin{proof}
Let $N_{\mathrm{ex},b}$ be the number of exchanges in panel $b$,
and let $m_{b,h}$ be the number of tail columns updated at exchange
$h$ of panel $b$. Since each node of the stream lies in at most one
panel,
\begin{equation}
\begin{aligned}\sum_{b}m_{b} & \le r,\qquad0\le m_{b,h}\le b_{\max},\\
N_{\mathrm{ex}} & =\sum_{b}N_{\mathrm{ex},b}\le r,\qquad Q_{\mathrm{ref}}=\sum_{b}\mathbf{1}_{\left\{ N_{\mathrm{ex},b}>0\right\} }\le N_{\mathrm{ex}}.
\end{aligned}
\label{eq:panel-work-counts}
\end{equation}
The arithmetic account is
\begin{equation}
\begin{array}{ll}
\text{Panel coefficients }H_{B}U_{\cdot J_{b}} & {\displaystyle \sum_{b}O\left(k^{2}m_{b}\right)=O\left(k^{2}r\right),}\\[3pt]
\text{Positivity, prefixes and ratio tests} & O\left(kr+kb_{\max}N_{\mathrm{ex}}\right),\\[3pt]
\text{Tail updates} & {\displaystyle \sum_{b,h}O\left(km_{b,h}\right)=O\left(kb_{\max}N_{\mathrm{ex}}\right),}\\[3pt]
\text{Initial and collective LU factorizations} & O\left(k^{3}\left(1+Q_{\mathrm{ref}}\right)\right).
\end{array}\label{eq:panel-work-account}
\end{equation}
Their sum gives \eqref{eq:panel-work-bound}. The stored basis, inverse
and panel contain $O\left(k^{2}+kb_{\max}\right)$ entries. For the
direct kernel update,
\begin{equation}
\begin{aligned}\sum^{r}_{i=1}\left(J-i+1\right)\left(r-i\right) & =\sum^{r-1}_{t=0}\left(k+t+1\right)t\\
 & =\frac{\left(k+1\right)r\left(r-1\right)}{2}+\frac{r\left(r-1\right)\left(2r-1\right)}{6},\\
\sum^{r}_{i=1}\left(J-i+1\right) & =kr+\frac{r\left(r+1\right)}{2}.
\end{aligned}
\label{eq:direct-kernel-work}
\end{equation}
The vector updates and weight scans thus take $\Theta\left(kr^{2}+r^{3}\right)$
operations, and $r=J-k$ gives $\Theta\left(J^{3}\right)$ for fixed
$k$.
\end{proof}

The divide-and-conquer reduction of \citet[Corollary~4]{litterer_high_2012}
also uses work linear in $J$ for a fixed number of moments.

\subsubsection{Numerical moment error}

Let the weights of $\mathbb{P}^{R,+}_{i,\kappa}$ be those computed
by Algorithm~\ref{alg:bbe}. Let $h$ have continuous derivatives up
to order $G+1$ on an open subset of $\mathcal{O}$ that contains $\mathrm{C.H.}(\operatorname{supp}(\mathbb{P}^{R,-}_{i}))$,
and let $\operatorname{Tay}_{G}h$ be its Taylor polynomial of degree
$G$ about the center $c_{i,\kappa}$ of \eqref{eq:residual-tensor}.
Then $e_{\mathrm{mom}}\left(h\right)=\left\langle \mathbb{D}^{R}_{i,\kappa},\operatorname{Tay}_{G}h\right\rangle $
is a linear combination of the components of $d^{\mathrm{org}}$, and
\begin{equation}
\left\langle \mathbb{D}^{R}_{i,\kappa},h\right\rangle =e_{\mathrm{mom}}\left(h\right)+\left\langle \mathbb{D}^{R}_{i,\kappa},h-\operatorname{Tay}_{G}h\right\rangle .\label{eq:numerical-moment-split}
\end{equation}
Let $r^{{M}}_{i,\kappa}$ be the centered metric radius of
\eqref{eq:residual-tensor}
and let $x$ range over $\mathrm{C.H.}(\operatorname{supp}(\mathbb{P}^{R,-}_{i}))$.
Since both weight vectors are nonnegative, $\left\lVert \mathbb{D}^{R}_{i,\kappa}\right\rVert _{\mathrm{TV}}\le w_{i,\kappa}+\mathbb{P}^{R,+}_{i,\kappa}\left(1\right)$,
and the argument for \eqref{eq:ver1-039} gives
\begin{equation}
\begin{aligned}\left|\left\langle \mathbb{D}^{R}_{i,\kappa},h\right\rangle \right| & \le\left|e_{\mathrm{mom}}\left(h\right)\right|\\
 & \quad+\frac{\left(w_{i,\kappa}+\mathbb{P}^{R,+}_{i,\kappa}\left(1\right)\right)\left(r^{{M}}_{i,\kappa}\right)^{G+1}}{\left(G+1\right)!}\sup_{x,\,\left\lVert v\right\rVert _{2}=1}\left|D^{G+1}h\left(x\right)\left[M^{-1}_{i}v,\ldots,M^{-1}_{i}v\right]\right|.
\end{aligned}
\label{eq:numerical-moment-remainder}
\end{equation}
The tolerances of the acceptance test thus control $e_{\mathrm{mom}}\left(h\right)$,
and the centered metric radius controls the remainder of order $G+1$.
With $h=h^{\alpha}_{i}$, \eqref{eq:numerical-moment-remainder} takes
the place of \eqref{eq:ver1-039} in the proof of Theorem~\ref{thm:residual}
for the computed measures.

\FloatBarrier

\section{Patch counts and retained populations}
\label{sec:work}
\noindent
This section bounds $K_{i}$ and $J^{+}_{i}$.
\begin{thm}
\label{thm:resolution-work}Let $i\in\mathcal{R}$. For $\Gamma\succ0$
and $\lambda>0$, a partition of the support of $\mathbb{P}^{R,-}_{i}$
is called admissible for $(\Gamma,\lambda)$ if each of its patches
passes \eqref{eq:partition-rule} with $\Gamma$ in place of $\Gamma_{i}$,
and $K^{\min}_{i}(\Gamma,\lambda)$ denotes the smallest number of patches
of such a partition. Then
\begin{equation}
K^{\min}_{i}\left(I_{N},\sqrt{\lambda_{\max}\left(\Gamma\right)}\right)\le K^{\min}_{i}\left(\Gamma,1\right)\le K^{\min}_{i}\left(I_{N},\sqrt{\lambda_{\min}\left(\Gamma\right)}\right).\label{eq:partition-sandwich}
\end{equation}
This remains true when the right side of \eqref{eq:partition-rule}
is replaced by $\lambda a(w_{i,\kappa})$ for a positive function $a$.
If $L=1$ and the mixed block $\overline{\Gamma}_{i,1}$ is nonzero,
then, with $\gamma^{+}_{i}=\gamma^{+}_{i,1}$ and the golden ratio $\varphi$,
\begin{equation}
K^{\min}_{i}\left(I_{N},\sqrt{\varphi\gamma^{+}_{i}}\right)\le K^{\min}_{i}\left(\Gamma^{\mathrm{AWLL}}_{i},1\right)\le K^{\min}_{i}\left(I_{N},\sqrt{\gamma^{+}_{i}/\varphi}\right).\label{eq:forward-work-scale}
\end{equation}
Moreover, the bound $2L^{(G+1)}_{i}(\Gamma)(T-t_{i})^{G/2}s^{(G+1)/2}_{i}/(G+1)!$
on the certificate of a patch that passes \eqref{eq:partition-rule}
with $\Gamma$ in place of $\Gamma_{i}$ and $\lambda=1$ lies between
the corresponding bounds for the WLL test with the scalars $\lambda_{\min}(\Gamma)^{1/2}$
and $\lambda_{\max}(\Gamma)^{1/2}$.
\end{thm}
\begin{proof}
By item~2 of Lemma~\ref{lem:metric-scalar}, applied pair by pair
with a common right side, every partition admissible for $(I_{N},\lambda_{\min}(\Gamma)^{1/2})$
is admissible for $(\Gamma,1)$, and every partition admissible for
$(\Gamma,1)$ is admissible for $(I_{N},\lambda_{\max}(\Gamma)^{1/2})$.
The partition into single nodes is admissible for every pair, since
its patches have pairwise radius zero. Taking minima over these nested
nonempty sets of partitions proves \eqref{eq:partition-sandwich},
also for a general $a$. The same nesting shows that $K^{\min}_{i}(I_{N},\lambda)$
is nonincreasing in $\lambda$. By \eqref{eq:scale-transmission},
$\lambda_{\min}(\Gamma^{\mathrm{AWLL}}_{i})\ge\gamma^{+}_{i}/\varphi$
and $\lambda_{\max}(\Gamma^{\mathrm{AWLL}}_{i})\le\varphi\gamma^{+}_{i}$,
so \eqref{eq:partition-sandwich} gives \eqref{eq:forward-work-scale}.
The last assertion follows from \eqref{eq:envelope-sandwich}.
\end{proof}

Theorem~\ref{thm:resolution-work} concerns admissible partitions of
one frontier. The recursive subdivision of Section~\ref{sec:subdivision}
also stops, by the population rule, at patches with at most $M_{G}$
nodes that fail \eqref{eq:partition-rule}, and Theorem~\ref{thm:finite-scale-work}
and Corollary~\ref{cor:refinement-exponent} include this rule. Runs
with different matrices select different patches, and their frontiers
at later dates differ. For the patches actually selected, the positive
reduction gives the following bounds.
\begin{thm}
\label{thm:population}Let $i\in\mathcal{R}$, and let $J_{i,\kappa}$
be the number of nodes of $U_{i,\kappa}$ and $k_{i,\kappa}$ the rank
of its matrix \eqref{eq:original-matrix}. If each $\mathbb{P}^{R,+}_{i,\kappa}$
is an extreme point of the set of solutions of \eqref{eq:feasible-reduction},
then
\begin{equation}
\begin{aligned}J^{+}_{i} & \le\sum^{K_{i}}_{\kappa=1}k_{i,\kappa}\le\sum^{K_{i}}_{\kappa=1}\min\left\{ J_{i,\kappa},M_{G}\right\} ,\\
J^{-}_{i}-J^{+}_{i} & \ge\sum^{K_{i}}_{\kappa=1}\left(J_{i,\kappa}-M_{G}\right)_{+}.
\end{aligned}
\label{eq:population-rank}
\end{equation}
Consequently, the ratio $\upsilon_{i}$ below, called the nominal polynomial
capacity utilization, satisfies
\begin{equation}
\upsilon_{i}=\frac{J^{+}_{i}}{M_{G}K_{i}},\qquad0<\upsilon_{i}\le1.\label{eq:capacity-utilization}
\end{equation}
If, moreover, $n\ge2$, $\mathcal{R}=\{1,\ldots,n-1\}$, this condition
holds at every $i\in\mathcal{R}$, and every node $[j]$ in the support
of $\mathbb{P}^{R,+}_{i-1}$, $1\le i\le n$, has exactly $b_{\mathrm{br}}\ge2$
labels $m\in\mathcal{B}_{i}$ with $\pi_{i,m}(\widehat{Y}_{i-1}(j))>0$,
then
\begin{equation}
\begin{aligned}J^{-}_{n} & =b_{\mathrm{br}}M_{G}K_{n-1}\upsilon_{n-1},\\
J^{-}_{n} & \le\min\left\{ b^{n}_{\mathrm{br}},\min_{1\le i<n}b^{n-i}_{\mathrm{br}}M_{G}K_{i}\right\} .
\end{aligned}
\label{eq:propagated-capacity}
\end{equation}
\end{thm}
\begin{proof}
By the argument after \eqref{eq:feasible-reduction}, an extreme point
of the set of solutions of \eqref{eq:feasible-reduction} has at most
$k_{i,\kappa}\le\min\{J_{i,\kappa},M_{G}\}$ positive entries. Summing
over $\kappa$ proves the first line of \eqref{eq:population-rank},
and subtracting it from $J^{-}_{i}=\sum_{\kappa}J_{i,\kappa}$ gives
the second. Since $J^{+}_{i}\ge1$, the first line gives \eqref{eq:capacity-utilization}.
Branching gives $J^{-}_{i+1}=b_{\mathrm{br}}J^{+}_{i}$, while every
reduction satisfies $J^{+}_{i}\le J^{-}_{i}$. Consequently
\begin{equation}
J^{-}_{n}\le b^{n-i}_{\mathrm{br}}J^{+}_{i}\le b^{n-i}_{\mathrm{br}}M_{G}K_{i}\label{eq:ver1-108}
\end{equation}
for $1\le i<n$, and $J^{-}_{n}\le b^{n}_{\mathrm{br}}$ since $J^{+}_{0}=1$.
The definition of $\upsilon_{n-1}$ gives the equality in \eqref{eq:propagated-capacity}.
\end{proof}

The retention fraction is $\eta_{i}=J^{+}_{i}/J^{-}_{i}$ for $i\in\mathcal{R}$
and $\eta_{i}=1$ at the other interior dates. If every node $[j]$
in the support of $\mathbb{P}^{R,+}_{i-1}$ has exactly $b^{\mathrm{Alg}}_{i}$
labels $m\in\mathcal{B}_{i}$ with $\pi_{i,m}(\widehat{Y}_{i-1}(j))>0$,
$1\le i\le n$, then $J^{-}_{i}=b^{\mathrm{Alg}}_{i}J^{+}_{i-1}$, and
repeated substitution gives
\begin{equation}
J^{-}_{n}=\prod^{n}_{i=1}b^{\mathrm{Alg}}_{i}\prod^{n-1}_{i=1}\eta_{i}.\label{eq:retention-history}
\end{equation}

The ranks in \eqref{eq:population-rank} are exact. If the computed
output of every patch of date $i$ passes the acceptance tests of
Algorithm~\ref{alg:bbe}, the first line of \eqref{eq:population-rank}
holds with the numerical ranks $k_{\tau}\le\min\{J_{i,\kappa},M_{G}\}$
of \eqref{eq:rank-scaling} in place of $k_{i,\kappa}$, so \eqref{eq:capacity-utilization}
still holds, although these ranks can exceed the exact ones (Section~\ref{sec:global-lambda}).

\subsection{Recursive subdivision}\label{sec:subdivision}
\noindent
Let $i\in\mathcal{R}$ and $J=J^{-}_{i}$. For $\Gamma\succ0$, the
frontier radius $R_{i}(\Gamma)$ and the resolution load $\chi_{i}(\Gamma)$
are defined by
\begin{equation}
\begin{aligned}\rho_{i} & =\left(T-t_{i}\right)^{G/\left(2\left(G+1\right)\right)}s^{1/2}_{i},\\
R_{i}\left(\Gamma\right) & =\operatorname{rad}_{\Gamma^{-1/2}}\left(\operatorname{supp}\mathbb{P}^{R,-}_{i}\right),\\
\chi_{i}\left(\Gamma\right) & =\left(\frac{R_{i}\left(\Gamma\right)}{\rho_{i}}\right)^{G+1},
\end{aligned}
\label{eq:resolution-load}
\end{equation}
where $\rho_{i}$ is the allowance of a patch of unit mass with $\lambda=1$.

The recursive subdivision starts from the support of $\mathbb{P}^{R,-}_{i}$,
with $\Gamma_{i}$ fixed. A candidate
is a leaf when it passes \eqref{eq:partition-rule} or, by the population
rule of Section~\ref{sec:directional}, has at most $M_{G}$ nodes,
and it is split otherwise. A split orders the candidate by projection
onto the principal axis of the current subset, recomputed for every
subset in the coordinates $M_{i}y$, with $M_{i}=I_{N}$ for WLL,
and divides it at the median of this projection into two sets whose
populations differ by at most one, as in \citet[Section~\mbox{3.2}, Step~\mbox{3 (i)}]{ninomiya_high-order_2025}.
This splitter is referred to as population bisection.

The upper bounds in Theorem~\ref{thm:finite-scale-work} depend on the
frontier only through $\chi_{i}$ and $J$, and Proposition~\ref{prop:compression-regime}
applies them at every date of a run.

The initial value of an asset price is usually normalized, although
the problem does not depend on this choice. Scaling it by a factor
$c>0$ maps the states by $y'=D(c)y$ with the unit map \eqref{eq:ver1-096},
whose exponents $\zeta_{j}$ follow \eqref{eq:jet-units}, and changes
the allowances in $y$ of the WLL test (Section~\ref{sec:directional}).
The last two assertions of Theorem~\ref{thm:finite-scale-work} give
the effect of this scaling on the number of patches of a frontier
whose weights and states $\widehat{Y}_{i}$ are fixed, with
$R_{i}$, $\chi_{i}$ and $\Gamma^{\mathrm{AWLL}}_{i}(c)$ computed after
the transport of Appendix~\ref{sec:transport}.
\begin{thm}
\label{thm:finite-scale-work}Let $K_{i}(\Gamma)$ be the number of
leaves that population bisection returns with the matrix test \eqref{eq:partition-rule}
for $\Gamma_{i}=\Gamma$ and $\lambda=1$, and put $L_{J}=\lceil\log_{2}J\rceil$.
Then
\begin{equation}
K_{i}\left(\Gamma\right)\le\min\left\{ J,\ 1+\sum^{L_{J}-1}_{k=0}\min\left\{ 2^{k},\chi_{i}\left(\Gamma\right)\right\} \right\} .\label{eq:recursive-work-upper}
\end{equation}
If each input weight is at most $\omega/J$, where $\omega\ge1$, then
\begin{equation}
K_{i}\left(\Gamma\right)\le\min\left\{ J,\max\left\{ 1,4\omega\chi_{i}\left(\Gamma\right)\right\} \right\} .\label{eq:balanced-weight-upper}
\end{equation}
Let $c>0$, and let $\Pi_{s}$ be the orthogonal projection onto the
span of the coordinate vectors $e_{j}$ with $\zeta_{j}=1$, so that
$\Pi_{s}D(c)=c\Pi_{s}$. Assume that, for constants $A>0$,
$d>0$ and $0<\omega_{\mathrm{at}}<1/2$,
\begin{equation}
\sup_{z}\mathbb{P}^{R,-}_{i}\left(\left\lVert \Pi_{s}\widehat{Y}_{i}-z\right\rVert _{2}\le r\right)\le Ar^{d}+\omega_{\mathrm{at}},\qquad r>0.\label{eq:finite-concentration}
\end{equation}
Then, for $\lambda>0$, every partition of the support of $\mathbb{P}^{R,-}_{i}$
whose patches pass the WLL test \eqref{eq:raw-rule} in the coordinates
$y'=D(c)y$ or have at most $M_{G}$ nodes has at least
\begin{equation}
\min\left\{ \frac{1}{M_{G}\omega_{\mathrm{at}}},\left(2A\right)^{-\left(G+1\right)/\left(G+1+d\right)}\left(\frac{c}{\lambda\rho_{i}}\right)^{d\left(G+1\right)/\left(G+1+d\right)}\right\} \label{eq:concentration-work-lower}
\end{equation}
patches. Under the hypotheses of Lemma~\ref{lem:unit-response},
$\chi_{i}(\Gamma^{\mathrm{AWLL}}_{i}(c))$
is bounded independently of $c\ge1$.
\end{thm}
\begin{proof}
If $R_{i}(\Gamma)=0$, the root is accepted and the bound is immediate.
Otherwise, a rejected subset of mass $w$ has radius at most $R_{i}(\Gamma)$,
and its failure of \eqref{eq:partition-rule} implies
\begin{equation}
R_{i}\left(\Gamma\right)>\rho_{i}w^{-1/\left(G+1\right)},\qquad w>\chi_{i}\left(\Gamma\right)^{-1}.\label{eq:ver1-115}
\end{equation}
At depth $k$, rejected subsets are disjoint and there are at most
$2^{k}$ of them. Their masses sum to at most one, so their number
is at most $\chi_{i}(\Gamma)$. A subset at depth $k$ has at most
$\lceil J/2^{k}\rceil$ nodes, so every subset at depth $L_{J}$ is
a single node and is accepted. A finite full binary tree has one more
leaf than split nodes, so summing the bounds by depth proves \eqref{eq:recursive-work-upper},
the bound $J$ holding since the leaves are disjoint. If the weights
satisfy the stated upper bound, a subset at depth $k<L_{J}$ has mass
at most
\begin{equation}
\frac{\omega}{J}\left\lceil \frac{J}{2^{k}}\right\rceil \le2\omega\,2^{-k}.\label{eq:ver1-116}
\end{equation}
If $\chi_{i}(\Gamma)\le1$, the root passes. Otherwise, by \eqref{eq:ver1-115}
and \eqref{eq:ver1-116}, every subset at a depth $k<L_{J}$ with $2^{k}\ge2\omega\chi_{i}(\Gamma)$
passes. If such a depth exists, the smallest one has $2^{k}<4\omega\chi_{i}(\Gamma)$,
and there are at most $2^{k}$ leaves. Otherwise $J\le2^{L_{J}}<4\omega\chi_{i}(\Gamma)$.
This proves \eqref{eq:balanced-weight-upper}. Let a patch of mass
$w$ pass the WLL test in the coordinates $y'=D(c)y$, and let $x$
be the state of one of its nodes. Since $\Pi_{s}D(c)=c\Pi_{s}$
and $\Pi_{s}$ is an orthogonal projection, the projected states of
the patch lie in the ball about $\Pi_{s}x$ of radius
\begin{equation}
\frac{\lambda\rho_{i}}{c}w^{-1/\left(G+1\right)}.\label{eq:ver1-117}
\end{equation}
Hence \eqref{eq:finite-concentration} implies
\begin{equation}
w\le A\left(\frac{\lambda\rho_{i}}{c}\right)^{d}w^{-d/\left(G+1\right)}+\omega_{\mathrm{at}}.\label{eq:ver1-118}
\end{equation}
If $w>2\omega_{\mathrm{at}}$, moving the last term to the left and
using $w-\omega_{\mathrm{at}}>w/2$ yields
\begin{equation}
w\le\left(2A\right)^{\left(G+1\right)/\left(G+1+d\right)}\left(\frac{\lambda\rho_{i}}{c}\right)^{d\left(G+1\right)/\left(G+1+d\right)}.\label{eq:ver1-119}
\end{equation}
A patch with at most $M_{G}$ nodes has at most $M_{G}$ distinct projected
states, each of mass at most $\omega_{\mathrm{at}}$ by \eqref{eq:finite-concentration}
as $r\to0$. Since $M_{G}\ge2$, every patch mass is at most the larger
of the right side of \eqref{eq:ver1-119} and $M_{G}\omega_{\mathrm{at}}$,
and since the patch masses sum to one, \eqref{eq:concentration-work-lower}
follows. For AWLL, \eqref{eq:normalized-unit-bound} bounds the normalized
displacement $\Gamma^{\mathrm{AWLL}}_{i}(c)^{-1/2}D(c)(x-y)$ of any
two states $x$ and $y$ by a constant independent of $c\ge1$
times $\lVert x-y\rVert_{2}$. Hence $R_{i}(\Gamma^{\mathrm{AWLL}}_{i}(c))$,
and with it $\chi_{i}(\Gamma^{\mathrm{AWLL}}_{i}(c))$, is bounded independently
of $c\ge1$.
\end{proof}
The upper bounds hold for every choice of split directions. The term
$\omega_{\mathrm{at}}$ accounts for atoms of the projected frontier,
including nodes with equal projections, and the exponent $d$
describes the frontier and is not determined by $N$ or $d_{B}$. For
AWLL, \eqref{eq:recursive-work-upper} and \eqref{eq:balanced-weight-upper}
therefore give bounds independent of $c\ge1$, of order $1+\log J$
and, for fixed $\omega$, of order one, whereas the lower bound \eqref{eq:concentration-work-lower}
for the WLL test grows like $c^{d(G+1)/(G+1+d)}$
until its atomic saturation level $1/(M_{G}\omega_{\mathrm{at}})$.
\begin{rem}
\label{rem:harmonic-exponent}Write $\psi(d)=d(G+1)/(d+G+1)$, so
that $1/\psi(d)=1/d+1/(G+1)$. With a fixed radius allowance $\rho$,
a $d$-dimensional measure needs of order $\rho^{-d}$ patches. The
factor $w^{-1/(G+1)}_{i,\kappa}$ in \eqref{eq:partition-rule} instead
turns the test with $\lambda=1$ into $w_{i,\kappa}\operatorname{rad}_{M_{i}}(U_{i,\kappa})^{G+1}\le\rho^{G+1}$,
where $\rho$ is its right side at unit mass, so an accepted patch
has local $L^{G+1}$ distortion at most $\rho^{G+1}$ about each of
its nodes. Low-mass regions are therefore covered coarsely,
and the exponent in \eqref{eq:concentration-work-lower} is
$\psi(d)<\min\{d,G+1\}$.
The same exponent
governs optimal quantization of order $G+1$. Replace the frontier
by an absolutely continuous probability $P$ on $\mathbb{R}^{d}$
with compactly supported density $h$, and let $K^{\min}(\rho)$ be
the least number of measurable patches $U$ with $P(U)\operatorname{rad}_{I_{d}}(U)^{G+1}\le\rho^{G+1}$,
$\operatorname{rad}_{I_{d}}(U)=\sup_{x,y\in U}\lVert x-y\rVert_{2}$,
partitioning its support. Such a partition into $K$ patches yields
a $K$-point quantizer with distortion at most $K\rho^{G+1}$,
and the asymptotic formula for the optimal quantization error
\citep[Theorem~6.2]{zador_1982,bucklew_wise_1982,graf_luschgy_2000}
gives
\begin{equation}
\liminf_{\rho\to0}\rho^{\psi(d)}K^{\min}(\rho)\ge Q_{G+1}\left([0,1]^{d}\right)^{d/(d+G+1)}\int h^{d/(d+G+1)}\,dx,\label{eq:zador-count}
\end{equation}
where $Q_{G+1}([0,1]^{d})$ is the quantization coefficient of the uniform
distribution on the unit cube. The matching order $\rho^{-\psi(d)}$
follows from item 1 of Proposition~\ref{prop:minimal-partition}
below with $N=d$, $B=1$ and $d'=d$, whose proof applies verbatim to $P$
since rejected cubes occur at only finitely many dyadic levels.
\end{rem}

\begin{prop}
\label{prop:minimal-partition}Let $i\in\mathcal{R}$ and $\Gamma\succ0$,
and write $z=\Gamma^{-1/2}\widehat{Y}_{i}$,
$R=R_{i}(\Gamma)$, $\mu$ for the image of $\mathbb{P}^{R,-}_{i}$
under $z$, and $\rho=\rho_{i}$.
\begin{enumerate}
\item Let a cube of side $R$ contain the normalized frontier, and
suppose that for every $k\ge0$ at most $B2^{kd'}$ cubes of side
$2^{-k}R$ in its dyadic subdivision meet the frontier, where $B\ge1$
and $d'>0$. Then
\begin{equation}
K^{\min}_{i}(\Gamma,1)\le1+(2^{N}-1)\,c_{d',G}\,B^{(G+1)/(d'+G+1)}\left(N^{(G+1)/2}\chi_{i}(\Gamma)\right)^{d'/(d'+G+1)},\label{eq:minimal-partition-upper}
\end{equation}
with $c_{d',G}=(1-2^{-d'})^{-1}+(1-2^{-(G+1)})^{-1}$, and the right
side also bounds the number of leaves of the dyadic construction in
the proof. Every frontier satisfies the hypothesis with $B=1$
and $d'=N$.
\item Suppose
$\sup_{x}\mu(\lVert z-x\rVert_{2}\le r)\le A(r/R)^{d}+\omega_{\mathrm{at}}$
for $r>0$, with $A>0$, $d>0$ and $0<\omega_{\mathrm{at}}<1/2$. Then
every admissible partition has at least
\begin{equation}
\min\left\{ \frac{1}{2\omega_{\mathrm{at}}},(2A)^{-(G+1)/(d+G+1)}\chi_{i}(\Gamma)^{d/(d+G+1)}\right\} \label{eq:minimal-partition-lower}
\end{equation}
patches.
\end{enumerate}
When $d'=d$, the minimal count is therefore of order
$\chi_{i}(\Gamma)^{d/(d+G+1)}$ above the atomic level, whereas the
bounds of Theorem~\ref{thm:finite-scale-work} for population bisection
are linear in $\chi_{i}(\Gamma)$.
\end{prop}

\begin{proof}
For item 1, accept a cube when its nodes satisfy \eqref{eq:partition-rule}
with $\lambda=1$, and otherwise split it into its $2^{N}$ children,
assigning each shared face to one child and discarding empty ones.
Distinct normalized states are separated, so
the construction stops, and its leaves form an admissible partition.
A tree whose internal nodes have at most $2^{N}$ children has at
most $1+(2^{N}-1)\sum_{k}n_{k}$ leaves, where $n_{k}$ counts the
rejected cubes of side $2^{-k}R$. The nodes of such a cube have pairwise
radius at most $\sqrt{N}2^{-k}R$, so its mass exceeds
$(\rho2^{k}/(\sqrt{N}R))^{G+1}$. Disjointness gives
$n_{k}\le\bar{\chi}\,2^{-k(G+1)}$ with $\bar{\chi}=N^{(G+1)/2}\chi_{i}(\Gamma)$, and
the hypothesis gives $n_{k}\le B2^{kd'}$. The two geometric sequences
cross at the value $V=B^{(G+1)/(d'+G+1)}\bar{\chi}^{d'/(d'+G+1)}$, and their
minima sum to at most $c_{d',G}V$. If $\bar{\chi}<B$, then $\bar{\chi}\le V$ and the
same bound holds. The level $k$ grid has $2^{kN}$ cubes, which gives
the case $B=1$, $d'=N$. For item 2, an accepted patch of mass $w$
lies in the ball of radius $\rho w^{-1/(G+1)}$ about any of its nodes,
so $w\le A(\rho/R)^{d}w^{-d/(G+1)}+\omega_{\mathrm{at}}$. If
$w>2\omega_{\mathrm{at}}$, this gives
$w\le(2A)^{(G+1)/(d+G+1)}\chi_{i}(\Gamma)^{-d/(d+G+1)}$. The patch
masses sum to one.
\end{proof}

\begin{prop}\label{prop:discrete-zador}
In the setting of Proposition~\ref{prop:minimal-partition}, let a
cube of side $S$ contain the normalized frontier, and for $k\ge0$
let $\mathcal{Q}_{k}$ be its dyadic subcubes of side $\delta_{k}=2^{-k}S$,
with the face convention of item~1 of Proposition~\ref{prop:minimal-partition}. Let $\mathcal{Q}^{\ast}_{k}\subseteq\mathcal{Q}_{k}$
consist of the cubes that contain at least two distinct normalized
states, let $\bar{\mu}_{k}$ be the largest mass $\mu(Q)$ over $Q\in\mathcal{Q}^{\ast}_{k}$,
with $\bar{\mu}_{k}=0$ if $\mathcal{Q}^{\ast}_{k}=\varnothing$, and
put
\[
\Lambda=\left\{ k\ge0:\left(\frac{\rho}{\sqrt{N}\delta_{k}}\right)^{G+1}<\bar{\mu}_{k}\right\} .
\]
For $d>0$, with $\psi(d)$ as in Remark~\ref{rem:harmonic-exponent},
define the discrete Zador functional
\begin{equation}
\mathcal{Z}_{k}(d)=\delta^{\psi(d)}_{k}\sum_{Q\in\mathcal{Q}^{\ast}_{k}}\mu(Q)^{d/(d+G+1)}.\label{eq:discrete-zador}
\end{equation}
Then $\Lambda$ is finite, $\delta_{k}>\rho/\sqrt{N}$ for every $k\in\Lambda$,
and for every $d>0$ the dyadic construction in the proof of item~1
of Proposition~\ref{prop:minimal-partition}, started from this cube, returns an admissible partition with at most
\begin{equation}
1+\left(2^{N}-1\right)\sum_{k\in\Lambda}\min\left\{ \left|\mathcal{Q}^{\ast}_{k}\right|,\,N^{\psi(d)/2}\rho^{-\psi(d)}\mathcal{Z}_{k}(d)\right\} \label{eq:discrete-zador-bound}
\end{equation}
patches. In particular, \eqref{eq:discrete-zador-bound} bounds $K^{\min}_{i}(\Gamma,1)$.
\end{prop}

\begin{proof}
As in the proof of item~1 of Proposition~\ref{prop:minimal-partition}, the number of leaves is at most $1+(2^{N}-1)\sum_{k}n_{k}$,
where $n_{k}$ counts the rejected cubes of side $\delta_{k}$. A cube
whose nodes have at most one distinct normalized state has radius
zero and is accepted, so the rejected cubes of side $\delta_{k}$ belong
to $\mathcal{Q}^{\ast}_{k}$. The nodes of such a cube $Q$ have pairwise
radius at most $\sqrt{N}\delta_{k}$, so its rejection gives $\sqrt{N}\delta_{k}>\rho\mu(Q)^{-1/(G+1)}$,
that is, $\mu(Q)>\underline{w}_{k}=(\rho/(\sqrt{N}\delta_{k}))^{G+1}$.
Hence $n_{k}=0$ for $k\notin\Lambda$, and $\bar{\mu}_{k}\le1$ gives
$\delta_{k}>\rho/\sqrt{N}$ for $k\in\Lambda$, so $\Lambda$ is finite.
For $k\in\Lambda$, each rejected cube contributes more than one to
the sum below, and the other terms are nonnegative, so
\[
n_{k}\le\min\left\{ \left|\mathcal{Q}^{\ast}_{k}\right|,\sum_{Q\in\mathcal{Q}^{\ast}_{k}}\left(\frac{\mu(Q)}{\underline{w}_{k}}\right)^{d/(d+G+1)}\right\} ,\qquad\underline{w}^{-d/(d+G+1)}_{k}=N^{\psi(d)/2}\rho^{-\psi(d)}\delta^{\psi(d)}_{k},
\]
since $(G+1)d/(d+G+1)=\psi(d)$. Summing over $k\in\Lambda$ proves
\eqref{eq:discrete-zador-bound}.
\end{proof}

\begin{rem}\label{rem:discrete-zador}
(1) The two bounds in the proof of item~1 of Proposition~\ref{prop:minimal-partition}
are limiting cases. As $d\to\infty$,
the second entry of the minimum tends to a limit at most $N^{(G+1)/2}(S/\rho)^{G+1}2^{-k(G+1)}$,
which is the mass bound of that proof when $S=R$, and as $d\to0$ it tends
to $\left|\mathcal{Q}^{\ast}_{k}\right|$. By H\"older's inequality,
$\mathcal{Z}_{k}(d)\le\delta^{\psi(d)}_{k}\left|\mathcal{Q}^{\ast}_{k}\right|^{(G+1)/(d+G+1)}$,
so under the covering hypothesis of item~1 with $S=R$ and $d=d'$,
each term of \eqref{eq:discrete-zador-bound} is at most the crossing
value $B^{(G+1)/(d'+G+1)}(N^{(G+1)/2}\chi_{i}(\Gamma))^{d'/(d'+G+1)}$
of that proof. Unlike $\chi_{i}(\Gamma)$, the functional weighs each
cube by its mass. Low-mass nodes far from the bulk of the frontier
therefore enter \eqref{eq:discrete-zador-bound} only through the
counts $\left|\mathcal{Q}^{\ast}_{k}\right|$ at coarse levels and
through small powers of their masses.

(2) Replace the frontier by a probability with bounded density $h$
and compact support in $\mathbb{R}^{N}$, as in Remark~\ref{rem:harmonic-exponent}.
Then $\mathcal{Q}^{\ast}_{k}$ is the set of cubes of positive mass,
$\mathcal{Z}_{k}(N)\to\int h^{N/(N+G+1)}$ as $k\to\infty$ by Lebesgue's
differentiation theorem, and $k\in\Lambda$ requires $\delta^{N+G+1}_{k}>\rho^{G+1}N^{-(G+1)/2}\left\Vert h\right\Vert ^{-1}_{\infty}$.
For a fixed $m\ge1$, the levels $k\le\max\Lambda-m$ contribute at
most the geometric sum of the bounds $\left|\mathcal{Q}^{\ast}_{k}\right|\le2^{kN}$,
which is $O(2^{-mN}\rho^{-\psi(N)})$ by the last condition, and the
second entry bounds each of the other $m$ levels. Hence $\limsup_{\rho\to0}\rho^{\psi(N)}K^{\min}$
is at most a multiple of $\int h^{N/(N+G+1)}$, with a factor that
depends on $N$, $G$ and, logarithmically, on $h$. By Remark~\ref{rem:harmonic-exponent},
the same functional bounds $\liminf_{\rho\to0}\rho^{\psi(N)}K^{\min}$
from below.

(3) For an invertible linear map $W$, the image density $h_{W}=\left|\det W\right|^{-1}h\circ W^{-1}$
satisfies $\int h^{N/(N+G+1)}_{W}=\left|\det W\right|^{\psi(N)/N}\int h^{N/(N+G+1)}$.
For the WLL test, the unit map $D(c)$ of \eqref{eq:ver1-096} therefore
multiplies this integral by $\left|\det D(c)\right|^{\psi(N)/N}=c^{\psi(N)\sum_{j}\zeta_{j}/N}$.
For AWLL the relevant map is $W=\Gamma^{\mathrm{AWLL}}_{i}(c)^{-1/2}D(c)$,
and \eqref{eq:normalized-unit-bound} bounds $\left|\det W\right|\le\lVert W\rVert^{N}_{2}$
uniformly in $c\ge1$. This is the Zador form, as $\rho\to0$ at fixed
$c$, of the unit response in Theorem~\ref{thm:finite-scale-work}.
\end{rem}

\begin{cor}
\label{cor:scale-count}Let $i\in\mathcal{R}$, let $L=1$ and the mixed
block $\overline{\Gamma}_{i,1}$ be nonzero, and let $\gamma^{+}_{i}=\gamma^{+}_{i,1}$.
Suppose that, for $\Gamma=I_{N}$, the hypotheses of items 1 and 2 of
Proposition~\ref{prop:minimal-partition} hold with $d'=d$. Put
$X_{i}=(R_{i}(I_{N})/(\sqrt{\gamma^{+}_{i}}\rho_{i}))^{\psi(d)}$.
Then, with the golden ratio $\varphi$,
\begin{equation}
\min\left\{ \frac{1}{2\omega_{\mathrm{at}}},c_{1}\varphi^{-\psi(d)/2}X_{i}\right\} \le K^{\min}_{i}\left(\Gamma^{\mathrm{AWLL}}_{i},1\right)\le1+C_{1}\varphi^{\psi(d)/2}X_{i},\label{eq:scale-count}
\end{equation}
with $c_{1}=(2A)^{-(G+1)/(d+G+1)}$ and
$C_{1}=(2^{N}-1)c_{d,G}B^{(G+1)/(d+G+1)}N^{\psi(d)/2}$.
\end{cor}
\begin{proof}
A scalar $\lambda$ in \eqref{eq:partition-rule} replaces $\rho_{i}$
by $\lambda\rho_{i}$, and hence $\chi_{i}(I_{N})$ by $(R_{i}(I_{N})/(\lambda\rho_{i}))^{G+1}$.
Apply Proposition~\ref{prop:minimal-partition} with $\Gamma=I_{N}$
at $\lambda=\sqrt{\varphi \gamma^{+}_{i}}$ for the lower bound and at
$\lambda=\sqrt{\gamma^{+}_{i}/\varphi}$ for the upper bound, and combine
with \eqref{eq:forward-work-scale}, whose scalars are these two.
\end{proof}

\subsection{Repeated compression}
\begin{prop}
\label{prop:compression-regime}Let $\mathcal{R}=\{1,\ldots,n-1\}$,
let every frontier be subdivided by population bisection with $\lambda=1$,
let each $\mathbb{P}^{R,+}_{i,\kappa}$ be an extreme point of the set
of solutions of \eqref{eq:feasible-reduction}, and let every node
$[j]$ in the support of $\mathbb{P}^{R,+}_{i-1}$, $1\le i\le n$,
have exactly $b_{\mathrm{br}}\ge2$ labels $m\in\mathcal{B}_{i}$ with
$\pi_{i,m}(\widehat{Y}_{i-1}(j))>0$. For $i\in\mathcal{R}$ and $0<\eta_{\mathrm{tar}}<1$,
the retention fraction satisfies $\eta_{i}\le\eta_{\mathrm{tar}}$ whenever
\begin{equation}
M_{G}\left[1+\sum^{\lceil\log_{2}J^{-}_{i}\rceil-1}_{k=0}\min\left\{ 2^{k},\chi_{i}\left(\Gamma_{i}\right)\right\} \right]\le\eta_{\mathrm{tar}}J^{-}_{i}.\label{eq:certified-retention}
\end{equation}
For runs with unit factors $c$, suppose that the frontier of every
$i\in\mathcal{R}$ satisfies $\chi_{i}(\Gamma_{i}(c))\le\Xi_{n}$, with
$\Xi_{n}\ge1$ independent of $i$ and $c$. Put
\begin{equation}
a_{n}=M_{G}\left[1+\Xi_{n}\left(1+\log_{2}b_{\mathrm{br}}+\frac{1}{\log2}\right)\right],\qquad u_{n}=2a_{n}\left(1+\log\left(2a_{n}\right)\right).\label{eq:ver1-126}
\end{equation}
Then
\begin{equation}
\sup_{c}\max_{1\le i<n}J^{+}_{i}\left(c\right)\le u_{n},\qquad\sup_{c}J^{-}_{n}\left(c\right)\le b_{\mathrm{br}}u_{n}.\label{eq:multistep-load-bound}
\end{equation}
If $\Xi_{n}\le C_{\chi}n^{p_{\chi}}$ with constants $C_{\chi}>0$ and
$p_{\chi}\ge0$, then $\sup_{c}J^{-}_{n}(c)=O(n^{p_{\chi}}\log(n+1))$
for fixed $N$, $G$, $b_{\mathrm{br}}$ and $C_{\chi}$. If the input weights
at date $n-1$ are at most $\omega/J^{-}_{n-1}$, with $\omega$ independent
of $c$, then
\begin{equation}
\sup_{c}J^{-}_{n}\left(c\right)\le b_{\mathrm{br}}M_{G}\max\left\{ 1,4\omega\Xi_{n}\right\} ,\label{eq:ver1-128}
\end{equation}
which is $O(n^{p_{\chi}})$ in the same case.
\end{prop}
\begin{proof}
By \eqref{eq:recursive-work-upper} and $J^{+}_{i}\le M_{G}K_{i}$
(Theorem~\ref{thm:population}), the left side of \eqref{eq:certified-retention}
bounds $J^{+}_{i}$, which proves the first assertion. For the second,
let $x=J^{+}_{i-1}\ge1$. Since $J^{-}_{i}=b_{\mathrm{br}}x$ and $\chi_{i}(\Gamma_{i}(c))\le\Xi_{n}$,
the same bounds give
\begin{equation}
\begin{aligned}J^{+}_{i} & \le M_{G}\left[1+\Xi_{n}\lceil\log_{2}\left(b_{\mathrm{br}}x\right)\rceil\right]\\
 & \le M_{G}\left[1+\Xi_{n}\left(1+\log_{2}b_{\mathrm{br}}\right)+\frac{\Xi_{n}}{\log2}\log x\right]\le a_{n}\left(1+\log x\right).
\end{aligned}
\label{eq:ver1-129}
\end{equation}
Since $a_{n}\ge1$ and $\log t<t$ for $t=1+\log(2a_{n})$,
\begin{equation}
a_{n}\left(1+\log u_{n}\right)=a_{n}\left(1+\log\left(2a_{n}\right)+\log\left(1+\log\left(2a_{n}\right)\right)\right)\le u_{n}.\label{eq:ver1-130}
\end{equation}
Since $J^{+}_{0}=1\le u_{n}$ and the right side of \eqref{eq:ver1-129}
increases with $x$, induction gives $J^{+}_{i}\le u_{n}$ for $1\le i<n$,
and $J^{-}_{n}=b_{\mathrm{br}}J^{+}_{n-1}$ gives the second bound
of \eqref{eq:multistep-load-bound}.
If $\Xi_{n}\le C_{\chi}n^{p_{\chi}}$,
then $a_{n}=O(n^{p_{\chi}})$ uniformly in $c$. Finally, \eqref{eq:balanced-weight-upper}
at date $n-1$ gives $K_{n-1}\le\max\{1,4\omega\Xi_{n}\}$, and $J^{-}_{n}=b_{\mathrm{br}}J^{+}_{n-1}\le b_{\mathrm{br}}M_{G}K_{n-1}$
gives \eqref{eq:ver1-128}.
\end{proof}

\begin{cor}
\label{cor:refinement-exponent}Let $s_{i}=T/n$ for $1\le i\le n$
and $\mathcal{R}=\{1,\ldots,n-1\}$, let every patch of every date
in $\mathcal{R}$ pass \eqref{eq:partition-rule} with $\lambda=1$ or
have at most $M_{G}$ nodes, and let every node $[j]$ in the support
of $\mathbb{P}^{R,+}_{i-1}$, $1\le i\le n$, have exactly $b_{\mathrm{br}}$
labels $m\in\mathcal{B}_{i}$ with $\pi_{i,m}(\widehat{Y}_{i-1}(j))>0$.
Write $z=\Gamma^{-1/2}_{n-1}\widehat{Y}_{n-1}$ and $\mu_{n}$ for the
image of $\mathbb{P}^{R,-}_{n-1}$ under $z$, and note that $\rho_{n-1}=(T/n)^{(2G+1)/(2(G+1))}$
by \eqref{eq:resolution-load}.
\begin{enumerate}
\item Suppose $\sup_{x}\mu_{n}(\lVert z-x\rVert_{2}\le r)\le Ar^{d}+\omega_{\mathrm{at}}$
for $r>0$, with $A>0$, $d>0$ and $0<\omega_{\mathrm{at}}<1/2$. Then
\begin{equation}
J^{-}_{n}\ge b_{\mathrm{br}}\min\left\{ \frac{1}{M_{G}\omega_{\mathrm{at}}},(2A)^{-(G+1)/(d+G+1)}\rho^{-\psi(d)}_{n-1}\right\} .\label{eq:refinement-lower}
\end{equation}
\item Let each $\mathbb{P}^{R,+}_{n-1,\kappa}$ be an extreme point of
the set of solutions of \eqref{eq:feasible-reduction}. Suppose $R_{i}(\Gamma_{i})\le\bar{R}n^{q}$
for $1\le i<n$, with $\bar{R}>0$ and $q\ge0$, and put $p_{\chi}=(2G+1)/2+q(G+1)$
and $\Xi_{n}=\max\{1,\bar{R}^{G+1}T^{-(2G+1)/2}n^{p_{\chi}}\}$. If
the date $n-1$ is partitioned by the dyadic construction in item 1
of Proposition~\ref{prop:minimal-partition}, then
\begin{equation}
J^{-}_{n}\le b_{\mathrm{br}}M_{G}\left[1+(2^{N}-1)c_{N,G}\left(N^{(G+1)/2}\Xi_{n}\right)^{N/(N+G+1)}\right].\label{eq:refinement-dyadic}
\end{equation}
If instead it is partitioned by population bisection and its input
weights are at most $\omega/J^{-}_{n-1}$, then $J^{-}_{n}\le b_{\mathrm{br}}M_{G}\max\{1,4\omega\Xi_{n}\}$.
\end{enumerate}
\end{cor}
\begin{proof}
Each patch has positive mass and keeps at least one node, so $J^{-}_{n}=b_{\mathrm{br}}J^{+}_{n-1}\ge b_{\mathrm{br}}K_{n-1}$.
Each patch at date $n-1$ either passes \eqref{eq:partition-rule} or
has at most $M_{G}$ nodes. In the first case, the proof of item 2 of
Proposition~\ref{prop:minimal-partition}, with $AR^{d}$ in place of
$A$ and $R=R_{n-1}(\Gamma_{n-1})$, bounds its mass by $\max\{2\omega_{\mathrm{at}},(2A)^{(G+1)/(d+G+1)}\rho^{\psi(d)}_{n-1}\}$.
In the second case, its at most $M_{G}$ distinct normalized states
each carry mass at most $\omega_{\mathrm{at}}$. Since $M_{G}\ge2$
and the masses sum to one, $K_{n-1}\ge\min\{1/(M_{G}\omega_{\mathrm{at}}),(2A)^{-(G+1)/(d+G+1)}\rho^{-\psi(d)}_{n-1}\}$,
which gives \eqref{eq:refinement-lower}. For item 2, $T-t_{i}\ge T/n$
gives $\rho^{G+1}_{i}\ge(T/n)^{(2G+1)/2}$, hence $\chi_{i}(\Gamma_{i})\le\Xi_{n}$
for every $i\in\mathcal{R}$. Item 1 of Proposition~\ref{prop:minimal-partition}
with $B=1$ and $d'=N$, and $J^{-}_{n}=b_{\mathrm{br}}J^{+}_{n-1}\le b_{\mathrm{br}}M_{G}K_{n-1}$
from Theorem~\ref{thm:population}, give \eqref{eq:refinement-dyadic},
and \eqref{eq:balanced-weight-upper} gives the last bound in the same way.
\end{proof}
If $A$, $d$, $\omega_{\mathrm{at}}$, $\bar{R}$, $q$ and $\omega$ do not
depend on $n$, the lower bound \eqref{eq:refinement-lower} grows like $n^{\psi(d)(2G+1)/(2(G+1))}$
until it reaches $b_{\mathrm{br}}/(M_{G}\omega_{\mathrm{at}})$, whereas
the upper bound \eqref{eq:refinement-dyadic} is of order $n^{Np_{\chi}/(N+G+1)}$,
which is $n^{\psi(N)(2G+1)/(2(G+1))}$ when $q=0$, and the upper bound
for population bisection is of order $n^{p_{\chi}}$.

\FloatBarrier

\section{Numerical evidence}

\label{sec:numerics}

\subsection{Test systems, initial values and settings}
\label{sec:models}
\noindent
The test problems are the Heston and sSABR systems below, with time
in years. The first coordinates of their initial values are the initial
price $S_{0}=100$ for Heston and the initial forward rate $F_{0}=0.01$
for sSABR. Each problem is also solved after a change of units by
two other factors $c$, which replaces $Y_{0}$ by $D(c)Y_{0}$, and the
first coordinate of the initial value, called the level, is then 0.01,
1 or 100 (Table~\ref{tab:equivalent-units}).

\begin{table}[!htbp]\centering\small
\caption{Parameters at the three levels. The unit factor \(c\) maps the row with \(c=1\) to the others. Initial volatility denotes \(\nu_0\) for Heston and \(\alpha_0\) for sSABR. }\label{tab:equivalent-units}
\begin{tabular}{llrrrr}\toprule
Model & Level & \(c\) & Initial volatility & Shift \(s\) & Sensitivity parameter\\\midrule
Heston & 0.01 & \(10^{-4}\) & 0.20 & n/a & \(\nu_0=0.20\)\\
Heston & 1 & \(10^{-2}\) & 0.20 & n/a & \(\nu_0=0.20\)\\
Heston & 100 & 1 & 0.20 & n/a & \(\nu_0=0.20\)\\
sSABR & 0.01 & 1 & 0.03 & 0.03 & \(\omega_0=0.006\)\\
sSABR & 1 & \(10^2\) & 0.30 & 3 & \(\omega_0=0.6\)\\
sSABR & 100 & \(10^4\) & 3.00 & 300 & \(\omega_0=60\)\\
\bottomrule\end{tabular}
\end{table}

\subsubsection{Heston}
\noindent
The Heston system \citep{heston_closed-form_1993} has initial volatility
$\nu_{0}=0.20$, spot drift $\mu=0$, mean-reversion rate $\kappa=3$,
long-run variance $\theta=0.04$, variance diffusion coefficient $\xi=0.4$
and Brownian correlation $\rho=-0.7$. Its six coordinates are the
spot $Y_{1}$, the variance $Y_{2}$, their $\nu_{0}$ tangents $Y_{3},Y_{4}$
(Remark~\ref{rem:greeks}), and the running integrals $Y_{5,t}=\int^{t}_{0}Y_{1,u}\,du$,
$Y_{6,t}=\int^{t}_{0}Y_{3,u}\,du$. The initial value is $Y_{0}=\prescript{t}{}{\left(S_{0},\nu^{2}_{0},0,2\nu_{0},0,0\right)}$.

The parameters satisfy the Feller condition $2\kappa\theta=0.24>\xi^{2}=0.16$
\citep{feller_1951}, so the
variance stays positive, and the solution of \eqref{eq:augmented}
stays in $\mathcal{O}=\{y:y_{2}>0\}$, on which the coefficients below
are smooth.

For the NV rules used below, positivity of the variance updates
is ensured under the weaker condition $4\kappa\theta>\xi^{2}$
\citep[Section~4.1]{ninomiya_high-order_2025}, which these parameters
also satisfy.
The closed-form update formulas admit extensions beyond the times
at which the classical solutions of the corresponding ordinary
differential equations reach $y_{2}=0$.
Boundary treatments of such CIR-type splitting schemes and their weak
approximation properties are discussed by \citet[Section~2.1]{alfonsi_2010}
and \citet{alfonsi_lombardo_2025}.
When the Feller condition fails, the variance can reach zero. The
NV scheme then remains well defined as long as $\xi^{2}\le4\kappa\theta$.
\citet[Section~2]{alfonsi_2010} also constructs modified weak second-order
schemes for unrestricted CIR parameters.
\citet{alfonsi_2005,alfonsi_2010,alfonsi_2015}
and \citet{alfonsi_lombardo_2024,alfonsi_lombardo_2025}
study the discretization
of such variance processes and of the Heston model, including
NV-type schemes and their high-order extensions on random grids.

On $y_{2}>0$, the It\^o matrix representation is 
\begin{equation}
\begin{aligned}dY_{t} & =b^{\mathrm{Heston}}\left(Y_{t}\right)dt+\Sigma^{\mathrm{Heston}}\left(Y_{t}\right)dB_{t},\\
b^{\mathrm{Heston}}\left(y\right) & =\left(\begin{array}{c}
\mu y_{1}\\
\kappa\left(\theta-y_{2}\right)\\
\mu y_{3}\\
-\kappa y_{4}\\
y_{1}\\
y_{3}
\end{array}\right),\qquad\Sigma^{\mathrm{Heston}}\left(y\right)=\left(\begin{array}{cc}
y_{1}\sqrt{y_{2}} & 0\\
\rho\xi\sqrt{y_{2}} & \sqrt{1-\rho^{2}}\,\xi\sqrt{y_{2}}\\
\sqrt{y_{2}}y_{3}+\dfrac{y_{1}y_{4}}{2\sqrt{y_{2}}} & 0\\[4pt]
\dfrac{\rho\xi y_{4}}{2\sqrt{y_{2}}} & \dfrac{\sqrt{1-\rho^{2}}\,\xi y_{4}}{2\sqrt{y_{2}}}\\[4pt]
0 & 0\\
0 & 0
\end{array}\right).
\end{aligned}
\label{eq:heston-ito}
\end{equation}
 Here $N=6$ and $d_{B}=2$. By the display after \eqref{eq:augmented},
$V^{\mathrm{Heston}}_{1}I_{N}$ and $V^{\mathrm{Heston}}_{2}I_{N}$ are the columns
of $\Sigma^{\mathrm{Heston}}$, and
\begin{equation}
V^{\mathrm{Heston}}_{0}I_{N}\left(y\right)=b^{\mathrm{Heston}}\left(y\right)-\prescript{t}{}{\left(\left(y_{2}/2+\rho\xi/4\right)y_{1},\xi^{2}/4,\left(y_{2}/2+\rho\xi/4\right)y_{3}+y_{1}y_{4}/2,0,0,0\right)}.\label{eq:heston-fields}
\end{equation}

A unit factor $c>0$ leaves the parameters unchanged and maps the
states and strikes by
\begin{equation}
Y_{t}'=D_{\mathrm{H}}\left(c\right)Y_{t},\qquad D_{\mathrm{H}}\left(c\right)=\operatorname{diag}\left(c,1,c,1,c,c\right),\qquad K'=cK.\label{eq:ver1-149}
\end{equation}
 The factors $c=10^{-4},10^{-2},1$ give the levels 0.01, 1 and 100.
Prices and $\nu_{0}$ model vegas are multiplied by $c$.

\subsubsection{sSABR}

The sSABR system, the SABR model of \citet{hagan_managing_2002} with
deterministic shift $s$ \citep[Section~2]{cui_kirkby_nguyen_2019},
has elasticity $\beta=0.50$, volatility of volatility $\nu=0.30$
and Brownian correlation $\rho=0.20$, and, for $c=1$, initial volatility
$\alpha_{0}=0.03$ and shift $s=0.03$. Its first two coordinates are
the forward and the volatility amplitude. With $\omega_{0}=\alpha_{0}\left(F_{0}+s\right)^{\beta}$,
the coordinates $Y_{3},Y_{4}$ are their $\omega_{0}$ tangents at fixed
$F_{0}$, $s$ and $\beta$, and $Y_{5},Y_{6}$ are the running integrals
of $Y_{1},Y_{3}$. The initial value is $Y_{0}=\prescript{t}{}{\left(F_{0},\alpha_{0},0,\left(F_{0}+s\right)^{-\beta},0,0\right)}$.

On $y_{1}+s>0$ and $y_{2}>0$, the It\^o matrix representation is 
\begin{equation}
\begin{aligned}dY_{t} & =b^{\mathrm{sSABR}}\left(Y_{t}\right)dt+\Sigma^{\mathrm{sSABR}}\left(Y_{t}\right)dB_{t},\\
b^{\mathrm{sSABR}}\left(y\right) & =\left(\begin{array}{c}
0\\
0\\
0\\
0\\
y_{1}\\
y_{3}
\end{array}\right),\qquad\Sigma^{\mathrm{sSABR}}\left(y\right)=\left(\begin{array}{cc}
y_{2}\left(y_{1}+s\right)^{\beta} & 0\\
\rho\nu y_{2} & \sqrt{1-\rho^{2}}\,\nu y_{2}\\
\left(y_{1}+s\right)^{\beta}y_{4}+\beta y_{2}\left(y_{1}+s\right)^{\beta-1}y_{3} & 0\\
\rho\nu y_{4} & \sqrt{1-\rho^{2}}\,\nu y_{4}\\
0 & 0\\
0 & 0
\end{array}\right).\end{aligned}
\label{eq:sabr-ito}
\end{equation}

By the display after \eqref{eq:augmented}, $V^{\mathrm{sSABR}}_{1}I_{N}$
and $V^{\mathrm{sSABR}}_{2}I_{N}$ are the columns of $\Sigma^{\mathrm{sSABR}}$,
and
\begin{equation}
\begin{aligned}V^{\mathrm{sSABR}}_{0}I_{N} & =b^{\mathrm{sSABR}}-\prescript{t}{}{\left(c_{1},c_{2},c_{3},c_{4},0,0\right)},\qquad c_{2}\left(y\right)=\tfrac{1}{2}\nu^{2}y_{2},\qquad c_{4}\left(y\right)=\tfrac{1}{2}\nu^{2}y_{4},\\
c_{1}\left(y\right) & =\tfrac{1}{2}\beta y^{2}_{2}\left(y_{1}+s\right)^{2\beta-1}+\tfrac{1}{2}\rho\nu y_{2}\left(y_{1}+s\right)^{\beta},\\
c_{3}\left(y\right) & =\left[\tfrac{1}{2}\beta\left(2\beta-1\right)y^{2}_{2}\left(y_{1}+s\right)^{2\beta-2}+\tfrac{1}{2}\rho\nu\beta y_{2}\left(y_{1}+s\right)^{\beta-1}\right]y_{3}\\
 & \quad+\left[\beta y_{2}\left(y_{1}+s\right)^{2\beta-1}+\tfrac{1}{2}\rho\nu\left(y_{1}+s\right)^{\beta}\right]y_{4}.
\end{aligned}
\label{eq:sabr-fields}
\end{equation}

For the constant elasticity of variance diffusion with $1/2\le\beta<1$,
zero is an absorbing boundary \citep[Section~2.1]{chen_oosterlee_weide_2012},
so $Y_{1}+s$ is absorbed at zero, and valuation
includes the probability mass accumulated there \citep[Sections~1 and~2]{gulisashvili_horvath_jacquier_2018}.
From the absorption time on, $Y_{1}=-s$, $Y_{3}=0$, $dY_{5}=-s\,dt$
and $dY_{6}=0$, while $Y_{2}$ and $Y_{4}$ continue their multiplicative
dynamics. The hypothesis after \eqref{eq:augmented} therefore holds
only up to the absorption time, and
\eqref{eq:semigroup} is the semigroup of the absorbed process. No tree node and no reference path reached
the face $\{y:y_{1}=-s\}$ in the reported runs, so every tree state lies in
$\mathcal{O}=\{y:y_{1}+s>0,\ y_{2}>0\}$.

A unit factor $c>0$ leaves $\beta$, $\nu$ and $\rho$ unchanged and
gives
\begin{equation}
\begin{aligned}F_{0}' & =cF_{0}, & s' & =cs, & \alpha_{0}' & =c^{1-\beta}\alpha_{0}, & K' & =cK,\\
D_{\mathrm{S}}\left(c\right) & =\operatorname{diag}\left(c,c^{1-\beta},1,c^{-\beta},c,1\right), & \omega_{0}' & =c\omega_{0}.
\end{aligned}
\label{eq:ver1-155}
\end{equation}
 The factors $c=1,10^{2},10^{4}$ give the levels 0.01, 1 and 100.
Prices are multiplied by $c$, and the $\omega_{0}$ model vega is invariant. At fixed $F_{0},s,\beta$, $\partial_{\alpha_{0}}C=\left(F_{0}+s\right)^{\beta}\partial_{\omega_{0}}C$.

For $c=1$, the initial diffusion amplitudes of the first coordinate
are $S_{0}\nu_{0}=20$ for Heston and $\omega_{0}=0.006$ for sSABR,
per square root year. All reported errors are converted to $c=1$,
and Deltas are invariant.

\subsubsection{Tasks and comparison settings}

The tasks are European and continuously monitored arithmetic Asian
calls and puts with maturity $T=1$, zero carry and discount factor
one. WLL uses the WLL test \eqref{eq:raw-rule}, with $M_{i}=I_{N}$
in the computational coordinates, and AWLL uses the matrix test \eqref{eq:partition-rule}
with $\Gamma_{i}=\Gamma^{\mathrm{AWLL}}_{i}$ and $L=1$, so that $M_{i}=(\Gamma^{\mathrm{AWLL}}_{i})^{-1/2}$
also sets the splitting directions, both with $\lambda=1$. Global
recombination, called Full, uses one patch containing the whole frontier.
All three methods use degree $G=5$ and uniform grids with $n=1,\ldots,9$
intervals and reduction at every interior date, which gives 162 cases.

All trees use $\mathrm{Alg}=\mathrm{NV}$, whose state operator $Q^{\mathrm{NV}}_{i}$
is \eqref{eq:ver1-007} with the branch labels and probabilities
\begin{equation*}
\mathcal{B}_{i}=\left\{ -\sqrt{3},0,\sqrt{3}\right\} ^{d_{B}}\times\left\{ +,-\right\} ,\qquad\pi_{i,\left(\eta,\pm\right)}=\frac{1}{2}6^{-\#\left\{ a:\eta_{a}\ne0\right\} }\left(2/3\right)^{\#\left\{ a:\eta_{a}=0\right\} }.
\end{equation*}
For Heston the maps and vector fields in the following
symmetric rule are restricted to the first four coordinates.
\begin{equation*}
\begin{aligned}\Phi_{i,\left(\eta,+\right)}\left(y\right) & =\exp\left(\frac{s_{i}}{2}V_{0}\right)\exp\left(\sqrt{s_{i}}\eta_{1}V_{1}\right)\cdots\exp\left(\sqrt{s_{i}}\eta_{d_{B}}V_{d_{B}}\right)\exp\left(\frac{s_{i}}{2}V_{0}\right)\left(y\right),\\
\Phi_{i,\left(\eta,-\right)}\left(y\right) & =\exp\left(\frac{s_{i}}{2}V_{0}\right)\exp\left(\sqrt{s_{i}}\eta_{d_{B}}V_{d_{B}}\right)\cdots\exp\left(\sqrt{s_{i}}\eta_{1}V_{1}\right)\exp\left(\frac{s_{i}}{2}V_{0}\right)\left(y\right).
\end{aligned}
\end{equation*}
The two running integrals are updated by the
trapezoidal rule using the price and its tangent at the step endpoints,
respectively. This update composes the exact flow over $s_{i}/2$ of
the vector field formed by the last two components of $V^{\mathrm{Heston}}_{0}$,
the rule for the first four coordinates and the same flow, a symmetric
composition that keeps the local weak error of order $s^{3}_{i}$
\citep[Theorem~1.17 and Section~4.2]{alfonsi_2010}. For sSABR the maps are those of the rule of
\citet[Eq.~(2.2)]{ninomiya_high-order_2025},
\begin{equation*}
\begin{aligned}\Phi_{i,\left(\eta,+\right)}\left(y\right) & =\exp\left(s_{i}V_{0}\right)\exp\left(\sqrt{s_{i}}\eta_{1}V_{1}\right)\cdots\exp\left(\sqrt{s_{i}}\eta_{d_{B}}V_{d_{B}}\right)\left(y\right),\\
\Phi_{i,\left(\eta,-\right)}\left(y\right) & =\exp\left(\sqrt{s_{i}}\eta_{d_{B}}V_{d_{B}}\right)\cdots\exp\left(\sqrt{s_{i}}\eta_{1}V_{1}\right)\exp\left(s_{i}V_{0}\right)\left(y\right),
\end{aligned}
\end{equation*}
where $\exp\left(V\right)\left(y\right)$ is the value at time one
of the solution of $dz/dt=VI_{N}\left(z\right)$, $z_{0}=y$, and,
for Heston, of the extension of its closed-form formula beyond the
time at which it reaches $y_{2}=0$. For
sSABR, these flows are computed with the forward replaced by its Lamperti
coordinate $2\sqrt{y_{1}+s}$, by the explicit fifth-order Runge--Kutta
method of \citet{butcher_1987}, first applied to such flows by \citet{ninomiya_ninomiya_2009}
and also used in \citet{ma_ninomiya_2026}. For vector fields with
coefficients in $C^{\infty}_{b}$, this method, used in place of the
exact flows, changes every map $\Phi_{i,m}$ by $O(s^{3}_{i})$ uniformly
\citep[Theorems~1.3 and 4.15]{ninomiya_ninomiya_2009}, the same order
as the local weak error of the NV rules. The
coordinate values $\pm\sqrt{3}$ and $0$ have probabilities $1/6$
and $2/3$, matching the standard normal moments through degree five.
With $d_{B}=2$, every parent has 18 positive branches. The labels
$(\eta,+)$ and $(\eta,-)$ give the same state if $\eta=0$ and, for the
symmetric rule, if some component of $\eta$ vanishes.

Let $y_{n,j,h}=\left(\widehat{Y}_{n}\left(j\right)\right)_{h}$ be component
$h$ at terminal node $j$. The price and first model vega estimators
are 
\begin{equation}
\begin{aligned}\widehat{C}^{A} & =\sum_{j}p^{R,-}_{n,j}\left(y_{n,j,5}/T-K\right)_{+},\\
\widehat{\mathcal{V}}^{A} & =\sum_{j}p^{R,-}_{n,j}\mathbf{1}_{\left\{ y_{n,j,5}/T>K\right\} }\,y_{n,j,6}/T,
\end{aligned}
\label{eq:asian-estimators}
\end{equation}
 where the vega is taken with respect to $\nu_{0}$ for Heston and
$\omega_{0}$ for sSABR. The vega estimator is the target of Remark~\ref{rem:greeks}, since
the average equals $K$ with probability zero, and no terminal node
has $y_{n,j,5}/T=K$. European versions use
coordinates one and three instead of five and six divided by $T$.
Put versions use the corresponding put payoff and signed indicator.

For either terminal underlying
or arithmetic average $Z$, let $P\left(K\right)$ denote the corresponding
call or put price and define $q\left(K\right)=\mathbb{P}\left(Z>K\right)$
for a call and $q\left(K\right)=-\mathbb{P}\left(Z<K\right)$ for
a put, so that $\partial_{K}P(K)=-q(K)$ at zero carry and discount
factor one. Prices are homogeneous of degree one under the unit maps
\eqref{eq:ver1-149} and \eqref{eq:ver1-155}, and an sSABR price depends
on $F_{0}$, $K$ and $s$ only through $F_{0}+s$ and $K+s$, so Euler's
theorem on homogeneous functions gives
\begin{equation}
\Delta^{\mathrm{Heston}}\left(K\right)=\frac{P\left(K\right)+Kq\left(K\right)}{S_{0}},\qquad\Delta^{\mathrm{sSABR}}\left(K\right)=\frac{P\left(K\right)+\left(K+s\right)q\left(K\right)-\left(1-\beta\right)\omega_{0}\partial_{\omega_{0}}P\left(K\right)}{F_{0}+s}.\label{eq:numerical-delta}
\end{equation}
The Heston derivative holds $\nu_{0}$ and the other model parameters
fixed. The sSABR forward derivative holds $\alpha_{0},s,\beta,\nu,\rho$
fixed, so that $\omega_{0}=\alpha_{0}\left(F_{0}+s\right)^{\beta}$
varies with $F_{0}$. Both the tree and the quasi-Monte Carlo benchmarks
evaluate the Deltas by \eqref{eq:numerical-delta} from prices, signed
exercise probabilities and, for sSABR, the $\omega_{0}$ model vega
of \eqref{eq:asian-estimators}.

Each plotted error is the maximum absolute difference over 401 equally
spaced strikes for the named option side. For $c=1$, the strike ranges
are $\left[0,300\right]$ for Heston and $\left[0,0.03\right]$ for
sSABR. Summary tables also maximize over call and put. Heston European
benchmarks use characteristic-function prices and analytical derivatives
\citep{cui_heston_2017}. Heston Asian and both sSABR products use
the symmetric NV scheme on $2^{10}$ equal time intervals, integrated with $2^{20}$
points of an unscrambled Sobol sequence with the direction numbers
of \citet{joe_kuo_2008}, omitting its initial zero point, and with
Brownian increments built by the Brownian bridge.
The reported discrepancies are relative to these finite-resolution
benchmarks.

\subsection{Support growth and the compression transition}\label{sec:support-growth}

\begin{figure}[!htbp]
\centering
\includegraphics[width=0.98\columnwidth]{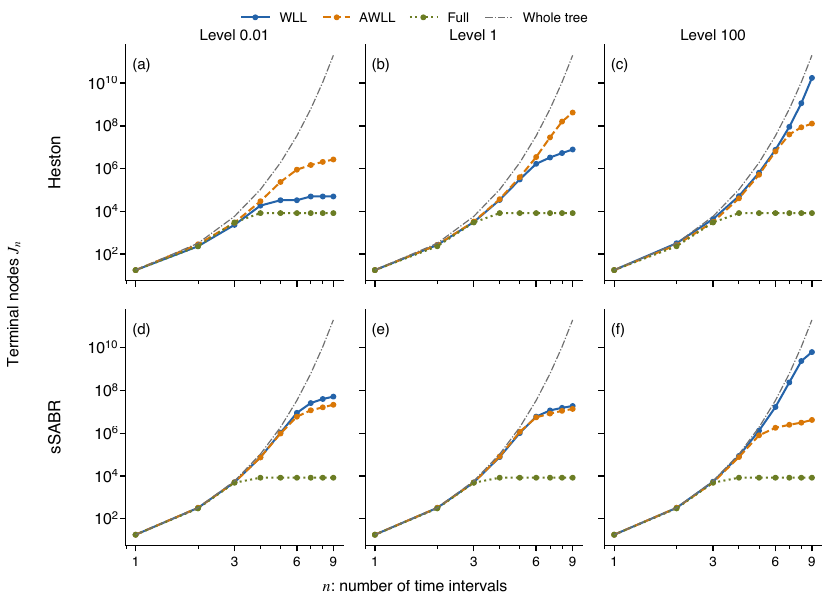}

\caption{The number \(J_{n}\) of terminal nodes on logarithmic axes.
Panels (a) to (c) show Heston and panels (d) to (f) show sSABR at
the levels 0.01, 1 and 100. Whole tree denotes the unrecombined count
\(18^{n}\), and Full denotes global recombination. \label{fig:all-supports}}
\end{figure}

Figure~\ref{fig:all-supports} shows the number $J_{n}=J^{-}_{n}$
of terminal nodes for the grid with $n$ intervals.
With $\lambda=1$, WLL has very different supports at the three
levels, whereas AWLL has smaller supports than WLL at level 100 and
larger ones for Heston at the two smaller levels.
Table~\ref{tab:growth-ols-tls} gives estimates of the exponent
$x$ in $J_{n}\approx an^{x}$ on the windows $m{:}9=\{m,\ldots,9\}$
of $n$, $m=1,4,6,7$. Ordinary least squares (OLS) regresses $\log J_{n}$
on $\log n$, and total least squares (TLS) minimizes orthogonal squared
distances in these unstandardized log coordinates. An exponent that
falls on the later windows marks the onset of compression. At level
100, WLL keeps exponents above 13 on the three later windows in both
models. For Heston WLL at level 100, the OLS fit of $\log J_{n}$ over all
nine grids gives $R^{2}=0.99966$ against $n$, with slope $k$ satisfying
$e^{k}=12.767$, and $R^{2}=0.91858$ against $\log n$, with slope $x=9.296$.
Its large exponents on the later windows therefore describe nearly
geometric growth on a finite range of grids.

By \eqref{eq:propagated-capacity}, $J_{n}=18M_{G}K_{n-1}\upsilon_{n-1}$,
where $M_{G}\upsilon_{n-1}$ is the average number of parents retained
per patch at date $n-1$. At nine intervals and level 100, this average
is 147.114 for Heston WLL and 461.922 for AWLL, against the nominal
capacity $M_{G}=\binom{11}{5}=462$, and 432.813 and 462.000 in sSABR.
The larger WLL supports at level 100 therefore result from
a greater number of patches
(Table~\ref{tab:work}) rather than from more retained nodes
per patch.
Full reaches $18M_{G}=8316$
terminal nodes at four intervals and remains at that population. On
the grid with two intervals, the first reduction keeps 13 of the 18
nodes for Heston and 17 for sSABR, so that $J_{2}=234$ and $306$,
and the Full rows of Table~\ref{tab:growth-ols-tls} differ between
the two models.

Sustained half-tree suppression starts at the first $n$
from which $J_{m}/18^{m}\le1/2$ for every observed $m\ge n$, and
near-capacity operation starts at the first $n$ from which
$\upsilon_{m-1}\ge0.95$ for every observed $m\ge n$. Once utilization
changes little, $J_{n}$ is nearly proportional to $K_{n-1}$.

\begin{table}[!htbp]\centering\small
\caption{Fitted exponents \(x\) in \(J_n\approx a n^x\). Each cell gives the OLS/TLS estimates on the indicated window. The four windows contain nine, six, four and three observations. Full has the same values at all three levels. }\label{tab:growth-ols-tls}
\begin{tabular}{lllrrrr}\toprule Model&Level&Method&\(1{:}9\)&\(4{:}9\)&\(6{:}9\)&\(7{:}9\)\\\midrule
Heston&0.01&WLL&3.86/4.09 & 1.21/1.32 & 0.94/1.20 & 0.00/0.00\\
Heston&0.01&AWLL&5.95/6.07 & 5.43/5.87 & 2.70/2.72 & 2.38/2.38\\
Heston&1&WLL&6.48/6.62 & 6.65/7.01 & 3.80/3.82 & 3.43/3.43\\
Heston&1&AWLL&8.11/8.52 & 11.87/11.90 & 11.97/12.09 & 10.67/10.82\\
Heston&100&WLL&9.30/10.11 & 15.52/15.79 & 19.00/19.13 & 20.92/20.99\\
Heston&100&AWLL&7.89/8.23 & 10.40/10.69 & 7.35/7.88 & 4.67/4.77\\
Heston&all&Full&2.82/3.27 & 0.00/0.00 & 0.00/0.00 & 0.00/0.00\\
\midrule
sSABR&0.01&WLL&7.54/7.75 & 8.19/8.73 & 4.27/4.54 & 2.78/2.81\\
sSABR&0.01&AWLL&7.09/7.28 & 6.92/7.55 & 3.11/3.19 & 2.38/2.38\\
sSABR&1&WLL&7.06/7.25 & 6.72/7.42 & 2.79/2.92 & 1.89/1.90\\
sSABR&1&AWLL&6.91/7.11 & 6.09/6.90 & 2.25/2.27 & 2.01/2.01\\
sSABR&100&WLL&9.43/10.03 & 14.39/14.47 & 14.96/15.24 & 13.09/13.60\\
sSABR&100&AWLL&6.21/6.39 & 4.53/5.21 & 2.05/2.05 & 2.08/2.09\\
sSABR&all&Full&2.73/3.28 & 0.00/0.00 & 0.00/0.00 & 0.00/0.00\\
\bottomrule\end{tabular}
\end{table}

At level 100, sustained half-tree suppression starts at $n=4$ for
both Heston methods, and at $n=5$ and $n=6$ for sSABR AWLL and WLL,
respectively. Near-capacity operation starts later, at $n=8$ for
Heston AWLL and $n=6$ for sSABR AWLL, whereas the utilization of WLL
stays below 0.95 up to $n=9$ in both models. Heston AWLL at level
one reaches near-capacity operation only at $n=9$, so its large exponents
on the later windows still include the transition to capacity.
For sSABR, every node lies on the hyperplane
$\omega_{0}y_{4}=y_{2}$, so the exact moment rank is at most
$252\approx0.545M_{G}$, and utilization above this value reflects the
numerical rank kept by batched basis exchange (Section~\ref{sec:global-lambda}).

For sSABR AWLL at level 100, both nested windows $6{:}9$ and $7{:}9$
lie within its near-capacity segment, and the OLS/TLS estimates $2.048/2.054$
and $2.075/2.093$ on them are close to two. For Heston AWLL
at level 100, the estimates fall from $7.353/7.883$ to $4.674/4.775$,
since its near-capacity operation starts only at $n=8$. These fitted exponents describe
the observed grids $n\le9$ and do not determine an asymptotic growth
rate. Heston WLL
has a constant support on $7{:}9$ at level 0.01, and estimates between
three and four on $6{:}9$ and $7{:}9$ at level one.

For sSABR at level one, AWLL has the smaller support but a slightly
larger exponent on $7{:}9$. On $n=7,8,9$, its support is respectively
$0.710$, $0.723$ and $0.730$ times that of WLL, and the difference
$2.008-1.895=0.113$ of the OLS exponents equals the OLS slope of the
logarithm of this ratio. Both methods reach exponents near two from
larger exponents on the earlier windows. The support and its exponent
measure different features of this transition.

\subsection{Accuracy and computational work}

Tables~\ref{tab:accuracy} and~\ref{tab:work} give the comparison
at level 100 with $\lambda=1$, with WLL as the baseline. The support
ratios WLL/AWLL are 137.05 for Heston and 1488.18 for sSABR, and the
storage ratios are about 50 and 590. For $n=5,\ldots,9$, the AWLL/WLL
ratios of the four price and model vega reference errors are at most
1.088 for Heston and 1.063 for sSABR. Thus AWLL reduces the tree greatly
at reference errors close to those of WLL. Full reduces the tree further,
with larger reference errors in every entry of Table~\ref{tab:accuracy}.

\begin{table}[!htbp]
\centering\small
\caption{Price and model vega reference errors at level 100 and nine intervals. Each value is the maximum absolute difference over 401 strikes and both option sides. }\label{tab:accuracy}
\begin{tabular}{llrrrr}
\toprule
Model&Method&European price&European vega&Asian price&Asian vega\\
\midrule
Heston&Full&\(5.144\times10^{-2}\)&\(8.501\times10^{-1}\)&\(1.075\times10^{-1}\)&\(1.137\times10^{0}\)\\
Heston&WLL&\(3.694\times10^{-2}\)&\(3.058\times10^{-1}\)&\(5.068\times10^{-2}\)&\(4.881\times10^{-1}\)\\
Heston&AWLL&\(3.603\times10^{-2}\)&\(3.096\times10^{-1}\)&\(5.513\times10^{-2}\)&\(4.758\times10^{-1}\)\\
sSABR&Full&\(3.467\times10^{-5}\)&\(2.634\times10^{-2}\)&\(2.045\times10^{-5}\)&\(1.793\times10^{-2}\)\\
sSABR&WLL&\(3.028\times10^{-5}\)&\(1.806\times10^{-2}\)&\(5.604\times10^{-6}\)&\(2.954\times10^{-3}\)\\
sSABR&AWLL&\(2.869\times10^{-5}\)&\(1.785\times10^{-2}\)&\(5.955\times10^{-6}\)&\(3.051\times10^{-3}\)\\
\bottomrule
\end{tabular}
\end{table}

\begin{table}[!htbp]\centering\small
\caption{Final patches (the patches of date 8), terminal nodes and storage demand at level 100 and nine intervals. Storage demand adds uncompressed node payload and direct stage and case metadata over all dates. }\label{tab:work}
\begin{tabular}{llrrr}\toprule
Model & Method & Final patches & Terminal nodes & Storage (GiB)\\\midrule
Heston & WLL & 6,570,690 & 17,399,545,704 & 897.592\\
Heston & AWLL & 15,269 & 126,955,656 & 17.781\\
Heston & Full & 1 & 8,316 & 0.006\\
sSABR & WLL & 799,036 & 6,224,993,316 & 491.797\\
sSABR & AWLL & 503 & 4,182,948 & 0.833\\
sSABR & Full & 1 & 8,316 & 0.006\\
\bottomrule\end{tabular}
\end{table}

At level 100 with nine intervals, the elapsed times
were 477 and 43 hours for Heston WLL and AWLL,
respectively, and 560 and 5.9 hours for sSABR
WLL and AWLL.
These times include scheduling and recovery between
attempts and depend on the computing environment.
An AWLL reduction takes
more CPU time on average than a WLL reduction, since its patches are
larger and retain more parents, 462 against 147 on average for Heston
(Section~\ref{sec:support-growth}), and the cost in Proposition~\ref{prop:kernel-work}
grows with the patch population and rank.

\begin{figure}[!htbp]
\centering
\includegraphics[width=0.98\columnwidth]{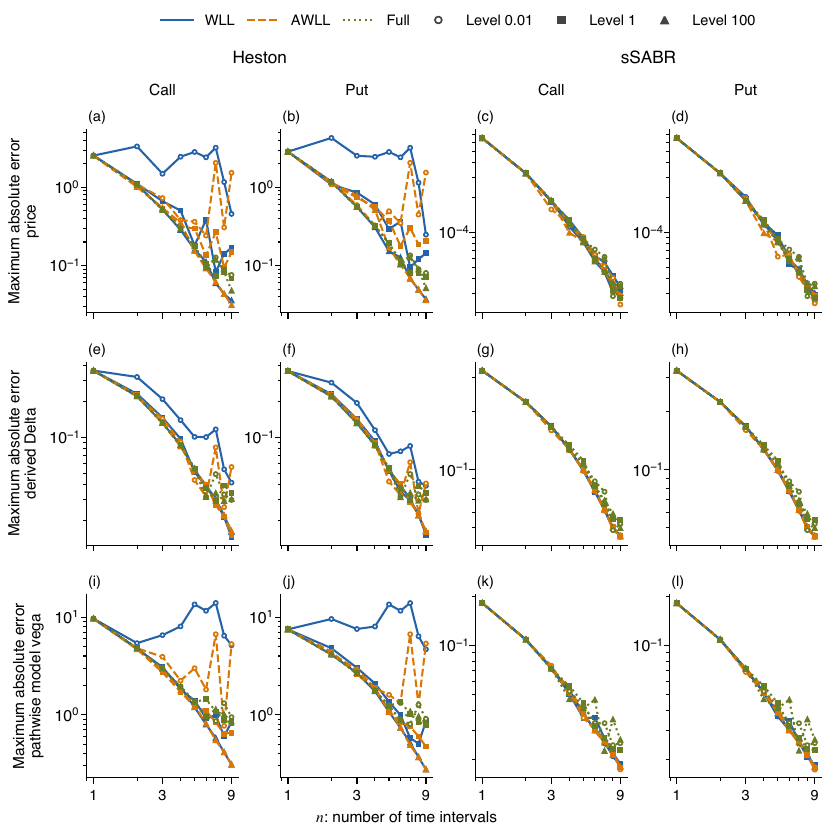}

\caption{European task errors.
Rows show maximum absolute errors for price, derived Delta and pathwise
model vega. Columns show Heston call, Heston put, sSABR call and sSABR
put. Each curve shows the maximum error over 401 strikes for the indicated option side. Both
axes are logarithmic. All errors are converted to $c=1$. Line color
and style identify the method, and marker shape identifies the level. Heston
uses spot Delta and $\nu_{0}$ vega, while sSABR uses forward Delta
at fixed $\alpha_{0}$ and $\omega_{0}$ vega. \label{fig:european-errors}}
\end{figure}

\begin{figure}[!htbp]
\centering
\includegraphics[width=0.98\columnwidth]{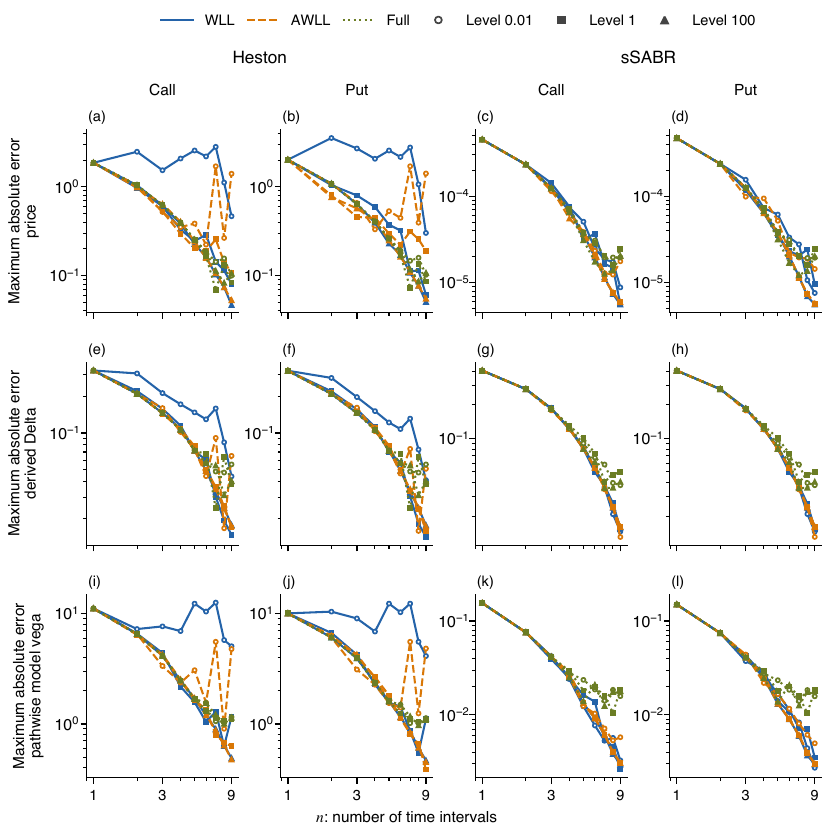}

\caption{Arithmetic Asian task errors. The panel
arrangement, units, maxima over strikes, curve encodings and Greeks
follow Figure~\ref{fig:european-errors}. \label{fig:asian-errors}}
\end{figure}

Figures~\ref{fig:european-errors} and~\ref{fig:asian-errors} show
all three levels and both option sides. In both figures, the errors
of WLL and AWLL decrease with $n$ at level 100, whereas for Heston
at the two smaller levels they are larger and depend more on the level.
For the arithmetic Asian tasks, the sSABR price and model vega errors
of Full stop decreasing, whereas those of WLL and AWLL are close at
all levels. For one model,
level, method, product and task, let $E_{n}$ be the maximum error
over strikes and sides. The effective order is the negative OLS slope
of $\log E_{n}$ against $\log n$, and Table~\ref{tab:temporal-orders}
gives it on the windows $1{:}9$ and $6{:}9$ at level 100.

\begin{table}[!tbp]\centering\small\setlength{\tabcolsep}{3pt}
\caption{Effective orders at level 100. Each entry gives the orders on \(1{:}9\) and \(6{:}9\). }\label{tab:temporal-orders}
\begin{tabular}{llrrrrrr}\toprule
&&\multicolumn{3}{c}{European}&\multicolumn{3}{c}{Arithmetic Asian}\\
\cmidrule(lr){3-5}\cmidrule(lr){6-8}
Model & Method & Price & Delta & Model vega & Price & Delta & Model vega\\\midrule
Heston & WLL & \(2.04\,/\, 3.01\) & \(1.45\,/\, 2.23\) & \(1.59\,/\, 2.37\) & \(1.68\,/\, 2.95\) & \(1.33\,/\, 2.51\) & \(1.46\,/\, 2.23\)\\
Heston & AWLL & \(2.08\,/\, 2.57\) & \(1.46\,/\, 2.11\) & \(1.59\,/\, 2.36\) & \(1.65\,/\, 2.67\) & \(1.33\,/\, 2.56\) & \(1.47\,/\, 2.28\)\\
Heston & Full & \(1.81\,/\, 1.76\) & \(1.29\,/\, 0.20\) & \(1.15\,/\, 0.32\) & \(1.50\,/\, 1.57\) & \(1.06\,/\, 1.20\) & \(1.17\,/\, 0.86\)\\
sSABR & WLL & \(1.43\,/\, 1.71\) & \(0.93\,/\, 1.40\) & \(1.08\,/\, 1.27\) & \(2.11\,/\, 3.20\) & \(1.48\,/\, 3.06\) & \(1.86\,/\, 2.85\)\\
sSABR & AWLL & \(1.45\,/\, 1.85\) & \(0.93\,/\, 1.37\) & \(1.08\,/\, 1.28\) & \(2.09\,/\, 3.10\) & \(1.47\,/\, 3.02\) & \(1.84\,/\, 2.77\)\\
sSABR & Full & \(1.38\,/\, 1.44\) & \(0.86\,/\, 1.23\) & \(0.90\,/\, 0.97\) & \(1.77\,/\, -0.35\) & \(1.23\,/\, 0.97\) & \(1.11\,/\, 0.14\)\\
\bottomrule\end{tabular}
\end{table}

On $6{:}9$, the orders of WLL and AWLL range from 2.11 to 3.01 for
Heston, and, for sSABR, from 2.77 to 3.20 for the Asian tasks and
from 1.27 to 1.85 for the European tasks. For Full, the orders on
$6{:}9$ of the sSABR Asian price and model vega are $-0.35$ and 0.14. The single patch of Full has unit mass and the whole frontier
as its support, so the recombination term of \eqref{eq:joint-error}
does not tend to zero as $n$ grows, and Theorem~\ref{thm:residual}
gives no convergence for Full. These slopes on finite grids describe
the recorded errors, which include the common propagation bias, the
recombination residual and the error of the benchmarks.

\subsection{AWLL behavior and global WLL equivalents}
\label{sec:covariance-evidence}
\label{sec:global-lambda}
\begin{figure}[!htbp]
\centering
\includegraphics[width=0.95\columnwidth]{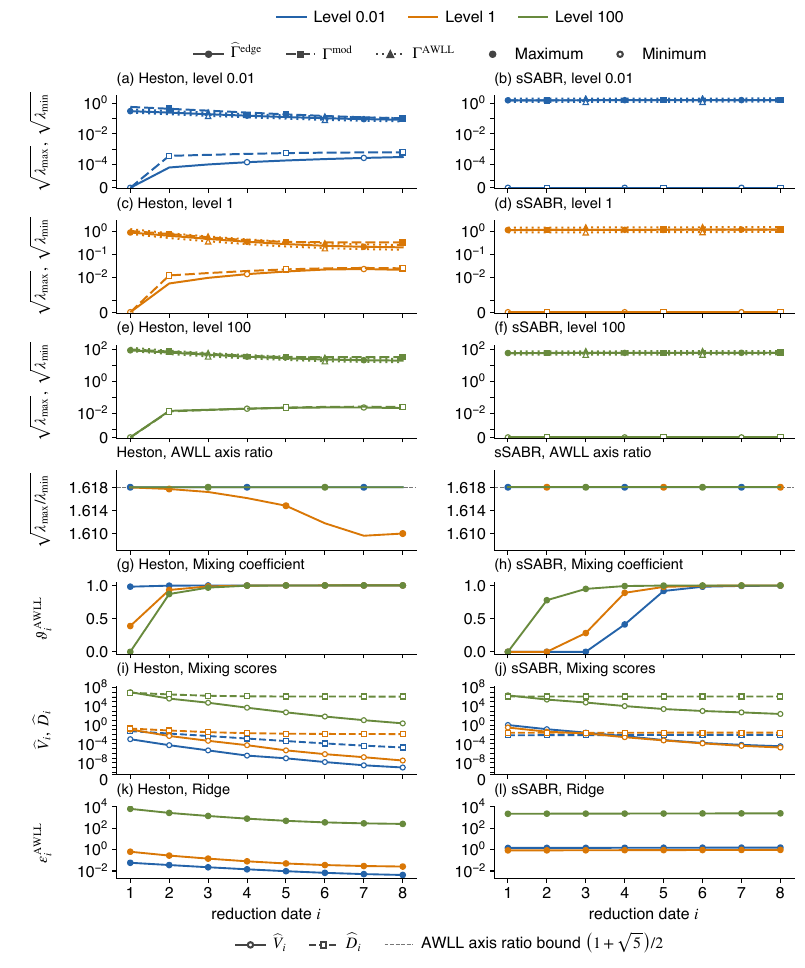}

\caption{AWLL matrix diagnostics and scalar coefficients evaluated on the six-state
parent measures for nine time intervals. Panels (a) through
(f) show the square roots of the largest and smallest eigenvalues
of $\widehat{\Gamma}^{\mathrm{edge}}_{i}$, $\Gamma^{\mathrm{mod}}_{i}$
and $\Gamma^{\mathrm{AWLL}}_{i}$, with separate rows for the
levels 0.01, 1 and 100. Line and marker styles distinguish the matrices.
Filled and open markers denote maxima and minima, respectively. Panels
(g) and (h) show $\vartheta^{\mathrm{AWLL}}_{i}$, (i) and (j) show
$\widehat{V}_{i},\widehat{D}_{i}$, and (k) and (l) show $\varepsilon^{\mathrm{AWLL}}_{i}$. The two unlettered panels between (e), (f) and (g), (h) show the AWLL axis ratio $\left(\lambda_{\max}\left(\Gamma^{\mathrm{AWLL}}_{i}\right)/\lambda_{\min}\left(\Gamma^{\mathrm{AWLL}}_{i}\right)\right)^{1/2}$ and its bound $\varphi$, the square root of the bound in \eqref{eq:ridge-condition}.
Colors identify the level. Spectra and scores use logarithmic
axes with linear neighborhoods of zero, and the ridge axis is logarithmic. \label{fig:wms-behavior}}
\end{figure}

\noindent
Figure~\ref{fig:wms-behavior} combines
the extreme eigenvalues of the two references and the selected AWLL matrix
with the scalar coefficients from \eqref{eq:mixing-coefficient}
and Theorem~\ref{thm:wms-derivation}.
Figure~\ref{fig:joint-states} in Appendix~\ref{app:metric-action}
compares the identity and AWLL metrics on the same six-state measures
at two dates and all three levels. Table~\ref{tab:mechanism} gives the final mixing scores,
regularized spectra and full-dimensional distance ratios.
Lemma~\ref{lem:metric-scalar}~(2) converts these eigenvalues into
WLL scalars (Section~\ref{sec:selected-scale}).

At level 100, the largest eigenvalue of the edge reference falls
from 7879 to 407 for Heston and rises from 3635 to 3959 for sSABR between
$i=1$ and $i=8$. It also decreases at both smaller Heston levels. For sSABR, the nodes satisfy $\omega_{0}\left(\widehat{Y}_{i}\right)_{4}=\left(\widehat{Y}_{i}\right)_{2}$
with maximum relative residual below $2.6\times10^{-15}$, so both references
have the null vector $e_{2}-\omega_{0}e_{4}$, whereas the AWLL matrix
remains positive definite. The first
three mixing coefficients at the smallest sSABR level are zero. By
the last date, the model weights $1-\vartheta^{\mathrm{AWLL}}_{i}$
at level 100 are $1.99\times10^{-6}$ for Heston and $1.91\times10^{-4}$
for sSABR.
By \eqref{eq:mixing-coefficient}, these weights equal \(\widehat{V}_{8}/\widehat{D}_{8}\),
and \(\widehat{D}_{8}\) is large because, by Proposition~\ref{prop:reference-discrepancy},
the running integral \(Y_5\) of \(Y_1\), which has no diffusion and
whose integrand has mean near 100, gives
\(\widehat{\Gamma}^{\mathrm{edge}}_8-\Gamma^{\mathrm{mod}}_8\) a \(\left(5,5\right)\)
entry of about \(-s_8\bar{b}_{8,5}^{2}\approx-100^{2}/9\), whose square
is 99.80\% of \(\widehat{D}_{8}\) for Heston and 92.22\% for sSABR.

The effect of the metric is measured on the full incoming probability
measure. With $p_{j}=p^{R,-}_{i,j}$ for the nodes $\left[j\right]$
in its support, the root mean square pair distance ratio is
\begin{equation}
\left(\frac{\sum_{j,j'}p_{j}p_{j'}\left\lVert M^{\mathrm{AWLL}}_{i}\left(\widehat{Y}_{i}\left(j\right)-\widehat{Y}_{i}\left(j'\right)\right)\right\rVert ^{2}_{2}}{\sum_{j,j'}p_{j}p_{j'}\left\lVert \widehat{Y}_{i}\left(j\right)-\widehat{Y}_{i}\left(j'\right)\right\rVert ^{2}_{2}}\right)^{1/2}=\left(\frac{\operatorname{tr}\left(M^{\mathrm{AWLL}}_{i}\operatorname{Cov}\left(\widehat{Y}_{i}\right)\prescript{t}{}{M^{\mathrm{AWLL}}_{i}}\right)}{\operatorname{tr}\left(\operatorname{Cov}\left(\widehat{Y}_{i}\right)\right)}\right)^{1/2},\label{eq:metric-distance-ratio}
\end{equation}
where $\operatorname{Cov}\left(\widehat{Y}_{i}\right)=\operatorname{Cov}_{\mathbb{P}^{R,-}_{i}}\left(\widehat{Y}_{i}\right)$.
Each double sum is twice the corresponding trace, so the factor two
cancels.
The statistics of this subsection and of Appendix~\ref{app:metric-action}
are second moments. Positive reduction
preserves them on every patch, so they take the same values on the
incoming and reduced measures of a date, up to the numerical moment
residual.
Table~\ref{tab:mechanism}
compares the metric and identity on identical incoming states. At
Heston level 100, the weighted standard deviation of $\left(\widehat{Y}_{i}\right)_{3}$ increases
from 28.85 to 41.01 between $i=1$ and $i=8$, while the $(3,3)$ entry
of the edge reference decreases from 7492 to 93.71. The resulting distance
ratio increases from 0.00775 to 0.0491. The states keep their earlier
spread while the edge reference contracts.

\begin{table}[!htbp]\centering\small\setlength{\tabcolsep}{3pt}
\caption{AWLL mixing and its action on the complete incoming six-state measures for \(n=9\). Panel (a) gives the final mixing scores. Panel (b) gives the ratio \eqref{eq:metric-distance-ratio} of root mean square pair distances under AWLL and \(I_6\) on the same incoming measure at each indicated date. The eigenvalues are those of \(\Gamma_8^{\mathrm{AWLL}}\), and the terminal population ratio compares the two complete histories. }\label{tab:mechanism}
\begin{minipage}{0.5\linewidth}
\begin{tabular*}{\linewidth}{@{\extracolsep{\fill}}lrrr@{}}\toprule
\multicolumn{4}{l}{(a) Mixing scores}\\
Model&Level&\(\widehat V_8\)&\(\widehat D_8\)\\\midrule
Heston & 0.01 & \(9.77\times10^{-10}\) & \(1.8046\times10^{-5}\)\\
Heston & 1 & \(3\times10^{-8}\) & \(1.2744\times10^{-2}\)\\
Heston & 100 & \(2.45\times10^{0}\) & \(1.2348\times10^{6}\)\\
sSABR & 0.01 & \(3.29\times10^{-5}\) & \(8.8944\times10^{-3}\)\\
sSABR & 1 & \(1.77\times10^{-5}\) & \(2.6576\times10^{-2}\)\\
sSABR & 100 & \(2.44\times10^{2}\) & \(1.2798\times10^{6}\)\\
\bottomrule\end{tabular*}
\end{minipage}\par\medskip
\begin{tabular}{lrrrrrr}\toprule
\multicolumn{7}{l}{(b) Spectrum, pair distances and retained population}\\
Model&Level&\(\lambda_{\min}\)&\(\lambda_{\max}\)&RMS ratio, \(i=1\)&RMS ratio, \(i=8\)&\(J_9^{\mathrm{AWLL}}/J_9^{\mathrm{WLL}}\)\\\midrule
Heston & 0.01 & \(4.26\times10^{-3}\) & \(1.11\times10^{-2}\) & \(2.55\times10^{0}\) & \(1.04\times10^{1}\) & \(5.37\times10^{1}\)\\
Heston & 1 & \(2.66\times10^{-2}\) & \(6.9\times10^{-2}\) & \(8.04\times10^{-1}\) & \(4.87\times10^{0}\) & \(5.32\times10^{1}\)\\
Heston & 100 & \(2.52\times10^{2}\) & \(6.59\times10^{2}\) & \(7.75\times10^{-3}\) & \(4.91\times10^{-2}\) & \(7.3\times10^{-3}\)\\
sSABR & 0.01 & \(1.54\times10^{0}\) & \(4.04\times10^{0}\) & \(5.59\times10^{-1}\) & \(5.56\times10^{-1}\) & \(4.16\times10^{-1}\)\\
sSABR & 1 & \(9.39\times10^{-1}\) & \(2.46\times10^{0}\) & \(6.88\times10^{-1}\) & \(7.1\times10^{-1}\) & \(7.3\times10^{-1}\)\\
sSABR & 100 & \(2.45\times10^{3}\) & \(6.41\times10^{3}\) & \(1.31\times10^{-2}\) & \(1.38\times10^{-2}\) & \(6.72\times10^{-4}\)\\
\bottomrule\end{tabular}
\end{table}

At the last reduction date, $\Gamma^{\mathrm{AWLL}}_{8}\prec I_{6}$
at both smaller Heston levels, so AWLL enlarges every pair distance
relative to WLL. Its root mean square distance ratios are 10.43 and
4.87. At Heston level 100, $\Gamma^{\mathrm{AWLL}}_{8}\succ I_{6}$
and the ratio is 0.0491. The sSABR ratios are 0.556, 0.710 and 0.0138. At
sSABR level one, 98.3\% of weighted state variance lies in the eigendirections
that AWLL contracts. The terminal AWLL/WLL population ratios in Table~\ref{tab:mechanism} reflect repeated differences in division and
retention. At Heston level one, where AWLL retains 53.2 times as many
terminal nodes as WLL, its Asian model vega error is 0.62499 against
1.16277, whereas its Asian price error is 0.19046 against 0.08037.
The larger retained population therefore does not lower every error.

For each $n$, the post hoc WLL equivalent $\widehat{\lambda}^{\mathrm{AWLL}}_{n}$
of AWLL is the scalar of a complete WLL run from the initial value
whose history best fits the AWLL patch size distributions and retained
populations at $i=1,\ldots,n-1$. The response of WLL to the scalar
is taken from its complete histories at the three levels, which include
the effects of earlier reductions on later frontiers.

Let $\mathcal{A}_{n}\left(c\right)=(\operatorname{tr}\Gamma^{\mathrm{mod}}_{1})^{1/2}$
for the grid with $n$ intervals at the factor $c$, with the model reference
evaluated at the common initial state. The
fitting model represents effective isotropic amplitude by $\mathcal{A}_{n}\left(c\right)/\lambda$. It is exact only for a uniform unit map
$D(c)=cI_{N}$, under which the WLL test at $c$ and $\lambda$ equals
that at $c=1$ and $\lambda/c$.
Piecewise linear interpolation in log amplitude of the log WLL response
profiles, with endpoint extrapolation, predicts each complete history.
One scalar minimizes the joint squared log discrepancies in patch
counts, retained populations and patch size quantiles at probabilities
$0.1,\ldots,0.9$. Dates have equal weight, as do patch counts, retained
populations and the mean squared log discrepancy over the nine quantiles.
The population prediction preserves \eqref{eq:retention-history}.
Table~\ref{tab:global-lambda} gives the resulting estimates.

\begin{table}[!htbp]\centering\small
\caption{Post hoc WLL equivalents \(\widehat{\lambda}^{\mathrm{AWLL}}_{n}\) of AWLL. Each scalar is shared by the reduction dates \(1,\ldots,{n-1}\). The column heads 0.01, 1 and 100 are the levels. }\label{tab:global-lambda}
\begin{tabular}{crrrrrr}\toprule
&\multicolumn{3}{c}{Heston}&\multicolumn{3}{c}{sSABR}\\
\cmidrule(lr){2-4}\cmidrule(lr){5-7}
\(n\)&0.01&1&100&0.01&1&100\\\midrule
2 & 0.12 & 0.32 & 26 & 0.1 & 1 & 77\\
3 & 0.58 & 0.14 & 14 & 1.1 & 1 & 76\\
4 & 0.46 & 0.11 & 14 & 1.2 & 1 & 73\\
5 & 0.55 & 0.14 & 17 & 1.2 & 1.1 & 75\\
6 & 0.56 & 0.16 & 20 & 1.2 & 1.1 & 78\\
7 & 0.56 & 0.17 & 23 & 1.3 & 1.1 & 79\\
8 & 0.55 & 0.18 & 25 & 1.3 & 1.1 & 80\\
9 & 0.56 & 0.17 & 26 & 1.3 & 1.1 & 80\\
\bottomrule\end{tabular}
\end{table}

At $n=9$ and level 100, the largest multiplicative
discrepancy between fitted and observed retained populations over
the interior dates is 3.11 for Heston and 2.04 for sSABR.

For comparison, at level 100 and $n=9$, bounding the pairwise radius
$\operatorname{rad}_{I_{N}}$ of each complete Full frontier by its
coordinate box diagonal shows
that WLL with the same propagation, deterministic reducer and node
ordering accepts the whole frontier as one patch at every interior
date, and so reproduces the Full history, for every scalar of at
least 3260 for Heston and 7320 for sSABR.

To test the prediction \(\widehat{\lambda}^{\mathrm{AWLL}}_{9}\) of
Table~\ref{tab:global-lambda}, WLL is run at level 100 with \(\lambda=26\)
for Heston and \(\lambda=80\) for sSABR.
The remaining numerical parameters are those of Section~\ref{sec:models},
with \(n=9\). Table~\ref{tab:scalar-validation-work} compares terminal
populations, and
Table~\ref{tab:scalar-validation-errors} gives the separate call and put errors.

\begin{table}[!htbp]
\centering\small
\caption{Terminal populations at level 100 with nine intervals. }\label{tab:scalar-validation-work}
\begin{tabular}{llr}
\toprule
Model & Method & Terminal nodes\\
\midrule
Heston & WLL, \(\lambda=1\) & \num{17399545704}\\
Heston & AWLL, \(\lambda=1\) & \num{126955656}\\
Heston & Full & \num{8316}\\
Heston & WLL, \(\lambda=26\) & \num{74338578}\\
\midrule
sSABR & WLL, \(\lambda=1\) & \num{6224993316}\\
sSABR & AWLL, \(\lambda=1\) & \num{4182948}\\
sSABR & Full & \num{8316}\\
sSABR & WLL, \(\lambda=80\) & \num{3659040}\\
\bottomrule
\end{tabular}
\end{table}

\begin{table}[!htbp]
\centering\small\setlength{\tabcolsep}{4pt}
\caption{Maximum absolute reference errors at level 100 with nine intervals. }\label{tab:scalar-validation-errors}
\begin{tabular}{llrrrr}
\toprule
\multicolumn{6}{l}{(a) Heston}\\
Product & Task & WLL \(\lambda=1\) & AWLL \(\lambda=1\) & Full & WLL \(\lambda=26\)\\
\midrule
European call & Price & \(3.578\times10^{-2}\) & \(3.134\times10^{-2}\) & \(4.758\times10^{-2}\) & \(3.130\times10^{-2}\)\\
 & Delta & \(1.602\times10^{-2}\) & \(1.633\times10^{-2}\) & \(3.026\times10^{-2}\) & \(1.626\times10^{-2}\)\\
 & Model vega & \(3.058\times10^{-1}\) & \(3.096\times10^{-1}\) & \(8.501\times10^{-1}\) & \(3.082\times10^{-1}\)\\
\addlinespace[2pt]
European put & Price & \(3.694\times10^{-2}\) & \(3.603\times10^{-2}\) & \(5.144\times10^{-2}\) & \(3.600\times10^{-2}\)\\
 & Delta & \(1.609\times10^{-2}\) & \(1.638\times10^{-2}\) & \(3.021\times10^{-2}\) & \(1.630\times10^{-2}\)\\
 & Model vega & \(2.741\times10^{-1}\) & \(2.741\times10^{-1}\) & \(8.146\times10^{-1}\) & \(2.727\times10^{-1}\)\\
\addlinespace[2pt]
Asian call & Price & \(4.677\times10^{-2}\) & \(5.290\times10^{-2}\) & \(1.052\times10^{-1}\) & \(5.290\times10^{-2}\)\\
 & Delta & \(1.771\times10^{-2}\) & \(1.735\times10^{-2}\) & \(4.040\times10^{-2}\) & \(1.736\times10^{-2}\)\\
 & Model vega & \(4.881\times10^{-1}\) & \(4.758\times10^{-1}\) & \(1.137\times10^{0}\) & \(4.763\times10^{-1}\)\\
\addlinespace[2pt]
Asian put & Price & \(5.068\times10^{-2}\) & \(5.513\times10^{-2}\) & \(1.075\times10^{-1}\) & \(5.513\times10^{-2}\)\\
 & Delta & \(1.767\times10^{-2}\) & \(1.733\times10^{-2}\) & \(4.038\times10^{-2}\) & \(1.733\times10^{-2}\)\\
 & Model vega & \(4.652\times10^{-1}\) & \(4.510\times10^{-1}\) & \(1.112\times10^{0}\) & \(4.515\times10^{-1}\)\\
\midrule
\multicolumn{6}{l}{(b) sSABR}\\
Product & Task & WLL \(\lambda=1\) & AWLL \(\lambda=1\) & Full & WLL \(\lambda=80\)\\
\midrule
European call & Price & \(3.028\times10^{-5}\) & \(2.866\times10^{-5}\) & \(3.463\times10^{-5}\) & \(2.871\times10^{-5}\)\\
 & Delta & \(4.510\times10^{-2}\) & \(4.538\times10^{-2}\) & \(4.990\times10^{-2}\) & \(4.533\times10^{-2}\)\\
 & Model vega & \(1.806\times10^{-2}\) & \(1.784\times10^{-2}\) & \(2.634\times10^{-2}\) & \(1.777\times10^{-2}\)\\
\addlinespace[2pt]
European put & Price & \(3.006\times10^{-5}\) & \(2.869\times10^{-5}\) & \(3.467\times10^{-5}\) & \(2.874\times10^{-5}\)\\
 & Delta & \(4.510\times10^{-2}\) & \(4.538\times10^{-2}\) & \(4.989\times10^{-2}\) & \(4.533\times10^{-2}\)\\
 & Model vega & \(1.799\times10^{-2}\) & \(1.785\times10^{-2}\) & \(2.633\times10^{-2}\) & \(1.778\times10^{-2}\)\\
\addlinespace[2pt]
Asian call & Price & \(5.596\times10^{-6}\) & \(5.955\times10^{-6}\) & \(2.045\times10^{-5}\) & \(5.956\times10^{-6}\)\\
 & Delta & \(1.584\times10^{-2}\) & \(1.598\times10^{-2}\) & \(4.068\times10^{-2}\) & \(1.615\times10^{-2}\)\\
 & Model vega & \(2.954\times10^{-3}\) & \(3.051\times10^{-3}\) & \(1.787\times10^{-2}\) & \(3.087\times10^{-3}\)\\
\addlinespace[2pt]
Asian put & Price & \(5.604\times10^{-6}\) & \(5.724\times10^{-6}\) & \(2.022\times10^{-5}\) & \(5.725\times10^{-6}\)\\
 & Delta & \(1.584\times10^{-2}\) & \(1.598\times10^{-2}\) & \(4.068\times10^{-2}\) & \(1.615\times10^{-2}\)\\
 & Model vega & \(2.952\times10^{-3}\) & \(2.983\times10^{-3}\) & \(1.793\times10^{-2}\) & \(3.019\times10^{-3}\)\\
\bottomrule
\end{tabular}
\end{table}

At these scalars, the terminal populations are 58.6\% and 87.5\% of
their AWLL counterparts. WLL with \(\lambda=1\) retains 234 and 1701
times as many terminal nodes as these runs. Across both products and option sides, the price,
Delta and model vega reference errors computed from the values in
Table~\ref{tab:scalar-validation-errors} differ from those of AWLL
by at most 0.51\% for Heston and 1.21\% for sSABR, relative to AWLL.
The direct computations therefore
recover the AWLL resolution scale at comparable accuracy.
These runs do not compare the two tests as methods. AWLL is adaptive
and selects its scale during the run, whereas the scalars 26 and 80
are post hoc values, fitted to the complete AWLL runs, which WLL by
itself could find only by a sweep over $\lambda$.
Full has the smallest tree but larger errors in every reported task.

In the numerical implementation, a patch whose moment matrix has
numerical rank equal to its node count, \(k_{\tau}=J\), is carried
forward unchanged, as described after \eqref{eq:moments}. Such a patch
has at most \(M_{G}\) nodes, since \(k_{\tau}\le M_{G}\). Table~\ref{tab:identity-forwarding} reports these
events over the eight interior reduction dates of the experiments with
\(n=9\). The completion records show no fallback to the original patch
measure after a failed reduction in the 112 stages covered by the table.
Every patch that the population rule stopped without passing
\eqref{eq:partition-rule}
and with \(k_{\tau}<J\) had \(k_{\tau}\) equal to its number of distinct
states. Its reduction therefore only redistributed weight among nodes
with equal states, which leaves its image on \(\mathcal{O}\) unchanged,
as assumed before \eqref{eq:joint-wms-budget}.

\begin{table}[!htbp]
\centering\small\setlength{\tabcolsep}{4pt}
\caption{Patches carried forward unchanged under the moment rank criterion for \(n=9\). Counts are summed over \(i=1,\ldots,8\), and shares are percentages of all patches processed at those dates. }\label{tab:identity-forwarding}
\begin{tabular}{lrlrrrr}
\toprule
Model & Level & Method & \(\lambda\) & Total patches & Forwarded unchanged & Share (\%)\\
\midrule
Heston & 0.01 & WLL & 1 & 39 & 0 & 0\\
Heston & 0.01 & AWLL & 1 & 678 & 1 & 0.1475\\
Heston & 1 & WLL & 1 & 1,746 & 3 & 0.1718\\
Heston & 1 & AWLL & 1 & 69,797 & 4 & 0.005731\\
Heston & 100 & WLL & 1 & 8,448,395 & 1,137,727 & 13.47\\
Heston & 100 & AWLL & 1 & 23,497 & 2 & 0.008512\\
Heston & 100 & WLL & 26 & 15,223 & 8 & 0.05255\\
\midrule
sSABR & 0.01 & WLL & 1 & 12,436 & 127 & 1.021\\
sSABR & 0.01 & AWLL & 1 & 5,136 & 29 & 0.5646\\
sSABR & 1 & WLL & 1 & 4,550 & 15 & 0.3297\\
sSABR & 1 & AWLL & 1 & 3,414 & 7 & 0.2050\\
sSABR & 100 & WLL & 1 & 1,321,134 & 35,031 & 2.652\\
sSABR & 100 & AWLL & 1 & 1,100 & 6 & 0.5455\\
sSABR & 100 & WLL & 80 & 974 & 6 & 0.6160\\
\bottomrule
\end{tabular}
\end{table}

By the relation \(\omega_{0}(\widehat{Y}_{i})_{4}=(\widehat{Y}_{i})_{2}\),
the moment matrices of degree five for sSABR have exact rank at most
\(\binom{10}{5}=252\).
The rank threshold of Appendix~\ref{app:algorithm} does not separate
the singular values that this relation leaves at the rounding level
from the others, and it retained up to \(M_{G}=462\) numerical constraints.
The sSABR supports of all three methods therefore exceed those of
an exact rank reduction, and the sSABR capacity utilizations in
Section~\ref{sec:support-growth} refer to this numerical rank. With the
factor \(10^{-12}\) in place of \(10^{-18}\), the sSABR Full patch
has rank 252, while the Heston Full patch keeps rank 462. The method
comparisons share one solver setting, and every output still satisfies
the moment equations of degree at most \(G\).

\subsection{A synthetic test of directional shape}
\label{sec:shape-test}
\noindent
In Sections~\ref{sec:support-growth} to \ref{sec:covariance-evidence}
the comparison is dominated by the resolution scale, and the continuation
derivatives that determine the effect of the shape through
\eqref{eq:joint-error}
are not available. This subsection uses a system in which the AWLL
scale is close to one and the continuations are explicit, and compares
WLL and AWLL on two tasks whose continuations vary in different directions.

Let $N=3$ and $d_{B}=2$, and let the state consist of a planar Brownian
motion and its L\'evy area,
\begin{equation}
V_{0}=0,\qquad V_{1}I_{3}\left(y\right)=\prescript{t}{}{\left(1,0,-\frac{y_{2}}{2}\right)},\qquad V_{2}I_{3}\left(y\right)=\prescript{t}{}{\left(0,1,\frac{y_{1}}{2}\right)},\qquad Y_{0}=0,\label{eq:heisenberg-fields}
\end{equation}
so that $(Y_{1,t},Y_{2,t})=(B^{1}_{t},B^{2}_{t})$ and $Y_{3,t}=\frac{1}{2}\int^{t}_{0}\left(Y_{1,v}\,dY_{2,v}-Y_{2,v}\,dY_{1,v}\right)$.
The flows are explicit, and the two NV rules of Section~\ref{sec:models}
give the same operator $Q^{\mathrm{NV}}_{i}$ since $V_{0}=0$. On a uniform grid they move $(Y_{1},Y_{2})$ on the lattice $(3T/n)^{1/2}\mathbb{Z}^{2}$, and
$Y_{3}$ by an increment affine in the parent values of $(Y_{1},Y_{2})$.
The two task functions are
\begin{equation}
f^{\mathrm{area}}\left(y\right)=\cos\left(4y_{3}\right),\qquad f^{\mathrm{end}}\left(y\right)=\cos\left(2y_{1}\right)+\cos\left(2y_{2}\right),\label{eq:shape-tasks}
\end{equation}
with targets $\operatorname{sech}2=0.2658$ and $2e^{-2}=0.2707$ at
$T=1$. The characteristic function
of the area conditional on the Brownian increment \citep{levy_1951},
followed by a Gaussian integral over that increment, gives
\begin{equation}
h^{\mathrm{area}}_{t}\left(y\right)=\frac{\cos\left(4y_{3}\right)}{\cosh\left(2\left(T-t\right)\right)}e^{-\left(y^{2}_{1}+y^{2}_{2}\right)\tanh\left(2\left(T-t\right)\right)},\qquad h^{\mathrm{end}}_{t}\left(y\right)=e^{-2\left(T-t\right)}\left(\cos\left(2y_{1}\right)+\cos\left(2y_{2}\right)\right).\label{eq:shape-continuations}
\end{equation}
The area continuation varies along $y_{3}$ on the scale $1/4$ and
along $(y_{1},y_{2})$ on the scale $(2\tanh(2(T-t)))^{-1/2}$, which
is more than $8^{1/2}$ times larger. The endpoint continuation does
not depend on $y_{3}$. The discrete continuations $Q^{\mathrm{NV}}_{(i,n]}f$ follow from a
backward recursion on the lattice, so \eqref{eq:discrete-residual}
gives the recombination error $e_{n}(f)=\mathbb{P}^{R,-}_{n}(f)-\mathbb{P}^{\mathrm{WT}}_{n}(f)$
of each tree and its date contributions $c_{i}(f)=\langle\mathbb{D}^{R}_{i},(Q^{\mathrm{NV}}_{(i,n]}f)\circ\widehat{Y}_{i}\rangle$
exactly, separately from the common propagation error. The other settings
are those of Section~\ref{sec:models}, with $G=5$, so that $M_{G}=56$.
Both tasks are evaluated on the same tree of each method.

In this system $\widehat{V}_{i}>\widehat{D}_{i}$ at every date, so
the mixing coefficient \eqref{eq:mixing-coefficient} vanishes. Since
$b=0$, the model reference \eqref{eq:model-proxy} is the parent mean
of $A$, close to $\operatorname{diag}(1,1,t_{i-1}/2)$, and $\Gamma^{\mathrm{AWLL}}_{i}=\Gamma^{\mathrm{mod}}_{i}+\varepsilon_{i}I_{3}$
is close to ${\operatorname{diag}(1+\varepsilon_{i},1+\varepsilon_{i},t_{i-1}/2+\varepsilon_{i})}$,
with $\varepsilon_{i}$ from \eqref{eq:closed-coefficients} decreasing
from 0.62 to 0.28. Its ellipsoid is narrower than the WLL ball along
$y_{3}$, by factors 0.79 to 0.88, and wider along $(y_{1},y_{2})$,
by factors 1.13 to 1.27. Since the endpoint continuation involves only
$(y_{1},y_{2})$, its AWLL envelope \eqref{eq:derivative-envelope} is
close to $(1+\varepsilon_{i})^{(G+1)/2}\in[2.1,4.2]$ times its WLL envelope.
The derivatives of order $G+1$ of the area continuation are dominated
by the $y_{3}$ direction, and a direct computation from \eqref{eq:shape-continuations}
puts its AWLL envelope between 0.32 and 0.48 times its WLL envelope.
A common scalar instead multiplies the envelopes of both tasks by the
same factor, and the two tasks share the patches of each tree. Smaller
area contributions together with larger endpoint contributions on the
same pair of trees, at nearly equal populations, therefore reflect
the shape and not a common scalar.

At $\lambda=1$ the populations of the two trees differ by at most
6\% for $n$ from 16 to 128. At $n=32$ the endpoint contributions
vanish to rounding up to the eighth date. From the ninth date on,
the AWLL contribution is the larger one at every date, and the absolute
values sum to $7.3\times10^{-12}$ for AWLL and $9.7\times10^{-13}$
for WLL. The AWLL sum is also the larger one at each of eight values
of $\lambda$ from $1/2$ to $2$ at $n=16$ and $n=32$, where the two
populations differ by at most 7\%. The area contributions are
ordered the other way. From the fourth date on at $n=32$, the AWLL
contribution is the smaller one at 25 of the 28 dates, and the absolute
values sum to $3.2\times10^{-9}$ for AWLL and $1.2\times10^{-8}$ for
WLL. At $n=128$ the area errors are $-3.7\times10^{-9}$ for AWLL and
$6.9\times10^{-9}$ for WLL.

At $n=32$ both trees also have a large area contribution at the third
date, $-8.9\times10^{-9}$ for WLL and $-1.07\times10^{-8}$ for AWLL,
which dominates their area errors, $-9.1\times10^{-9}$ and $-1.36\times10^{-8}$.
A possible cause is the population rule of Section~\ref{sec:directional}.
At that date the frontier is still small, and since $(Y_{1},Y_{2})$
lies on the lattice, a patch at which the subdivision stops by this
rule can carry more distinct states than independent moment columns.
Such a patch is reduced although it fails \eqref{eq:partition-rule},
contrary to the independence assumed before \eqref{eq:joint-wms-budget},
and its contribution is not controlled by the test. We leave the analysis
of this contribution to future work.

All recombination errors are far below the propagation errors, $1.7\times10^{-3}$
and $1.5\times10^{-4}$ for the two tasks at $n=32$. Reducing the original
row tolerances $\tau^{\mathrm{org}}_{\alpha}$ of Appendix~\ref{app:algorithm}
from $10^{-8}$ to $10^{-10}$ changes the AWLL errors at $n=32$ by at
most 12\%, while the smaller WLL endpoint error is at the level
of that tolerance. Apart from the third date at $n=32$, AWLL thus lowers
the contributions of the task that varies along the narrow axis of
its ellipsoid and raises those of the task that varies along its wide
axes, as the envelopes above predict.

\FloatBarrier

\section{Discussion}

\label{sec:discussion}

\subsection{The resolution scale selected by AWLL}
\label{sec:selected-scale}
\noindent
The scalar $\lambda$ of the WLL test \eqref{eq:raw-rule} is dimensional.
The left side of the test carries the state unit, while the time and
mass factors on the right side are the same at every level, so a
change of state units can be compensated only through $\lambda$.
By Section~\ref{sec:directional}, a fixed $\lambda=1$ divides the
allowance in $y$ along each axis by the factor that the computational
coordinates apply to that axis. The test is therefore matched to the increments
only when their diffusion amplitudes are of order one in those coordinates.
At level 100 the leading amplitudes at the root are $S_{0}=100$ for
the Heston price tangent, by \eqref{eq:ver1-156}, and 60 for the sSABR
price. The test with $\lambda=1$ is then far stricter than the increments
require, and the frontier fragments into patches with little moment
redundancy. At nine intervals the final Heston WLL patches retain
147 parents on average against the capacity $M_{G}=462$ (Section~\ref{sec:support-growth}),
and 13.47\% of its patches are carried forward unchanged (Table~\ref{tab:identity-forwarding}).

In \eqref{eq:partition-rule} the metric $M_{i}=\Gamma^{-1/2}_{i}$
carries the inverse of a diffusion amplitude, since $\Gamma_{i}$ is
a second moment of state increments per unit time. Under a uniform
unit map $D\left(c\right)=cI_{N}$ and with no zero mixed block,
the AWLL matrix acquires the factor $c^{2}$ and the AWLL acceptance decisions on a given frontier
are unchanged (Appendix~\ref{sec:transport}). Write $\widetilde{\Gamma}_{i}=\Gamma^{\mathrm{AWLL}}_{i}/\lambda_{\max}\left(\Gamma^{\mathrm{AWLL}}_{i}\right)$,
so that
\begin{equation}
M^{\mathrm{AWLL}}_{i}=\lambda_{\max}\left(\Gamma^{\mathrm{AWLL}}_{i}\right)^{-1/2}\widetilde{\Gamma}^{-1/2}_{i}.\label{eq:metric-shape-scale}
\end{equation}
By Lemma~\ref{lem:metric-scalar}~(1), the AWLL test at $\lambda=1$
is the matrix test with $\Gamma_{i}=\widetilde{\Gamma}_{i}$ and
$\lambda=\lambda_{\max}(\Gamma^{\mathrm{AWLL}}_{i})^{1/2}$. With the
single block $L=1$ of Section~\ref{sec:numerics},
\eqref{eq:ridge-condition}
gives $\varphi^{-2}I_{6}\preceq\widetilde{\Gamma}_{i}\preceq I_{6}$.
Lemma~\ref{lem:metric-scalar}~(2) then
places the AWLL test between the WLL tests with scalars $\lambda_{\min}(\Gamma^{\mathrm{AWLL}}_{i})^{1/2}$
and $\lambda_{\max}(\Gamma^{\mathrm{AWLL}}_{i})^{1/2}$, whose ratio
is at most $\varphi$, and \eqref{eq:scale-transmission} places this
interval inside $[(\gamma^{+}_{i}/\varphi)^{1/2},(\varphi\gamma^{+}_{i})^{1/2}]$,
which \eqref{eq:forward-work-scale} uses for minimal partitions. In terms of patch counts, Corollary~\ref{cor:scale-count}
places the minimal AWLL count above the atomic level between constant
multiples of $\varphi^{-\psi(d)/2}X_{i}$ and $\varphi^{\psi(d)/2}X_{i}$, where
$X_{i}$ is, by Proposition~\ref{prop:minimal-partition}, the order of the
minimal WLL count at the scalar $(\gamma^{+}_{i})^{1/2}$. AWLL therefore acts
as a WLL test whose scalar
is the measured amplitude $(\gamma^{+}_{i})^{1/2}$ up to a factor
$\varphi^{\pm1/2}$, combined with a directional adjustment whose extreme
stretch ratio is at most $\varphi$.

Theorem~\ref{thm:finite-scale-work} states the unit response behind
this mechanism for a fixed frontier. Expanding a set of coordinates
by $c\ge1$ forces at least of order $c^{\psi(d)}$ WLL patches
up to the atomic saturation level, by \eqref{eq:concentration-work-lower},
whereas under the hypotheses of Lemma~\ref{lem:unit-response} the
AWLL load $\chi_{i}(\Gamma^{\mathrm{AWLL}}_{i}(c))$ stays bounded uniformly
in $c$. In the Zador regime, Remark~\ref{rem:discrete-zador}~(3) gives
the corresponding order $\left|\det D(c)\right|^{\psi(N)/N}$ for
the WLL test. The experiments of Section~\ref{sec:numerics} compare
complete multistep runs, whose frontiers differ between levels, so
they do not test these fixed-frontier statements directly. The
subdivision of a stored frontier can be repeated
under sweeps of the scalar and of the unit factor, the discrete Zador
functional \eqref{eq:discrete-zador} can be computed from the same
frontier, and the growth exponents of Corollary~\ref{cor:refinement-exponent}
can be compared on longer grids. These numerical checks are left to
future work.

The measured amplitude follows the current increments. The mixing
coefficient $\vartheta^{\mathrm{AWLL}}_{i}$ is close to one from $i=6$
on at every level (Figure~\ref{fig:wms-behavior}(g), (h)).
By \eqref{eq:bounded-model-injection} and the Lipschitz bound after
Theorem~\ref{thm:wms-derivation}, $\Gamma^{\mathrm{AWLL}}_{i}$
is then close to the matrix that \eqref{eq:wms-covariance} gives with
$\vartheta^{\mathrm{AWLL}}_{i}=1$, and $\gamma^{+}_{i}$ is close to
the largest eigenvalue of the edge reference.
The model reference supplies the scale at the first dates, where the
edge reference rests on few branches. At level 100, the largest eigenvalue
of the edge reference falls for Heston, as the price tangent amplitude
decays by \eqref{eq:heston-tangent-damping}, and rises slightly for
sSABR (Section~\ref{sec:global-lambda}), and $\lambda_{\max}(\Gamma^{\mathrm{AWLL}}_{8})^{1/2}$
equals 25.7 and 80.1 (Table~\ref{tab:mechanism}(b)).

These values agree with the post hoc WLL equivalents 26 and 80 at
$n=9$ in Table~\ref{tab:global-lambda}. They lie at the upper end
of the interval, where the two tests coincide along the leading eigenvector
of $\Gamma^{\mathrm{AWLL}}_{8}$. The interval concerns one frontier,
whereas the fitted scalar summarizes a complete history, so the agreement
is a consistency check and not a consequence of \eqref{eq:forward-work-scale}.
The direct WLL runs with $\lambda=26$ and $\lambda=80$ in Tables~\ref{tab:scalar-validation-work}
and \ref{tab:scalar-validation-errors} confirm the scale.
At level 100 the reduction by factors 137 and 1488 relative to WLL
with $\lambda=1$ is therefore obtained by one isotropic scalar of
the size of the measured amplitude. What AWLL contributes is the value
of that scalar, which is not known in advance. The fitted equivalents
vary from 14 to 26 with $n$ for Heston at level 100, and from 0.17
to 80 across the two models and three levels at $n=9$. AWLL obtains its scale
during the recombination, without a sweep over $\lambda$.

The same mechanism governs the smaller levels, for which
the unit maps \eqref{eq:ver1-149} and \eqref{eq:ver1-155} are not
uniform. At Heston level one, the late spectrum is that of level 100
scaled by about $c^{2}=10^{-4}$, with $\lambda_{\max}(\Gamma^{\mathrm{AWLL}}_{8})=6.9\times10^{-2}$
against 659. At level 0.01 it is set by the unscaled variance and
its $\nu_{0}$ tangent, whose diffusion energies at the root are $\xi^{2}\nu^{2}_{0}=0.0064$
and $\xi^{2}=0.16$ by \eqref{eq:heston-ito}, and $\lambda_{\max}(\Gamma^{\mathrm{AWLL}}_{8})=1.11\times10^{-2}$.
At both levels $\Gamma^{\mathrm{AWLL}}_{8}\prec I_{6}$, so AWLL is stricter
than WLL with $\lambda=1$ by Lemma~\ref{lem:metric-scalar}~(2), and
its terminal populations are 53.2 and 53.7 times larger at levels
one and 0.01. In sSABR
the $\omega_{0}$ price tangent is unscaled by \eqref{eq:ver1-155}
and has diffusion amplitude one at the root, and at level 0.01 the
volatility tangent has amplitude 1.5 (Appendix~\ref{sec:augmented-amplitude}).
The late sSABR spectra at the two smaller levels are of order one,
with $\lambda_{\max}(\Gamma^{\mathrm{AWLL}}_{8})=4.04$ and 2.46.

A fixed $\lambda=1$ is not a neutral reference at the smaller levels
either. At Heston level 0.01, WLL with $\lambda=1$ does not show the
error decrease of the other cases before $n=9$ in Figures~\ref{fig:european-errors}
and \ref{fig:asian-errors}. The post hoc equivalents lie inside the
interval $[\lambda_{\min}(\Gamma^{\mathrm{AWLL}}_{8})^{1/2},\lambda_{\max}(\Gamma^{\mathrm{AWLL}}_{8})^{1/2}]$
at Heston level one and at sSABR levels 0.01 and one, and above it
at Heston level 0.01, where patches are few (Table~\ref{tab:identity-forwarding})
and the population rule governs the subdivision. Away from level 100
they are not tested by direct runs, and AWLL at Heston level 0.01 also
has large errors on some of the last grids. The forward amplitude
fixes the order of magnitude of the scalar. The accuracy of the result
remains tied to the continuation derivatives through \eqref{eq:joint-error}.

\subsection{Directional shape}
\label{sec:shape-scale}
\noindent
The shape factor $\widetilde{\Gamma}^{-1/2}_{i}$ in \eqref{eq:metric-shape-scale}
is the part of AWLL that a single scalar cannot reproduce. The certificate volume of Proposition~\ref{prop:task-aligned}
depends only on the shape of $\Gamma$.
By Remark~\ref{rem:shape-ceiling}, for every task family the shape factor
changes the certificate volume of the identity by a factor between
$\varphi^{-(N-1)}$ and $\varphi^{N-1}\approx11.1$ for $N=6$.
Since \eqref{eq:retention-history} multiplies the retention fractions
of the interior dates, these per-date bounds do not by themselves
exclude a large multistep effect. At level 100, the direct WLL runs at the post hoc scalars $\lambda=26$ and $\lambda=80$
reach errors comparable to AWLL with smaller terminal populations,
so the directional adjustment is not the source of the reduction
relative to WLL with $\lambda=1$. Since these runs use one scalar at all dates,
following the date dependence of $\gamma^{+}_{i}$ is not needed for
the reported populations and accuracy.
The visible effect of the shape factor is
the more nearly circular projections
in Figure~\ref{fig:joint-states}.
In the synthetic system of Section~\ref{sec:shape-test}, where the
continuations are explicit, the shape lowers the area contributions
at most dates and raises the endpoint contributions on the same trees, as
\eqref{eq:shape-ceiling} suggests,
since the area task varies mainly along
the narrowest axis and the endpoint task only along the widest axes.

The potential gain of shape alignment depends on the derivatives
$D^{G+1}h^{\alpha}_{i}$, which are not available for the tasks of
Sections~\ref{sec:support-growth} to \ref{sec:covariance-evidence}. When the continuation derivatives are strongly anisotropic,
the condition bound \eqref{eq:ridge-condition} leaves most of that gain unused. Raising it within the
AWLL construction requires either relaxing the ridge, which weakens
the lower bound in \eqref{eq:scale-transmission} that limits the subdivision
caused by a nearly singular mixed block, or using several blocks, whose relative scales are then set by forward
second moments rather than by the continuation derivatives.

A one-sided safeguard keeps the shape and never makes the matrix test
stricter than the WLL test at a prescribed scalar $\lambda$. It multiplies
$\lambda$ by $\max\{1,\lambda_{\min}(\Gamma^{\mathrm{AWLL}}_{i})^{-1/2}\}$,
so that, by Lemma~\ref{lem:metric-scalar}, the matrix test accepts
every patch that the WLL test accepts at $\lambda$. This safeguard
was not tested here, and at Heston level 0.01 the WLL test whose
acceptances it preserves, with $\lambda=1$, does not show the error
decrease of the other cases (Section~\ref{sec:selected-scale}).

\FloatBarrier

\section{Conclusion}
\label{sec:conclusion}

\noindent
We have introduced AWLL, a matrix-valued extension of the
weighted Litterer--Lyons test that selects its resolution scale
from the increments of the current tree and the model coefficients,
without requiring an
advance choice of the scalar parameter $\lambda$. With one
coordinate block, the selected scale is comparable to the
measured increment amplitude up to a factor determined by the
golden ratio $\varphi$. Under the assumptions of
Theorem~\ref{thm:finite-scale-work}, the AWLL resolution load
on a fixed frontier remains bounded when a set of coordinate
units is expanded, whereas the lower bound on the number of
WLL patches grows as a power of the expansion factor until
atomic saturation. In the Heston and sSABR experiments at
level 100, AWLL reduces the terminal supports relative to WLL
with $\lambda=1$ by factors of 137 and 1488, respectively,
at comparable reference errors. WLL achieves similar results
with scalars selected post hoc from the AWLL runs. This comparison identifies
the automatic selection of an appropriate scale, rather than
the directional shape, as the main source of the observed
reductions. At the smaller Heston levels, AWLL instead retains
larger supports than WLL with $\lambda=1$.

The matrix formulation also provides joint error bounds for
prices and sensitivities, together with bounds on certificate
volumes, patch counts and retained populations. For
one-block AWLL, the effect of the directional shape on the certificate
volume is established with attainable two-sided bounds, but the Heston and sSABR experiments show no
clear improvement attributable to it. The synthetic example
illustrates how this effect depends on the continuation
derivatives of the tasks. The recombination order estimates
remain conditional on the growth of patch counts and
continuation derivatives, and the fixed-frontier unit bounds
do not directly establish the complexity of a multistep run.

For positive reduction, we have developed a batched basis
exchange algorithm that produces an extreme point of the
moment-feasibility set. For fixed moment rank and panel width,
its elimination stage requires $O(J)$ arithmetic operations
for a patch of $J$ nodes, excluding rank assessment and
initialization, compared with $\Theta(J^3)$ for the direct
kernel update. Other linear-in-$J$ reduction methods exist,
and a numerical comparison remains to be made. The numerical
rank overestimation observed for sSABR also indicates the
importance of detecting exact relations among moment
constraints.

Further work includes selecting the resolution scale under
a joint error budget, incorporating continuation information
to obtain task-adapted directional shapes, and identifying
problems in which the directional shape yields measurable improvements.
Experiments on finer grids and broader volatility regimes
would help assess the observed growth patterns and the work
required to meet prescribed tolerances. For larger moment
ranks, improved basis updates and the removal of exact linear
invariants are additional directions for development.

\section*{Disclosure statement}
\noindent
The authors declare no competing financial or nonfinancial interests.

\section*{Funding}
\noindent
No funding was received for conducting this study.

\section*{Data and code availability}
\noindent
All experimental data were generated from the mathematical models
specified in this article. The data and code supporting the findings
are available from the corresponding author on reasonable request.
The code will be made available on GitHub after publication.

\section*{References}

\bibliographystyle{tfcad}
\begingroup\color{black}
\hypersetup{citecolor=black,urlcolor=black,linkcolor=black}
\bibliography{Bibliography_ver2}
\endgroup

\clearpage\appendix\pagestyle{wmsappendix}\thispagestyle{wmsappendix}

\FloatBarrier

\section{Unit response, amplitudes and metric action}

\label{app:geometry}

\subsection{Unit response of the AWLL matrix}
\label{sec:transport}
\noindent
Fix $i\in\mathcal{R}$, the parent measure $\mathbb{P}^{R,+}_{i-1}$ and the frontier
$\mathbb{P}^{R,-}_{i}$. The unit map
\begin{equation}
D\left(c\right)=\operatorname{diag}_{j=1,\ldots,N}\left(c^{\zeta_{j}}\right),\qquad c>0,\label{eq:ver1-096}
\end{equation}
transports them by mapping every state $y$ to $y'=D(c)y$ and keeping
the weights, and it transforms the coefficients
of \eqref{eq:augmented} into $b'(y')=D(c)b(y)$ and $\Sigma'(y')=D(c)\Sigma(y)$.
The numbers $\zeta_{j}$ are called the unit exponents. If a coordinate
$X_{j}$ of the base state of Remark~\ref{rem:greeks} and a parameter
$\theta_{h}$ change as $X_{j}'=c^{\zeta^{X}_{j}}X_{j}$ and $\theta_{h}'=c^{\zeta^{\theta}_{h}}\theta_{h}$,
then
\begin{equation}
\partial_{\theta_{h}'}=c^{-\zeta^{\theta}_{h}}\partial_{\theta_{h}},\qquad\partial_{\theta_{h}'}X_{j}'=c^{\zeta^{X}_{j}-\zeta^{\theta}_{h}}\partial_{\theta_{h}}X_{j},\label{eq:jet-units}
\end{equation}
so the tangent coordinate $\partial_{\theta_{h}}X_{j}$ has the exponent
$\zeta^{X}_{j}-\zeta^{\theta}_{h}$. A running integral has the exponent
of its integrand, since the time unit is fixed.

A quantity computed after the transport is written with the argument
$c$, and a tuple of blocks is identified with the block-diagonal matrix
it defines. Then $\widehat{\Gamma}^{\mathrm{edge}}_{i}(c)=D(c)\widehat{\Gamma}^{\mathrm{edge}}_{i}(1)D(c)$,
the same congruence holds for $\Psi_{e}$ and for the model reference
\eqref{eq:model-proxy}, and the reference discrepancy and the edge
dispersion \eqref{eq:sensitivity-objective} become
\begin{equation}
\begin{aligned}\widehat{D}_{i}\left(c\right) & =\sum^{N}_{j,k=1}c^{2\left(\zeta_{j}+\zeta_{k}\right)}\left(\widehat{\Gamma}^{\mathrm{edge}}_{i}\left(1\right)-\Gamma^{\mathrm{mod}}_{i}\left(1\right)\right)^{2}_{jk},\\
\widehat{V}_{i}\left(c\right) & =\sum^{N}_{j,k=1}c^{2\left(\zeta_{j}+\zeta_{k}\right)}\sum_{e}\varpi^{2}_{e}\left(\Psi_{e}\left(1\right)-\widehat{\Gamma}^{\mathrm{edge}}_{i}\left(1\right)\right)^{2}_{jk}.
\end{aligned}
\label{eq:score-unit-response}
\end{equation}
 For the uniform unit map $D(c)=cI_{N}$, \eqref{eq:score-unit-response}
leaves the mixing coefficient unchanged, and if no mixed block is
zero, \eqref{eq:closed-coefficients} gives $\Gamma^{\mathrm{AWLL}}_{i}(c)=c^{2}\Gamma^{\mathrm{AWLL}}_{i}(1)$,
so that the matrix test \eqref{eq:partition-rule} accepts the same
patches before and after the transport. For other exponents, the
mixing coefficient can vary with $c$ through \eqref{eq:score-unit-response}.

\begin{lem}
\label{lem:unit-response}Let $L=1$, let $\zeta_{j}\le1$ for
all $j$, let $\Pi_{s}\ne0$ be the projection of Theorem~\ref{thm:finite-scale-work},
and assume that
\begin{equation*}
\underline{\gamma}_{i}=\frac{1}{\operatorname{rank}\Pi_{s}}\min\left\{ \operatorname{tr}\left(\Pi_{s}\widehat{\Gamma}^{\mathrm{edge}}_{i}\left(1\right)\right),\operatorname{tr}\left(\Pi_{s}\Gamma^{\mathrm{mod}}_{i}\left(1\right)\right)\right\} >0.
\end{equation*}
Let $c\ge1$. Then, whatever the value of $\vartheta^{\mathrm{AWLL}}_{i}(c)$
in $\left[0,1\right]$, the largest eigenvalue $\gamma^{+}_{i}(c)$ of the
mixed block satisfies
\begin{equation}
c^{2}\underline{\gamma}_{i}\le\gamma^{+}_{i}\left(c\right)\le c^{2}\max\left\{ \lambda_{\max}\left(\widehat{\Gamma}^{\mathrm{edge}}_{i}\left(1\right)\right),\lambda_{\max}\left(\Gamma^{\mathrm{mod}}_{i}\left(1\right)\right)\right\} ,\label{eq:dominant-unit-scale}
\end{equation}
and, with the golden ratio $\varphi$, the AWLL matrix
\eqref{eq:wms-covariance} satisfies
\begin{equation}
\left\lVert \Gamma^{\mathrm{AWLL}}_{i}\left(c\right)^{-1/2}D\left(c\right)\right\rVert _{2}\le\left(\varphi/\underline{\gamma}_{i}\right)^{1/2}.\label{eq:normalized-unit-bound}
\end{equation}
\end{lem}
\begin{proof}
Write $\vartheta=\vartheta^{\mathrm{AWLL}}_{i}(c)$ and $H=\vartheta\widehat{\Gamma}^{\mathrm{edge}}_{i}(1)+(1-\vartheta)\Gamma^{\mathrm{mod}}_{i}(1)$,
so that the mixed block is $D(c)HD(c)$. Since $D(c)\Pi_{s}=c\Pi_{s}$, the mean of the Rayleigh quotients
of $D(c)HD(c)$ over an orthonormal basis of the range of $\Pi_{s}$
is $c^{2}\operatorname{tr}(\Pi_{s}H)/\operatorname{rank}\Pi_{s}$.
The largest eigenvalue is at least this mean, and $\operatorname{tr}(\Pi_{s}H)$,
a convex combination of the two traces, is at least the smaller one,
which proves the lower bound in \eqref{eq:dominant-unit-scale}.
Since $c\ge1$, $\zeta_{j}\le1$ and $\Pi_{s}\ne0$, $\lVert D(c)\rVert_{2}=c$.
Hence $D(c)\widehat{\Gamma}^{\mathrm{edge}}_{i}(1)D(c)\preceq c^{2}\lambda_{\max}(\widehat{\Gamma}^{\mathrm{edge}}_{i}(1))I_{N}$,
likewise for the model reference, and the convexity of the largest
eigenvalue gives the upper bound. By the lower bound, the mixed block is nonzero, and
\eqref{eq:scale-transmission} gives $\Gamma^{\mathrm{AWLL}}_{i}(c)\succeq\varphi^{-1}c^{2}\underline{\gamma}_{i}I_{N}$.
Since $\lVert\Gamma^{\mathrm{AWLL}}_{i}(c)^{-1/2}D(c)\rVert_{2}\le c\lVert\Gamma^{\mathrm{AWLL}}_{i}(c)^{-1/2}\rVert_{2}$,
this proves \eqref{eq:normalized-unit-bound}.
\end{proof}

 By \eqref{eq:dominant-unit-scale} and \eqref{eq:scale-transmission},
every eigenvalue of $\Gamma^{\mathrm{AWLL}}_{i}(c)$ is of order $c^{2}$.
A displacement along a coordinate with $\zeta_{j}=0$ therefore shrinks
by a factor of order $1/c$ in the normalized coordinates, and the effect
of this coarser resolution on the tasks enters \eqref{eq:joint-error}
through the envelope \eqref{eq:derivative-envelope}.

\subsection{Amplitudes of the two augmented models}\label{sec:augmented-amplitude}
\noindent
At the Heston root, $Y_{4,0}=2\nu_{0}$, $Y_{2,0}=\nu^{2}_{0}$, and
$Y_{3,0}=0$. The first diffusion coefficient of the price tangent
is consequently 
\begin{equation}
\sqrt{Y_{2,0}}Y_{3,0}+\frac{Y_{1,0}Y_{4,0}}{2\sqrt{Y_{2,0}}}=S_{0}.\label{eq:ver1-156}
\end{equation}
 Although the tangent starts at zero, the model reference \eqref{eq:model-proxy}
at $i=1$ therefore has the $(3,3)$ entry $S^{2}_{0}$, against the
$(1,1)$ entry $S^{2}_{0}\nu^{2}_{0}$ of the price and the $(5,5)$
entry $s_{1}S^{2}_{0}$ that the drift of its running integral $Y_{5}$
gives.
At level 100 this is the largest diagonal entry, and by the lower
bound in \eqref{eq:scale-transmission} every eigenvalue of $\Gamma^{\mathrm{AWLL}}_{1}$
is at least $\varphi^{-1}$ times the $(3,3)$ entry of the mixed
block, so this amplitude sets the scale in every direction.

Along the paths this contribution decays. Let $B^{\mathrm{v}}=\rho B^{1}+\sqrt{1-\rho^{2}}B^{2}$
be the Brownian motion that drives the variance in \eqref{eq:heston-ito}.
Since $Y_{2}$ stays positive, It\^o's formula and differentiation in
$\nu_{0}=\sqrt{Y_{2,0}}$ give
\begin{equation}
\begin{aligned}d\sqrt{Y_{2,t}} & =\left(\frac{4\kappa\theta-\xi^{2}}{8\sqrt{Y_{2,t}}}-\frac{\kappa}{2}\sqrt{Y_{2,t}}\right)dt+\frac{\xi}{2}dB^{\mathrm{v}}_{t},\\
\partial_{\nu_{0}}\sqrt{Y_{2,t}} & =\exp\left(-\frac{\kappa t}{2}-\frac{4\kappa\theta-\xi^{2}}{8}\int^{t}_{0}\frac{du}{Y_{2,u}}\right)\le e^{-\kappa t/2},
\end{aligned}
\label{eq:heston-tangent-damping}
\end{equation}
where the inequality uses $4\kappa\theta=0.48>\xi^{2}=0.16$. Since $Y_{4,t}/(2\sqrt{Y_{2,t}})=\partial_{\nu_{0}}\sqrt{Y_{2,t}}$,
the term $y_{1}y_{4}/(2\sqrt{y_{2}})$ of the price tangent coefficient
$\Sigma^{\mathrm{Heston}}_{31}=\sqrt{y_{2}}y_{3}+y_{1}y_{4}/(2\sqrt{y_{2}})$
in \eqref{eq:heston-ito} is at most $Y_{1,t}e^{-\kappa t/2}$ along
the path, while the other term $\sqrt{y_{2}}y_{3}$ starts at zero
with $Y_{3}$.
Consistently, at level 100 the $(3,3)$ entry and the largest eigenvalue
of the edge reference decrease between $i=1$ and $i=8$ (Section~\ref{sec:covariance-evidence}).

For sSABR, the volatility tangent starts at $Y_{4,0}=\left(F_{0}+s\right)^{-\beta}$,
so its diffusion amplitude at the root is $\nu\left(F_{0}+s\right)^{-\beta}$,
which is 1.5, 0.15 and 0.015 at the levels 0.01, 1 and 100. The diffusion
amplitude $\omega_{0}=\alpha_{0}\left(F_{0}+s\right)^{\beta}$ of
the price is 0.006, 0.6 and 60, and that of the price tangent is
one at every level. At the root, the volatility tangent therefore has the
largest diffusion amplitude at level 0.01, and the price at level
100.

\subsection{Metric action on the reduced measures}\label{app:metric-action}
\begin{figure}[!htbp]
\centering
\includegraphics[width=0.98\columnwidth]{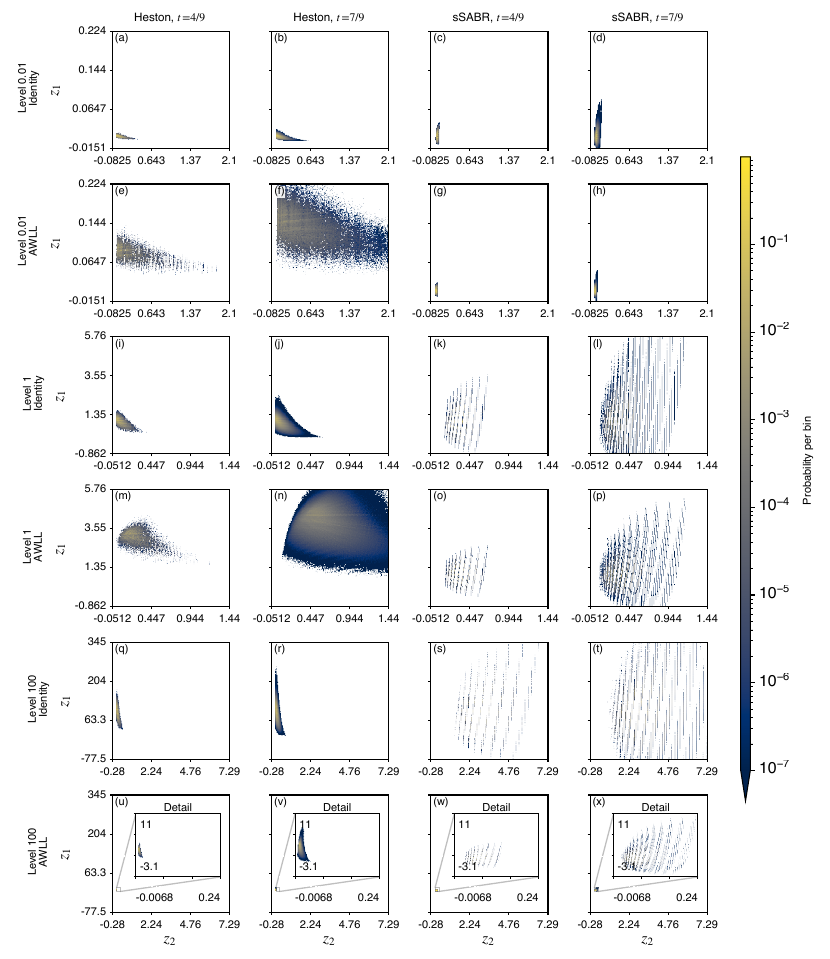}

\caption{Metric action at $n=9$. In each pair of rows, the row labeled
Identity applies $M=I_{6}$ and the row labeled AWLL applies the recorded
$M^{\mathrm{AWLL}}_{i}$, both to the same reduced states and probabilities
of AWLL. The displayed $z_{j}=\left(M\widehat{Y}_{i}\right)_{j}$, $j=1,2$,
follow the full six-dimensional matrix action. Row pairs give levels 0.01, 1 and
100. Columns give Heston and sSABR at $t=4/9,7/9$. The eight main
panels at each level share axis ranges and bins, as do the four level
100 detail windows. Color gives probability per bin. \label{fig:joint-states}}
\end{figure}

\noindent
Figure~\ref{fig:joint-states} shows how the identity and AWLL metrics
act on the reduced measures of the AWLL runs at the three levels.

At level one and $t=7/9$, the ratio of the major to the minor axis
of the probability-weighted covariance ellipse of $(z_{1},z_{2})$,
each divided by its common axis span in Figure~\ref{fig:joint-states},
falls from 1.98 under the identity to 1.47 under AWLL for Heston,
and from 1.59 to 1.14 for sSABR. Under the identity, $z_{1}$ and $z_{2}$ are
the two coordinates of the base state (Remark~\ref{rem:greeks}),
whereas under AWLL they also depend on the four augmented coordinates.

The ratio $\left(\widehat{Y}_{i}\right)_{2}/\left(Y_{0}\right)_{2}$
of the second coordinate, the variance in Heston and the volatility
amplitude in sSABR, to its initial value is not changed by the unit
maps \eqref{eq:ver1-149} and \eqref{eq:ver1-155} and does not involve
the metric. With the node probabilities of the reduced AWLL measure
of each date as weights, its standard deviation ranges from 0.775
to 0.796 for Heston at $t=4/9$ and from 0.802 to 0.815 at $t=7/9$,
and from 0.2020 to 0.2026 and from 0.2693 to 0.2697 for sSABR. These similar values across the three levels
contrast with the much larger changes of the distance ratio \eqref{eq:metric-distance-ratio}
and the terminal population ratio in Table~\ref{tab:mechanism}.

\clearpage

\FloatBarrier

\section{Batched basis exchange algorithm}

\label{app:algorithm}

\noindent
Algorithm~\ref{alg:bbe} implements the positive reduction from Section~\ref{sec:solver}. Its working dimension $k$ is the numerical estimate
$k_{\tau}$ from \eqref{eq:rank-scaling}. In the experiments, $\operatorname{rank}_{\tau}$
counts the singular values of $A^{\mathrm{eq}}$ above $10^{-18}\max\{M_{G},J\}\sigma_{\max}(A^{\mathrm{eq}})$,
with the consequence for sSABR discussed in Section~\ref{sec:global-lambda}.
The bound $b_{\max}$ limits
panel width. The tolerances $\tau_{\mathrm{pos}}\ge0$ and $\tau_{\mathrm{piv}}\ge0$
govern intermediate prefix feasibility and admissible pivot entries.
The upper bound $\kappa^{\mathrm{basis}}_{\max}$ limits $\kappa_{2}\left(U_{\cdot B}\right)$,
the ratio of the largest to smallest singular value. Tolerances for orthogonality
and reconstruction apply to $U\prescript{t}{}{U}-I_{k}$ and $A^{\mathrm{org}}-RU$.
Final acceptance requires nonnegative weights, the total mass within
its tolerance, and $\left|d^{\mathrm{org}}_{\alpha}\right|\le\tau^{\mathrm{org}}_{\alpha}$
for every original row $\alpha$. A failed test returns the original measure
of the patch.

\RestyleAlgo{ruled}\LinesNumbered\DontPrintSemicolon\SetAlgoVlined\SetAlFnt{\small}\SetAlgoNlRelativeSize{-1}

\SetAlFnt{\fontsize{8}{9}\selectfont}
\begin{algorithm*}[H]
\caption{Batched basis exchange on one patch.}\label{alg:bbe}
\KwIn{\(A^{\mathrm{org}}\), \(p^{R,-}\ge0\), panel width bound \(b_{\max}\), rank threshold \(\tau\), and validation tolerances.}
\KwOut{Validated \(p^{R,+}\), or the original measure.}
\(A^{\mathrm{eq}}\leftarrow\operatorname{SinkhornKnoppL2}\left(A^{\mathrm{org}}\right)\)\;
\(k\leftarrow\operatorname{rank}_{\tau}\left(A^{\mathrm{eq}}\right)\)\;
\(\left(R,U\right)\leftarrow\operatorname{MassAnchoredMGS}\left(A^{\mathrm{org}},k,J^{-1/2}\prescript{t}{}{\mathbf{1}_{J}}\right)\)\;
\If{\(k\notin\left\{1,\ldots,\min\left(M_G,J\right)\right\}\) or either the orthogonality check or the original row reconstruction check fails}{\Return{original measure}\;}
\(B\leftarrow\operatorname{RankRevealingBasis}\left(U,k,\kappa_{\max}^{\mathrm{basis}}\right)\)\;
\If{basis validation fails}{\Return{original measure}\;}
Factor \(U_{\cdot B}\) and set \(H_B\leftarrow U_{\cdot B}^{-1}\), \(\eta\leftarrow p_B^{R,-}\), \(\mathcal J_{\mathrm{rem}}\leftarrow B^c\)\;
\ForEach{\(j\in\mathcal J_{\mathrm{rem}}\)}{\(\tilde{p}_j\leftarrow p_j^{R,-}\)\;}
\While{\(\mathcal J_{\mathrm{rem}}\ne\varnothing\)}{
  take the next panel \(J_b\) of the stream, with \(m_b\le b_{\max}\)\;
  \(C_b\leftarrow H_BU_{\cdot J_b}\), \(\mathrm{changed}\leftarrow\mathrm{false}\)\;
  \While{\(J_b\ne\varnothing\)}{
    \(\eta^{\left(0\right)}\leftarrow\eta\), \(\eta^{\left(h\right)}\leftarrow\eta+\sum_{q=1}^{h}\tilde{p}_{j_q}c_{j_q}\)\;
    \(h_b^*\leftarrow\max\left\{0\le h\le m_b\middle|\eta^{\left(q\right)}\ge-\tau_{\mathrm{pos}}\mathbf1_k,\ 1\le q\le h\right\}\)\;
    absorb the prefix, set \(\eta\leftarrow\eta^{\left(h_b^*\right)}\), and zero its remaining weights\;
    remove the prefix indices from \(\mathcal J_{\mathrm{rem}}\), \(J_b\), and the corresponding columns of \(C_b\)\;
    \If{\(J_b\ne\varnothing\)}{
      let \(j\) and \(c_j\) be the first remaining index and coefficient vector\;
      \(\mathcal N_j\leftarrow\left\{a\middle|\left(c_j\right)_a<-\tau_{\mathrm{piv}}\right\}\)\;
      \If{\(\mathcal N_j=\varnothing\)}{\Return{original measure}\;}
      select \(a^{*}\in\operatorname*{arg\,min}_{a\in\mathcal N_j}\frac{-\eta_a}{\left(c_j\right)_a}\) with fixed tie order\;
      \(\delta_j\leftarrow-\eta_{a^{*}}/\left(c_j\right)_{a^{*}}\)\;
      \If{\(\delta_j<0\) or \(\delta_j>\tilde{p}_j\)}{\Return{original measure}\;}
      replace \(\beta_{a^{*}}\) by \(j\) in the ordered basis \(B\)\;
      \(\eta\leftarrow\eta+\delta_jc_j+\left(\tilde{p}_j-\delta_j\right)e_{a^{*}}\), \(\tilde{p}_j\leftarrow0\)\;
      remove \(j\) from \(\mathcal J_{\mathrm{rem}}\) and \(J_b\), and remove its coefficient vector\;
      \(C_b\leftarrow C_b-\frac{\left(c_j-e_{a^{*}}\right)\left(\prescript{t}{}{e_{a^{*}}} C_b\right)}{\left(c_j\right)_{a^{*}}}\)\;
      \(\mathrm{changed}\leftarrow\mathrm{true}\)\;
    }
    \(m_b\leftarrow\left|J_b\right|\)\;
  }
  \If{\(\mathrm{changed}\)}{
    \If{\(\kappa_2\left(U_{\cdot B}\right)>\kappa_{\max}^{\mathrm{basis}}\)}{\Return{original measure}\;}
    factor \(U_{\cdot B}\) and refresh \(H_B\leftarrow U_{\cdot B}^{-1}\)\;
  }
}
\(p_B^{R,+}\leftarrow\eta\), \(p_{B^c}^{R,+}\leftarrow0\)\;
\If{\(\min_j p_j^{R,+}<0\) or any validation check for the support, mass, or original moment rows fails}{\Return{original measure}\;}
\Return{\(p^{R,+}\)}\;
\end{algorithm*}

\clearpage
\fancyhead[C]{}

\FloatBarrier

\ifdefined\AMFsubmissionlayout
\pagestyle{fancy}
\processdelayedfloats
\let\processdelayedfloats\relax
\clearpage
\renewcommand{\listfigurename}{Figure captions}
{\hypersetup{linkcolor=black}\listoffigures}
\fi

\end{document}